\documentclass[review,3p]{IEEEtran}

\usepackage{mathrsfs,latexsym,float,enumerate,amssymb,amsmath,bm,color,lineno}
\usepackage{amscd,verbatim,cuted}
\usepackage{graphicx}
\usepackage{epstopdf}
\usepackage[colorlinks=true,linkcolor=blue,citecolor=blue]{hyperref}
\usepackage{booktabs}
\usepackage{subfigure}

\newtheorem{theorem}{Theorem}[section]
\newtheorem{lemma}{Lemma}[section]
\newtheorem{proposition}{Proposition}[section]

\newtheorem{definition}{Definition}[section]

\newtheorem{remark}{Remark}

\usepackage{pifont}

\begin{document}

\title{Picture Fuzzy Interactional Bonferroni Mean Operators via Strict Triangular Norms
and Applications to Multi-Criteria Decision Making}

\author{Lantian Liu, Xinxing Wu, Guanrong Chen,~\IEEEmembership{Life Fellow,~IEEE}
\thanks{Manuscript received xx 00, 201x; revised xx 00, 201x.}
\thanks{L. Liu is with the School of Sciences, Southwest Petroleum University,
Chengdu, Sichuan 610500, China. e-mail: (llt13207428373@163.com).}
\thanks{X. Wu is with (1) the School of Mathematics and Statistics, Guizhou University 
of Finance and Economics, Guiyang 550025, China; 
(2) the School of Sciences, Southwest Petroleum University, Chengdu, Sichuan 610500, China.
e-mail: (wuxinxing5201314@163.com).}
\thanks{G. Chen is with the Department of Electrical Engineering, City University of
Hong Kong, Hong Kong SAR, China.
e-mail: (eegchen@cityu.edu.hk).}
\thanks{All correspondences should be addressed to X. Wu.}
\thanks{This work was supported by the Natural Science Foundation of
Sichuan Province (No.~2022NSFSC1821) and the Young Scholars Science Foundation 
of Lanzhou Jiaotong University (No.~2020022).}
}


\maketitle

\begin{abstract}
Based on the closed operational laws in picture fuzzy numbers and strict triangular norms,
we extend the Bonferroni mean (BM) operator under the picture fuzzy environment to propose the
picture fuzzy interactional Bonferroni mean (PFIBM), picture fuzzy interactional weighted
Bonferroni mean (PFIWBM), and picture fuzzy interactional normalized weighted Bonferroni
mean (PFINWBM) operators. We prove the monotonicity, idempotency, boundedness, and
commutativity for the PFIBM and PFINWBM operators. We also establish a novel multi-criteria
decision making (MCDM) method under the picture fuzzy environment by applying the PFINWBM operator.
Furthermore, we apply our MCDM method to the enterprise resource planning (ERP) systems selection.
The comparative results for our MCDM method induced by six classes of well-known triangular norms
ensure that the best selection is always the same ERP system. Therefore, our MCDM method is
effective for dealing with the picture fuzzy MCDM problems.
\end{abstract}


\IEEEpeerreviewmaketitle

\section{Introduction}
\IEEEPARstart{S}INCE Zadeh's fuzzy sets~\cite{Za1965}, Atanassov intuitionistic fuzzy sets
(IFSs)~\cite{Ata1999}, and hesitant fuzzy sets (HFSs)~\cite{T2010} can not be used in the
evaluation of the situations that have more than three different aspects, Cuong
and Kreinovich~\cite{CK2014,CK2013} introduced the concept of picture fuzzy sets (PFSs) to
extend the evaluation dimension. The PFSs are characterized by a positive membership function,
a neutral membership function, and a negative membership function, whose sum at every point
is less than or equal to $1$. In analogy to IFSs, considering the uncertainty of three
membership degrees for PFSs, Cuong~\cite{CK2014} introduced the concept of interval-valued
PFSs (IvPFSs). As an extension of IFSs, when the neutral membership function of a PFS is
equal to $0$, the PFS reduces to an IFS. Combining PFSs with HFSs, Ullah et al.~\cite{UAJMM2018}
introduced the concept of picture hesitant fuzzy sets (PHFSs).

For some practical multi-criteria decision-making (MCDM) problems, the aggregation operator,
which aims to provide an effective single overall representation of the input values,
is a useful tool. Moreover, the interrelationships among multi-criteria in practical MCDM
problem are universal. To capture such interrelationship for each pair of criteria,
Bonferroni~\cite{Bo1950} proposed a mean-type aggregation operator for crisp numbers,
called the Bonferroni mean (BM). Then, it was successfully generalized and
applied to various decision-making environments, including fuzzy environment (\cite{LLC2012,ZLW2015}),
intuitionistic fuzzy environment (\cite{BJ2013,DGM2016,GA2018,HH2016,HHC2015,LCL2017,Tian2020,
XXZ2012,XXZ2013,XC2012,XC2011,XY2011,ZH2012}), q-Rung orthopair fuzzy environment (\cite{LL2018,LW2018-2,YP2019}),
(picture) hesitant fuzzy environment (\cite{Deli2021,HHWC2015,MAA2021,Xu2014,XZ2018,ZX2013,ZXX2012}),
type-2 fuzzy environment (\cite{Chiao2021,GHZLD2015,HKMCBPN2013,LG2021,LGM2019}), spherical fuzzy
environment (\cite{FSDGK2021}), and picture fuzzy environment (\cite{AA2020-1,Wei2017K}). For
example, Xu et al.~\cite{XXZ2012,XC2011,XY2011} extended the BM to the intuitionistic fuzzy
Bonferroni mean (IFBM), the weighted intuitionistic fuzzy Bonferroni mean (WIFBM), the interval-valued
intuitionistic fuzzy Bonferroni mean (IVIFBM), and the weighted interval-valued intuitionistic fuzzy
Bonferroni mean (WIVIFBM), for aggregating intuitionistic fuzzy and interval-valued intuitionistic
fuzzy information. They obtained the idempotency, monotonicity, commutativity, and boundedness for
IFBM, WIFBM, IVIFBM, and WIVIFBM. Then, Das et al.~\cite{DGM2016} extended IFBM and WIFBM operators
by using strict t-conorm and introduced the AIF-EBM operator and the weighted AIF-EBM (WAIFEBM)
operator, which were successfully applied to MCDM problems. Considering the interactions between
the membership and the nonmembership functions under the intuitionistic fuzzy environment,
He et al.~\cite{HHC2015} presented the intuitionistic fuzzy interaction BM (IFIBM) and the
weighted intuitionistic fuzzy interaction BM (WIFIBM) operators and extended them to EIFIBM
operator by using continuous Archimedean t-norm \cite{HH2016}. However, the three operators
does not have the monotonicity, which may lead to some unreasonable aggregation results,
{although their idempotency and commutativity} were proved by He et al.~\cite{HH2016,HHC2015}.
Liu and Gao~\cite{LG2021} developed a novel green supplier selection MCDM approach based on
the interval type-2 fuzzy prioritized Choquet normalized weighted BM operators. Chiao~\cite{Chiao2021}
presented some new models for MCDM problems by developed the quantifier guided ordered
weighted averaging (QGOWA) and BM operators under the interval type-2 fuzzy environment.
Based on the BM and the geometric mean, Zhu et al.~\cite{ZX2013,ZXX2012} proposed the
hesitant fuzzy geometric Bonferroni mean (HFGBM), hesitant fuzzy Choquet geometric Bonferroni
mean (HFCGBM), and hesitant fuzzy Bonferroni mean (HFBM) operators for hesitant fuzzy date, and
applied them to MCDM. Combining BM and the power mean in hesitant fuzzy environments,
He et al.~\cite{HHWC2015} obtained the hesitant fuzzy power geometric Bonferroni mean and
the hesitant fuzzy power Bonferroni mean, and proposed a new approach for hesitant fuzzy MAGDM.

To extend the evaluation dimension and {applied range} of IFSs, many scholars started
to investigate the picture fuzzy BM (PFBM), picture hesitant fuzzy BM (PHFBM),
and q-Rung orthopair fuzzy BM (q-ROFBM) operators and their applications to MCDM.
Liu and Wang~\cite{LW2018-2} introduced the q-ROF Archimedean BM (q-ROFABM) and the q-ROF weighted
Archimedean BM (q-ROFWABM) operators, and developed a novel MCDM by q-ROFWABM operator.
Wei~\cite{Wei2017K} proposed some picture uncertain linguistic BM (PULBM) operators and applied
them to choosing communications industry for their service outsourcing suppliers. Recently, Ate\c{s}
and Akay~\cite{AA2020-1} presented a new picture fuzzy MCDM method by introducing the PFBM,
the picture fuzzy normalized weighted BM (PFNWBM), and the picture fuzzy ordered weighted BM (PFOWBM) operators.
Mahmood et al.~\cite{MAA2021} provided a new method for MAGDM under the the picture hesitant fuzzy
environment based on the PHFBM, the picture hesitant fuzzy weighted BM (PHFWBM),
the picture hesitant fuzzy geometric BM (PHFGBM), and the picture
hesitant fuzzy weighted geometric BM (PHFWGBM) introduced by them. However, as Wu et al.
pointed out in ~\cite{WZCLZY2021}, the picture fuzzy aggregation operators in \cite{AA2020-1,Wei2017K}
have a common shortcoming: They may not satisfy the condition that the sum of the three degrees
cannot exceed $1$. Thus, such operators are not closed in PFNs. Meanwhile, it is easy to see
that the operations $\otimes$ and $A^{\lambda}$ in \cite{MAA2021} are not closed in PHFS,
because they are both defined by using two triangular conorms.

To obtain closed picture fuzzy BM operators, based on some basic operations of PFNs
introduced by Wu et al.~\cite{WZCLZY2021}, this paper introduces the picture fuzzy
interactional Bonferroni mean (PFIBM), the picture fuzzy interactional
weighted Bonferroni mean (PFIWBM), and the picture fuzzy interactional
normalized weighted Bonferroni mean (PFINWBM) operators for PFNs, which are proved to
be monotonous, idempotent, bounded, and commutative. Moreover, a novel MCDM method under
the picture fuzzy environment is proposed by using the PFINWBM operator.
Furthermore, this MCDM method is used for enterprise resource planning systems selection.
By using six classes of well-known triangular norms, including the algebraic product
$T_{\mathbf{P}}$, Schweizer-Sklar t-norm $T_{\gamma}^{\mathbf{SS}}$, Hamacher t-norm
$T_{\gamma}^{\mathbf{H}}$, Frank t-norm $T_{\gamma}^{\mathbf{F}}$, Dombi t-norm
$T_{\gamma}^{\mathbf{D}}$, and Acz\'{e}l-Alsina t-norm $T_{\gamma}^{\mathbf{AA}}$,
the best ERP system is always the same one, demonstrating the effectiveness of our method.

\section{Preliminaries}\label{Sec-2}

\subsection{Bonferroni mean (BM)}

Being an important aggregation operator for crisp date, Bonferroni
mean was originally introduced by Bonferroni~\cite{Bo1950} as follows.

\begin{definition}[{\textrm{\protect\cite{Bo1950}}}]
\label{Def 1}
Let $p, q \geq 0$ and $a_{i}$ ($i= 1, 2, \ldots, n$) be a collection
of crisp data with $a_{i} \geq 0$. The function $\mathrm{B}^{p,q}$ defined by
{\small $$
\mathrm{B}^{p,q} (a_{1}, \ldots ,a_{n})
=\left[\frac {1}{n(n\!-1)} \sum\limits_{i, j\!=1 \atop i\neq j}^{n}a_{i}^{p}\cdot
a_{j}^{q}\right]^{\frac {1}{p\!+q}},
$$}
is called the \textit{Bonferroni mean} (BM) operator.
\end{definition}

\subsection{Picture fuzzy sets}

\begin{definition}[{\textrm{\protect\cite[Definition 3.1]{CK2014}}}]
\label{PFS-Def}
Let $X$ be the universe of discourse. A \textit{picture fuzzy set} (PFS) $P$
on $X$ is defined as an object with the following form
\begin{equation}
\label{2}
P=\left\{\langle x; \mu_{P}(x), \eta_{P}(x), \nu_{P}(x)\rangle\mid x\in
X\right\},
\end{equation}
where $\mu_{P}(x)$, $\eta_{P}(x)$, $\nu_{P}(x)\in [0, 1]$,
$\mu_{P}(x)+\eta_{P}(x)+\nu_{P}(x)\leq 1$, and $\pi_{P}(x)=1-(\mu_{P}(x)+
\eta_{P}(x)+\nu_{P}(x))$ for any $x\in X$. $\mu_{P}(x)$, $\eta_{P}(x)$, and
$\nu_{P}(x)$ denote the \textit{degree of positive membership},
\textit{neutral membership}, and \textit{negative membership} of $x$ in $P$,
respectively. $\pi_{P}(x)$ is the \textit{degree of refusal membership} of
$x$ in $P$. In particular, the triplet $\langle \mu, \eta,
\nu\rangle$ satisfying that $\mu$, $\eta$, $\nu\in [0, 1]$ and
$\mu+\eta+\nu\leq 1$ is called a \textit{picture fuzzy number} (PFN).
For convenience, a PFN $\alpha$ is denoted by $\alpha=\langle \mu_{\alpha},
\eta_{\alpha}, \nu_{\alpha}\rangle$.
\end{definition}

In order to effectively distinguish all PFNs, Wu et al.~\cite{WZCLZY2021}
introduced the following linear order for PFNs by applying score function
and accuracy functions.

\begin{definition}[{\textrm{\protect\cite[Definition 2.11]{WZCLZY2021}}}]
\label{Wu-Order-Def}
Let $\alpha = \langle \mu_{\alpha}, \eta_{\alpha}, $ $\nu_{\alpha} \rangle$ be
a picture fuzzy number. Then, the \textit{score function} $S(\alpha)$ is defined as
$S(\alpha) = \mu_{\alpha}-\nu_{\alpha}$ and the \textit{first accuracy function}
$H_{1}(\alpha)$ and the \textit{second accuracy function} $H_{2}(\alpha)$ of $\alpha$
are defined as $H_{1}(\alpha) = \mu_{\alpha}+ \nu_{\alpha}$ and $H_{2}(\alpha) =
\mu_{\alpha}+\eta_{\alpha}+\nu_{\alpha}$, respectively. Then, for two PFNs $\alpha$ and $\beta$,
\begin{itemize}
  \item if $S(\alpha) < S(\beta)$, then $\alpha$ is smaller than $\beta$, denoted by $\alpha \prec_{_{\mathrm{W}}} \beta$;
  \item if $S(\alpha) = S(\beta)$, then
  \begin{itemize}
  \item if $H_{1}(\alpha) < H_{1}(\beta)$, then $\alpha$ is smaller than $\beta$, denoted by $\alpha \prec_{_\mathrm{{W}}} \beta$;
  \item if $H_{1}(\alpha) = H_{1}(\beta)$, then
  \begin{itemize}
  \item if $H_{2}(\alpha) = H_{2}(\beta)$, then $\alpha = \beta$;
  \item if $H_{2}(\alpha) < H_{2}(\beta)$, then $\alpha$ is smaller than $\beta$, denoted by $\alpha \prec_{_{\mathrm{W}}} \beta$.
\end{itemize}
\end{itemize}
\end{itemize}
If $\alpha \prec_{_{\mathrm{W}}} \beta$ or $\alpha =\beta$,
we will denote it by $\alpha \preccurlyeq_{_{\mathrm{W}}} \beta$.
\end{definition}

\subsection{Triangular norm}

\begin{definition}[{\textrm{\protect\cite[Definition 1.1]{KMP2000}}}]
A binary function $T:[0, 1]^{2} \longrightarrow [0, 1]$ is said to be
a \textit{triangular norm} (shortly, \textit{t-norm}) on $[0, 1]$ if,
for any $x$, $y$, $z\in [0, 1]$, the following conditions are fulfilled:
\begin{itemize}
\item[(T$_{1}$)] (\textit{Commutativity}) $T(x,y)=T(y,x)$;

\item[(T$_{2}$)] (\textit{Associativity}) $T(T(x, y)
,z) =T(x, T(y, z)) $;

\item[(T$_{3}$)] (\textit{Monotonicity}) If $x\leq y$, then $T(x, z)\leq T(y, z$);

\item[(T$_{4}$)] (\textit{Neutrality}) $T(1, x)=x$.
\end{itemize}
\end{definition}

Because of the associativity of a t-norm $T$, we can
extend $T$ to an $n$-ary function $T^{(n)}: [0, 1]^{n} \longrightarrow [0, 1]$ as follows
(see \cite[Remark 1.10]{KMP2000}):
$$
T^{(n)}(x_{1}, \ldots, x_{n-1}, x_{n}) \triangleq T\left(T^{(n-1)}
\left(x_{1}, \ldots, x_{n-1}\right), x_{n}\right).
$$
In particular, when $x_{1} = x_{2} = \cdots = x_{n} = x$, we briefly write
$$
x_{T}^{(n)} = T^{(n)}(x, x,\ldots, x) \ (n\geq 2),\ x_{T}^{(0)} = 1, \text{ and } x_{T}^{(1)}= x.
$$

\begin{definition}[{\textrm{\protect\cite[Definition 1.13]{KMP2000}}}]
A binary function $S: [0, 1]^{2}\longrightarrow [0, 1]$ is said to be a
\textit{triangular conorm} (shortly, \textit{t-conorm}) on $[0, 1]$ if,
for any $x$, $y$, $z\in [0, 1]$, it satisfies (T$_{1}$)--(T$_{3}$) and
\begin{itemize}
\item[(S$_{4}$)] (\textit{Neutrality})  $S(x,0)=x$.
\end{itemize}
\end{definition}

The duality expression of t-norms and t-conorms is the following result.
\begin{proposition}[{\textrm{\protect\cite[Proposition 1.15]{KMP2000}}}]
\label{Dual-Prop}
A binary function $S: [0, 1]^2\longrightarrow [0, 1]$ is a t-conorm if and only
if there exists a t-norm $T$ such that, for any $(x, y)\in [0, 1]^{2}$,
\begin{equation}
\label{3}
S(x,y)=1-T(1-x,1-y).
\end{equation}
\end{proposition}
The t-norm given by formula~\eqref{3} is called the \textit{dual t-conorm} of $T$, analogously,
we can give the definition of the dual t-norm of a t-conorm $S$.

\begin{definition}[{\textrm{\protect\cite[Definition 3.2]{KMP2000}}}]
\label{Pseudo-Inverse-Def}
Let $[a, b]$ and $[c, d]$ be two closed subintervals of the extended real line
$[-\infty, +\infty]$ and $f : [a, b] \longrightarrow [c, d]$ be a monotone function.
Then the \textit{pseudo-inverse} $f^{(-1)} : [c, d] \longrightarrow [a, b]$ of $f$ is defined by
$$
f^{(\!-1)}=\sup\left\{x \in [a, b] \mid (f(x)\!-y)(f(b)\!-f(a))\!<0\right\}.
$$
\end{definition}

\begin{definition}[{\textrm{\protect\cite[Definition 3.25]{KMP2000}}}]
\label{Add-Generator-Def}
An \textit{additive generator} $\tau: [0, 1] \longrightarrow [0, +\infty]$ of a t-norm $T$ is a
strictly decreasing function which is also right-continuous in $0$ and satisfies $\tau(1) = 0$,
such that for any
$(x, y) \in [0, 1]^{2}$, we have
$$
\tau(x) + \tau(y) \in \mathrm{Ran}(\tau)\cup [\tau(0), +\infty],
$$
and
$$
T(x, y) = \tau^{(-1)}(\tau(x) + \tau(y)).
$$
\end{definition}

\begin{definition}[{\textrm{\protect\cite[Definitions 2.9 and 2.13]{KMP2000}}}]
\label{Def 8}
A t-norm $T$ is said to be
\begin{enumerate}[(i)]
  \item \textit{strictly monotone} if
$T(x, y) < T(x, z)$ whenever $x > 0$ and $y < z$.
  \item \textit{strict} if it is continuous and strictly monotone.
\end{enumerate}
\end{definition}

\begin{lemma}[{\textrm{\protect\cite[Proposition~2.15, Theorem 5.1]{KMP2000}}}]
\label{Strict-T-Norm-Char}
For a function $T: [0, 1]^{2}\longrightarrow [0, 1]$, the following are equivalent:
\begin{enumerate}[{\rm (i)}]
\item $T$ is a strict t-norm;

\item $T$ has a continuous additive generator $\tau$ with $\tau(0)=+\infty$.
\end{enumerate}
\end{lemma}

\section{Basic operations for PFNs}\label{Sec-3}
This section presents some basic operations for PFNS, which originated from
the work of Wu et al.~\cite{WZCLZY2021}.

\begin{definition}[{\textrm{\protect\cite[Definition~5.1]{WZCLZY2021}}}]
\label{Baisc-Oper-T-Def}
\label{Baisc-Oper-T-Def}
Let $\alpha_{1} =\langle \mu_{1}, \eta_{1}, $ $\nu_{1}\rangle$,
$\alpha_{2} =\langle \mu_{2}, \eta_{2}, \nu_{2}\rangle$,
and $\alpha =\langle \mu, \eta, \nu\rangle$ be three PFNs,
$T$ be a strict t-norm with an additive generator $\tau$, $S$ be
the dual t-conorm of $T$. Define
\begin{enumerate}[(i)]
  \item
$\alpha_{1}\oplus_{_{T}}\alpha_{2} = \langle S(\mu_{1},\mu_{2}),
T(\eta_{1}+\nu_{1},\eta_{2}+\nu_{2})
-T(\nu_{1},\nu_{2}), T(\nu_{1},\nu_{2})\rangle$ ;

  \item
$\alpha_{1}\otimes_{_{T}}\alpha_{2} = \langle T(\mu_{1},\mu_{2}),
  T(\eta_{1}+\mu_{1},\eta_{2}+\mu_{2})-T(\mu_{1},\mu_{2}),
  S(\nu_{1},\nu_{2}) \rangle$;

  \item
$\lambda_{_{T}}\cdot\alpha = \langle \zeta^{-1}(\lambda\cdot\zeta(\mu)),
\tau^{-1}(\lambda\cdot\tau(\eta+\nu))-\tau^{-1}(\lambda\cdot\tau(\nu)),
\tau^{-1}(\lambda\cdot\tau(\nu))\rangle$,  $\lambda > 0$;

  \item
$\alpha^{\lambda_{_{T}}} =\langle \tau^{-1}(\lambda\cdot\tau(\mu)),
\tau^{-1}(\lambda\cdot\tau(\eta+\mu))-\tau^{-1}(\lambda\cdot\tau(\mu)),
\zeta^{-1}(\lambda\cdot\zeta(\nu))\rangle$, $\lambda > 0$;
\end{enumerate}
where $\zeta(x)=\tau(1-x)$, which is an additive generator of $S$.
If the t-norm $T$ can be clearly known from the context,
$\oplus_{_{T}}$, $\otimes_{_{T}}$, $\lambda_{_{T}}\cdot\alpha$, and $\alpha^{\lambda_{_{T}}}$ are
denoted as $\oplus$, $\otimes$, $\lambda\cdot\alpha$, and $\alpha^{\lambda}$, respectively.
\end{definition}

\begin{theorem}[{\textrm{\protect\cite[Theorem~5.1]{WZCLZY2021}}}]
\label{Closed-Thm}
Let  $\alpha_{1}=\langle\mu_{1},\eta_{1},\nu_{1}\rangle$, $\alpha_{2}=\langle\mu_{2},\eta_{2},\nu_{2}\rangle$,
and $\alpha=\langle\mu,\eta,\nu\rangle$ be three PFNs, and $T$ be a strict t-norm. Then, for any $\lambda > 0$,
all $\alpha_{1}\oplus_{_{T}}\alpha_{2}$, $\alpha_{1}\otimes_{_{T}}\alpha_{2}$, $\lambda_{_{T}}\alpha$
and $\alpha^{\lambda_{_{T}}}$ are PFNs.
\end{theorem}

\begin{theorem}[{\textrm{\protect\cite[Theorem~5.2]{WZCLZY2021}}}]
\label{Operation-Properties-Thm}
Let $\alpha$, $\beta$, $\gamma$ be three PFNs, and $T$ be a strict t-norm.
Then, for any $\lambda$, $\xi > 0$, we have
\begin{enumerate}[{\rm (i)}]
  \item $\alpha\oplus_{_{T}}\beta=\beta\oplus_{_{T}}\alpha$;
  \item $\alpha\otimes_{_{T}}\beta=\beta\otimes_{_{T}}\alpha$;
  \item $\left(\alpha\oplus_{_{T}}\beta\right)\oplus_{_{T}}\gamma
  =\alpha\oplus_{_{T}}\left(\beta\oplus_{_{T}}\gamma\right)$;
  \item $\left(\alpha\otimes_{_{T}}\beta\right)\otimes_{_{T}}\gamma
  =\alpha\otimes_{_{T}}\left(\beta\otimes_{_{T}}\gamma\right)$;
  \item $\left(\xi_{_{T}}\alpha\right)\oplus_{_{T}}\left(\lambda_{_{T}}
  \alpha\right)=\left(\xi+\lambda\right)_{_{T}}\alpha$;
  \item $\left(\alpha^{\xi_{_{T}}}\right)\otimes_{_{T}}\left(\alpha^{\lambda_{_{T}}}\right)
  =\alpha^{\left(\xi+\lambda\right)_{_{T}}}$;
  \item $\lambda_{_{T}}\cdot\left(\alpha\oplus_{_{T}}\beta\right)=\left(\lambda_{_{T}}\cdot\alpha\right)
      \oplus_{_{T}}\left(\lambda_{_{T}}\cdot\beta\right)$;
  \item $\left(\alpha\otimes_{_{T}}\beta\right)^{\lambda_{_{T}}}=\alpha^{\lambda_{_{T}}}
  \otimes_{_{T}}\beta^{\lambda_{_{T}}}$;
  \item $\xi_{_{T}}\cdot\left(\lambda_{_{T}}\cdot\alpha\right)=\left(\lambda\cdot\xi\right)_{_{T}}\cdot\alpha$;
  \item $\left(\alpha^{\lambda_{_{T}}}\right)^{\xi_{_{T}}}=\alpha^{\left(\lambda\cdot\xi\right)_{_{T}}}$.
\end{enumerate}
\end{theorem}

\begin{theorem}[{\textrm{\protect\cite[Theorem~5.3]{WZCLZY2021}}}]
\label{Formula-Thm}
Let $\alpha_{i}=\langle \mu_{i}, \eta_{i}, \nu_{i} \rangle$
($i = 1, 2, \ldots, n$) be a collection of PFNs and $T$ be a strict t-norm. Then,
\begin{equation}\label{4}
\begin{aligned}
&\alpha_{1}\oplus_{_{T}}\alpha_{2}\oplus_{_{T}} \cdots \oplus_{_{T}}\alpha_{n}\\
=& \left\langle S^{(n)}(\mu_{1}, \ldots , \mu_{n}), T^{(n)}
(\eta_{1}\!+\nu_{1}, \ldots ,\eta_{n}\!+\nu_{n})\right.
\\
&\left.\quad \!-T^{(n)}(\nu_{1}, \ldots , \nu_{n}), T^{(n)}
(\nu_{1}, \ldots , \nu_{n})\right\rangle,
\end{aligned}
\end{equation}
and
\begin{equation}\label{5}
\begin{aligned}
&\alpha_{1}\otimes_{_{T}}\alpha_{2}\otimes_{_{T}} \cdots \otimes_{_{T}}\alpha_{n}\\
=&\left\langle T^{(n)}(\mu_{1}, \ldots , \mu_{n}), T^{(n)}
(\eta_{1}\!+\mu_{1}, \ldots , \eta_{n}\!+\mu_{n})\right.
\\
&\left.\quad \!-T^{(n)}(\mu_{1}, \ldots , \mu_{n}), S^{(n)}
(\nu_{1},\ldots , \nu_{n})\right\rangle,
\end{aligned}
\end{equation}
where $S$ is the dual t-conorm of $T$.
\end{theorem}

\begin{remark}\label{Remark 1}
Let $\alpha_{i}=\langle \mu_{i}, \eta_{i}, \nu_{i} \rangle$
($i = 1, 2, \ldots, n$) be a collection of PFNs and $T$ be a strict t-norm with an
additive generator $\tau$. By Theorem~\ref{Formula-Thm},
it can be verified that
{\small \begin{equation}\label{6}
\begin{aligned}
 \bigoplus_{i=1}^{n}\alpha_{i}
 \triangleq& \alpha_{1}\oplus_{_{T}}\alpha_{2}\oplus_{_{T}} \cdots \oplus_{_{T}}\alpha_{n}\\
 =&\Bigg\langle \zeta^{\!-1}\left(\sum\limits_{i\!=1}^{n}\zeta\left(\mu_{\alpha_{i}}\right)\right), \tau^{-1}\left(\sum\limits_{i\!=1}^{n}\tau\left(\eta_{\alpha_{i}}\!+\nu_{\alpha_{i}}\right)\right)\\
 &\quad \!-\tau^{\!-1}\left(\sum\limits_{i\!=1}^{n}\tau\left(\nu_{\alpha_{i}}\right)\right),
 \tau^{\!-1}\left(\sum\limits_{i\!=1}^{n}\tau\left(\nu_{\alpha_{i}}\right)\right)\Bigg\rangle,
 \end{aligned}
 \end{equation}
 }
 and
{\small \begin{equation}\label{7}
\begin{aligned}
 \bigotimes_{i=1}^{n}\alpha_{i}
 \triangleq& \alpha_{1}\otimes_{_{T}}\alpha_{2}\otimes_{_{T}} \cdots \otimes_{_{T}}\alpha_{n}\\
 =&\Bigg\langle \tau^{-1}\left(\sum\limits_{i=1}^{n}\tau\left(\mu_{\alpha_{i}}\right)\right), \tau^{-1}\left(\sum\limits_{i=1}^{n}\tau\left(\eta_{\alpha_{i}}\!+\mu_{\alpha_{i}}\right)\right)\\
 &\quad \!-\tau^{-1}\left(\sum\limits_{i=1}^{n}\tau\left(\mu_{\alpha_{i}}\right)\right),
 \zeta^{-1}\left(\sum\limits_{i=1}^{n}\zeta\left(\nu_{\alpha_{i}}\right)\right)\Bigg\rangle.
 \end{aligned}
 \end{equation}
 }
 In particular, when $\alpha_{i}=\alpha=\left\langle \mu, \eta, \nu \right\rangle$
 ($i = 1, 2,\ldots, n$),
 we have
$\alpha\oplus_{_{T}}\alpha\oplus_{_{T}} \cdots \oplus_{_{T}}\alpha
 = \langle \zeta^{-1}(n\cdot\zeta(\mu_{\alpha})), \tau^{-1}(n\cdot\tau(\eta_{\alpha}+\nu_{\alpha}))
-\tau^{-1}(n\cdot\tau(\nu_{\alpha})),
 \tau^{-1}(n\cdot\tau(\nu_{\alpha}))\rangle
 =n\cdot \alpha$
 and
$\alpha\otimes_{_{T}}\alpha\otimes_{_{T}} \cdots \otimes_{_{T}}\alpha=
\langle \tau^{-1}(n\cdot\tau(\mu_{\alpha})), \tau^{-1}(n\cdot\tau(\eta_{\alpha}+\mu_{\alpha}))
-\tau^{-1}(n\cdot\tau(\mu_{\alpha})), \zeta^{-1}(n\cdot\zeta(\nu_{\alpha}))\rangle=\alpha^{n}$.
\end{remark}

\section{Picture fuzzy interactional Bonferroni mean (PFIBM) operator}\label{Sec-4}

Combining the basic operational laws defined in Definition~\ref{Baisc-Oper-T-Def}
with BM operator in \cite{Bo1950}, this section introduces an aggregation operator called
Picture fuzzy interactional Bonferroni mean (PFIBM) for aggregating the picture fuzzy
information and proves the monotonicity, idempotency, boundedness, and
commutativity for the PFIBM operator.

\begin{definition}
\label{PFIBM-Def}
Let $\alpha_{i}$ ($i=1,2,\ldots,n$) be a collection of PFNs and $T$ be a t-norm. For $p$, $q > 0$, define
the \textit{picture fuzzy interactional Bonferroni mean} (PFIBM) operator induced by $T$ as follows:
{\small$$
\mathrm{PFIBM}_{T}^{p,q}(\alpha_{1}, \ldots ,\alpha_{n})
=\left[\frac {1}{n(n\!-1)}
\bigoplus\limits_{i,j\!=1 \atop i\neq j}^{n}\left(\alpha_{i}^{p}
\otimes\alpha_{j}^{q}\right)\right]^{\frac{1}{p\!+q}}.
$$}
\end{definition}

\begin{theorem}
\label{PFIBM-Formula-Thm}
Let $\alpha_{i}=\left\langle\mu_{\alpha_{i}},\eta_{\alpha_{i}},\nu_{\alpha_{i}}\right\rangle$
($i=1,2,\ldots,n$) be a collection of PFNs and $T$ be a strict t-norm with an additive generator $\tau$. Then,
for $p, q > 0$, the aggregated value by using the PFIBM induced by $T$ is also a PFN, and
\begin{equation}\label{11}
\begin{aligned}
&\mathrm{PFIBM}_{T}^{p,q}\left(\alpha_{1}, \ldots ,\alpha_{n}\right)\\
=&\Bigg\langle\tau^{\!-1}\left(\frac{1}{p\!+q}\cdot\tau\left(\mu\right)\right),
\tau^{\!-1}\left(\frac{1}{p\!+q}\cdot\tau\left(\eta\!+\mu\right)\right)\\
&\quad\!-\tau^{\!-1}\left(\frac{1}{p+q}\cdot\tau\left(\mu\right)\right),
\zeta^{\!-1}\left(\frac{1}{p\!+q}\cdot\zeta\left(\nu\right)\right)\Bigg\rangle,
\end{aligned}
\end{equation}
where
{$\mu \!= \zeta^{\!-1}\bigg(\frac{1}{n(n\!-1)}\!\sum\limits_{i,j\!=1 \atop i\neq j}^{n}
\zeta\Big(\tau^{\!-1}\left(p\!\cdot\tau\left(\mu_{\alpha_{i}}\right)
\!+q\!\cdot\tau\left(\mu_{\alpha_{j}}\right)\right)\Big)\bigg)$,
$\eta=\tau^{\!-1}\bigg(\frac{1}{n(n\!-1)} \sum\limits_{i,j\!=1 \atop i\neq j}^{n}
\tau\Big(\tau^{\!-1}(p\cdot\tau(\eta_{\alpha_{i}}\!+\mu_{\alpha_{i}})
\!+q\cdot\tau(\eta_{\alpha_{j}}\!+\mu_{\alpha_{j}}))
\!-\tau^{\!-1}(p\cdot\tau(\mu_{\alpha_{i}})\!+q\cdot\tau(\mu_{\alpha_{j}}))
\!+\zeta^{\!-1}(p\cdot\zeta(\nu_{\alpha_{i}})\!+q\cdot\zeta(\nu_{\alpha_{j}}))\Big)\bigg)
\!-\tau^{\!-1}\bigg(\frac{1}{n(n\!-1)} \sum\limits_{i,j\!=1 \atop i\neq j}^{n}
\tau\Big(\zeta^{\!-1}(p\cdot\zeta(\nu_{\alpha_{i}})\!+q\cdot\zeta
(\nu_{\alpha_{j}}))\Big)\bigg)$,
$\nu= \tau^{\!-1}\bigg(\frac{1}{n(n\!-1)} \sum\limits_{i,j\!=1 \atop i\neq j}^{n}
\tau\Big(\zeta^{\!-1}\left(p\cdot\zeta\left(\nu_{\alpha_{i}}\right)
\!+q\cdot\zeta\left(\nu_{\alpha_{j}}\right)\right)\Big)\bigg)$},
and $\zeta(x)=\tau(1\!-x)$.
\end{theorem}

\begin{IEEEproof}
First, by Theorem~\ref{Closed-Thm}, we know that all four operations
$\oplus$, $\otimes$, $\lambda \alpha$, and $\alpha^{\lambda}$ are closed in PFNs.
Thus, $\mathrm{PFIBM}_{T}^{p,q}\left(\alpha_{1}, \alpha_{2},\ldots,
\alpha_{n}\right)$ is a PFN by Definition~\ref{PFIBM-Def}.

Second, we prove the formula~\eqref{11} holds. By the operational laws (ii) and (iv)
defined in Definition~\ref{Baisc-Oper-T-Def}, we have that, for $1\leq i, j\leq n$,
\begin{align*}
\alpha_{i}^{p}=&\big\langle \tau^{-1}(p\cdot\tau(\mu_{\alpha_{i}})),
\tau^{-1}(p\cdot\tau(\eta_{\alpha_{i}}+\mu_{\alpha_{i}}))\\
&\quad-\tau^{-1}(p\cdot\tau(\mu_{\alpha_{i}})),\zeta^{-1}(p\cdot\zeta(\nu_{\alpha_{i}}))\big\rangle,
\end{align*}
and
\begin{align*}
\alpha_{j}^{q}=&\big\langle \tau^{-1}(q\cdot\tau(\mu_{\alpha_{j}})), \tau^{-1}
(q\cdot\tau(\eta_{\alpha_{j}}+\mu_{\alpha_{j}}))\\
&\quad-\tau^{-1}(q\cdot\tau(\mu_{\alpha_{j}})),\zeta^{-1}(q\cdot\zeta(\nu_{\alpha_{j}}))\big\rangle,
\end{align*}
and thus
    {\small \begin{equation}\label{15}
    \begin{aligned}
    &\alpha_{i}^{p}\otimes\alpha_{j}^{q}\\
    =& \bigg\langle T\Big(\tau^{-1}(p\cdot\tau(\mu_{\alpha_{i}})),
    \tau^{-1}(q\cdot\tau(\mu_{\alpha_{j}}))\Big),\bigg. \\
    &\quad T\Big(\tau^{-1}\left(p\cdot\tau(\eta_{\alpha_{i}}+\mu_{\alpha_{i}})\right),
    \tau^{-1}\left(q\cdot\tau(\eta_{\alpha_{j}}+\mu_{\alpha_{j}})\right)\Big)\\
    &\quad-T\Big(\tau^{-1}\left(p\cdot\tau(\mu_{\alpha_{i}})\right),
    \tau^{-1}\left(q\cdot\tau(\mu_{\alpha_{j}})\right)\Big),\\
    &\quad\bigg. S\Big(\zeta^{-1}\left(p\cdot\zeta\left(\nu_{\alpha_{i}}\right)\right),
    \zeta^{-1}\left(q\cdot\zeta\left(\nu_{\alpha_{j}}\right)\right)\Big)\bigg\rangle\\
    =& \bigg\langle \tau^{-1}\Big(p\cdot\tau\left(\mu_{\alpha_{i}}\right)
    +q\cdot\tau\left(\mu_{\alpha_{j}}\right)\Big),\bigg.\\
     &\quad \tau^{-1}\Big(p\cdot\tau\left(\eta_{\alpha_{i}}+\mu_{\alpha_{i}}\right)
     +q\cdot\tau\left(\eta_{\alpha_{j}}+\mu_{\alpha_{j}}\right)\Big)\\
   &\quad-\tau^{-1}\Big(p\cdot\tau\left(\mu_{\alpha_{i}}\right)
   +q\cdot\tau\left(\mu_{\alpha_{j}}\right)\Big),\\
   &\quad \bigg.\zeta^{-1}\Big(p\cdot\zeta\left(\nu_{\alpha_{i}}\right)
   +q\cdot\zeta\left(\nu_{\alpha_{j}}\right)\Big)\bigg\rangle.
    \end{aligned}
    \end{equation}
    }

    Let $\beta_{ij}=\left\langle\mu_{\beta_{ij}},\eta_{\beta_{ij}},\nu_{\beta_{ij}}\right\rangle
    =\alpha_{i}^{p}\otimes\alpha_{j}^{q}$. Then, by Theorem~\ref{Formula-Thm} or formula~\eqref{6}, we have
    {\small \begin{equation}\label{16}
    \begin{aligned}
    \bigoplus_{i,j=1 \atop i\neq j}^{n}\beta_{ij}
    = &\left\langle \zeta^{-1}\Big(\sum\limits_{i,j=1 \atop i\neq j}^{n}\zeta(\mu_{\beta_{ij}})\Big),
     \tau^{-1}\Big(\sum\limits_{i,j=1 \atop i\neq j}^{n}\tau(\eta_{\beta_{ij}}+\nu_{\beta_{ij}})\Big)\right.\\
    &\quad\left.-\tau^{-1}\Big(\sum\limits_{i,j=1 \atop i\neq j}^{n}\tau(\nu_{\beta_{ij}})\Big),
     \tau^{-1}\Big(\sum\limits_{i,j=1 \atop i\neq j}^{n}\tau(\nu_{\beta_{ij}})\Big)\right\rangle.
    \end{aligned}
    \end{equation}
    }

    This, together with Definition~\ref{Baisc-Oper-T-Def} (iii), implies that
  {\small
    \begin{align*}
    & \frac{1}{n(n-1)}\bigoplus_{i,j=1 \atop i\neq j}^{n}\left(\alpha_{i}^{p}\otimes\alpha_{j}^{q}\right)
    =\frac{1}{n(n-1)}\bigoplus_{i,j=1 \atop i\neq j}^{n}\beta_{ij}\\
    =&\left\langle \zeta^{-1}\Bigg(\frac{1}{n(n-1)} \sum\limits_{i,j=1 \atop i\neq
    j}^{n}\zeta(\mu_{\beta_{ij}})\Bigg), \tau^{-1}\Bigg(\frac{1}{n(n-1)} \sum\limits_{i,j=1 \atop i\neq j}^{n}\right.\\
    & \quad \tau(\eta_{\beta_{ij}}+\nu_{\beta_{ij}})\Bigg)
    -\tau^{-1}\Bigg(\frac{1}{n(n-1)} \sum\limits_{i,j=1 \atop i\neq j}^{n}\tau(\nu_{\beta_{ij}})\Bigg),\\
    &\quad\left.\tau^{-1}\Bigg(\frac{1}{n(n-1)} \sum\limits_{i,j=1 \atop i\neq j}^{n}\tau(\nu_{\beta_{ij}})\Bigg)\right\rangle\\
    =&\left\langle \zeta^{-1}\Bigg(\frac{1}{n(n-1)} \sum\limits_{i,j=1 \atop i\neq j}^{n}\zeta\left(\tau^{-1}(p\cdot\tau(\mu_{\alpha_{i}})+q\cdot\tau(\mu_{\alpha_{j}}))\right)\Bigg)\right.,\\
    &\quad\tau^{-1}\Bigg(\frac{1}{n(n-1)} \sum\limits_{i,j=1 \atop i\neq j}^{n}
    \tau\Big(\tau^{-1}(p\cdot\tau(\eta_{\alpha_{i}}+\mu_{\alpha_{i}})\Big.\Bigg.\\
    &\quad \quad \quad +q\cdot\tau(\eta_{\alpha_{j}}+\mu_{\alpha_{j}}))-\tau^{-1}(p\cdot
    \tau(\mu_{\alpha_{i}})+q\cdot\tau(\mu_{\alpha_{j}}))\\
    &\Bigg.\Big.\quad\quad\quad+\zeta^{-1}(p\cdot\zeta(\nu_{\alpha_{i}})
    +q\cdot\zeta(\nu_{\alpha_{j}}))\Big)\Bigg)\\
    &\quad-\tau^{-1}\Bigg(\frac{1}{n(n-1)} \sum\limits_{i,j=1 \atop i\neq j}^{n}\tau\Big(\zeta^{-1}(p\cdot\zeta(\nu_{\alpha_{i}})+q\cdot\zeta(\nu_{\alpha_{j}}))\Big)\Bigg),\\
    &\left.\quad\tau^{-1}\Bigg(\frac{1}{n(n-1)} \sum\limits_{i,j=1 \atop i\neq j}^{n}\tau\Big(\zeta^{-1}
    \left(p\cdot\zeta(\nu_{\alpha_{i}})+q\cdot
    \zeta(\nu_{\alpha_{j}})\right)\Big)\Bigg)\right\rangle.
    \end{align*}}

    Let $\gamma=\langle\mu, \eta, \nu\rangle
    =\frac{1}{n(n-1)}\bigoplus\limits_{i,j=1 \atop i\neq j}^{n}\left(\alpha_{i}^{p}
    \otimes\alpha_{j}^{q}\right)$. Then, by Definition~\ref{Baisc-Oper-T-Def} (iv), we have
    \begin{align*}
& \mathrm{PFIBM}_{T}^{p,q}(\alpha_{1}, \alpha_{2},\ldots ,\alpha_{n})=\gamma^{\frac{1}{p+q}}\\
=&\left\langle\tau^{-1}\left(\frac{1}{p+q}\cdot\tau\left(\mu\right)\right)
\tau^{-1}\left(\frac{1}{p+q}\cdot\tau\left(\eta+\mu\right)\right)\right.\\
&\quad \left.-\tau^{-1}\left(\frac{1}{p+q}\cdot\tau\left(\mu\right)\right),
\zeta^{-1}\left(\frac{1}{p+q}\cdot\zeta\left(\nu\right)\right)\right\rangle,
\end{align*}
where
$\mu, \eta, \nu$ are equal to the values in Theorem~\ref{PFIBM-Formula-Thm}.
\end{IEEEproof}

{\begin{remark}
If we take $\eta_{\alpha_i}=0$ for each $\alpha_i$, i.e., each $\alpha_i$ reduces to an
intuitionistic fuzzy number, then
\begin{equation}
\label{eq-T-formula}
\begin{split}
&\mathrm{PFIBM}_{T}^{p,q}\left(\alpha_{1}, \alpha_{2},\ldots ,\alpha_{n}\right)\\
=&\left\langle\tau^{-1}\left(\frac{1}{p+q}\cdot\tau\left(\mu\right)\right), ~0,~
 \zeta^{-1}\left(\frac{1}{p+q}\cdot\zeta\left(\nu\right)\right)\right\rangle,
 \end{split}
\end{equation}
which is consistent with \cite[Theorem~1]{Ni2016}.
This means that \cite[Theorem~1]{Ni2016} is direct corollary
of Theorem~\ref{PFIBM-Formula-Thm}.
\end{remark}
}

If we use some specific t-norms to Theorem \ref{PFIBM-Formula-Thm},
we can obtain the following results.

(1) If $T$ = $T_{\mathbf{P}}$ (algebraic product),
  noting that the additive generator of $T_{\mathbf{P}}$
  is $-\log x$, then
{\small
\begin{align*}
&\mathrm{PFIBM}_{T_{\mathbf{P}}}^{p,q}
\left(\alpha_{1}, \alpha_{2},\ldots ,\alpha_{n}\right)\\
\!=&\left\langle\Bigg[1\!-\Bigg(\prod\limits_{i,j\!=1 \atop i\neq j}^{n}
\left(1\!-\mu_{\alpha_{i}}^{p}\mu_{\alpha_{j}}^{q}\right)
\Bigg)^{\frac{1}{n(n\!-1)}}\Bigg]^{\frac{1}{p+q}},\right. \\
&\Bigg[\Bigg(\prod\limits_{i,j\!=1 \atop i\neq j}^{n}
\left(\left(\eta_{\alpha_{i}}\!+\mu_{\alpha_{i}}\right)^{p}
\left(\eta_{\alpha_{j}}\!+\mu_{\alpha_{j}}\right)^{q}
\!-\mu_{\alpha_{i}}^{p}\mu_{\alpha_{j}}^{q}\right.\Bigg.\Bigg.\\
&\Bigg.\left.\!+1\!-\left(1-\nu_{\alpha_{i}}\right)^{p}\left(1\!-
\nu_{\alpha_{j}}\right)^{q}\right)\Bigg)^{\frac{1}{n(n\!-1)}}
\\
&\!-\Bigg(\prod\limits_{i,j\!=1 \atop i\neq j}^{n}
\left(1\!-\left(1\!-\nu_{\alpha_{i}}\right)^{p}
\left(1\!-\nu_{\alpha_{j}}\right)^{q}\right)\Bigg)^{\frac{1}{n(n\!-1)}}\\
&\Bigg.\!+1\!-\Bigg(\prod\limits_{i,j\!=1 \atop i\neq j}^{n}
\left(1\!-\mu_{\alpha_{i}}^{p}\mu_{\alpha_{i}}^{q}\right)
\Bigg)^{\frac{1}{n(n\!-1)}}\Bigg]^{\frac{1}{p\!+q}}
\\
&\Bigg.\!-\Bigg[1\!-\Bigg(\prod\limits_{i,j\!=1 \atop i\neq j}^{n}
\left(1\!-\mu_{\alpha_{i}}^{p}\mu_{\alpha_{j}}^{q}\right)
\Bigg)^{\frac{1}{n(n\!-1)}}\Bigg]^{\frac{1}{p\!+q}},\Bigg.
\\
&\left.1\!-\Bigg[1\!-\Bigg(\prod\limits_{i,j=1 \atop i\neq j}^{n}
\left(1\!-\left(1\!-\nu_{\alpha_{i}}\right)^{p}
\left(1\!-\nu_{\alpha_{j}}\right)^{q}\right)
\Bigg)^{\frac{1}{n(n\!-1)}}\Bigg]^{\frac{1}{p\!+q}} \right\rangle.
\end{align*}
}

{In particular, if we take $\eta_{\alpha_i}=0$
for each $\alpha_i$, i.e., each $\alpha_i$ reduces to an
intuitionistic fuzzy number, then
\begin{equation}
\label{eq-Tp}
\begin{split}
&\quad \mathrm{PFIBM}_{T_{\mathbf{P}}}^{p,q}
(\alpha_{1}, \alpha_{2},\ldots ,\alpha_{n})\\
& = \left\langle\Bigg[1\!-\Bigg(\prod\limits_{i,j=1 \atop i\neq j}^{n}
\left(1\!-\mu_{\alpha_{i}}^{p}\mu_{\alpha_{j}}^{q}\right)
\Bigg)^{\frac{1}{n(n\!-1)}}\Bigg]^{\frac{1}{p+q}}, ~0,~ 1\!- \right.
\\
&\left. \Bigg[1\!-\Bigg(\prod\limits_{i,j=1 \atop i\neq j}^{n}
\left(1\!-\left(1\!-\nu_{\alpha_{i}}\right)^{p}
\left(1\!-\nu_{\alpha_{j}}\right)^{q}\right)\Bigg)^{\frac{1}{n(n\!-1)}}
\Bigg]^{\frac{1}{p\!+q}} \right\rangle,
\end{split}
\end{equation}
which is consistent with \cite[Theorem~1]{XY2011} and \cite[Theorem~1.4.1]{XC2012}.
This means that \cite[Theorem~1]{XY2011} and \cite[Theorem~1.4.1]{XC2012} 
are direct corollaries of Theorem~\ref{PFIBM-Formula-Thm}.}

(2) If we take $T$ as Schweizer-Sklar t-norms
$T_{\gamma}^{\mathbf{SS}}$ ($\gamma \in (\!-\infty, 0)$), then
$\mathrm{PFIBM}_{T_{\gamma}^{\mathbf{SS}}}^{p,q}
(\alpha_{1}, \alpha_{2},\ldots ,\alpha_{n})=
\big \langle(1\!-\frac{1}{p\!+q} (1\!-\mu^{\gamma} ))^{\frac{1}{\gamma}},
(1\!-\frac{1}{p\!+q}(1\!-(\eta+\mu )^{\gamma} ) )^{\frac{1}{\gamma}}
\!-(1\!-\frac{1}{p\!+q} (1\!-\mu^{\gamma} ))^{\frac{1}{\gamma}},
1\!-(1\!-\frac{1}{p\!+q}(1\!-(1\!-\nu)^{\gamma}))^{\frac{1}{\gamma}}\big\rangle$,
where
$\mu=1\!-\bigg[\frac{1}{n(n\!-1)}\sum\limits_{i,j\!=1 \atop i\neq j}^{n}
 \bigg(1\!-\left(1\!-p\left(1\!-\mu_{\alpha_{i}}^{\gamma}\right)
 \!-q\left(1\!-\mu_{\alpha_{j}}^{\gamma}\right)\right)
^{\frac{1}{\gamma}}\bigg)^{\gamma}\bigg]^{\frac{1}{\gamma}}$,
{\small \begin{align*}
\eta
=& \Bigg\{\frac{1}{n(n\!-1)}\sum\limits_{i,j\!=1 \atop i\neq j}^{n}\bigg[
\left(1\!-p\left(1\!-\left(\eta_{\alpha_{i}}\!+\mu_{\alpha_{i}}\right)^{\gamma}\right)\right.\bigg.\Bigg.\\
&\quad\left.\!-q\left(1\!-\left(\eta_{\alpha_{j}}\!+\mu_{\alpha_{j}}\right)^{\gamma}\right)\right)^{\frac{1}{\gamma}}
\!-\left(1\!-p\left(1\!-\mu_{\alpha_{i}}^{\gamma}\right)\right.\\
&\quad\left.\!-q\left(1\!-\mu_{\alpha_{j}}^{\gamma}\right)\right)
^{\frac{1}{\gamma}}\!+1\!-\left(1\!-p\left(1\!-\left(1\!-\nu_{\alpha_{i}}\right)^{\gamma}\right)\right.\\
&\quad\Bigg.\bigg.\left.\!-q\left(1\!-\left(1\!-\nu_{\alpha_{j}}\right)
^{\gamma}\right)\right)^{\frac{1}{\gamma}}\bigg]^{\gamma}\Bigg\}^{\frac{1}{\gamma}}\\
&\!-\Bigg\{\frac{1}{n(n\!-1)}\sum\limits_{i,j\!=1 \atop i\neq j}^{n}\bigg[1\!-\left(1\!-p\left(1\!-\left(1\!-\nu_{\alpha_{i}}\right)^{\gamma}\right)\right.\bigg.\Bigg.\\
&\quad\quad\Bigg.\bigg.\left.\!-q\left(1\!-\left(1\!-\nu_{\alpha_{j}}\right)^{\gamma}\right)
\right)^{\frac{1}{\gamma}}\bigg]^{\gamma}\Bigg\}
^{\frac{1}{\gamma}},
\end{align*}
and
\begin{align*}
\nu \!=& \Bigg[\frac{1}{n(n\!-1)}\sum\limits_{i,j\!=1 \atop i\neq j}^{n}\bigg(1\!-\left(1\!-p\left(1\!-\left(1\!-\nu_{\alpha_{i}}\right)^{\gamma}\right)\right.\bigg.\Bigg.\\
&\quad\quad\Bigg.\bigg.\left.\!-q\left(1\!-\left(1-\nu_{\alpha_{j}}\right)^{\gamma}
\right)\right)^{\frac{1}{\gamma}}\bigg)^{\gamma}\Bigg]^{\frac{1}{\gamma}}.
\end{align*}
}

 (3) If we take $T$ as Hamacher t-norms $T_{\gamma}^{\mathbf{H}}$~($\gamma \in (0, +\infty)$),\\
then
$\mathrm{PFIBM}_{T_{\gamma}^{\mathbf{H}}}^{p,q}(\alpha_{1}, \alpha_{2},\ldots ,\alpha_{n})
\!=\Big\langle\frac{\gamma}{\left(\frac{\gamma}{\mu}
\!+1\!-\gamma\right)^{\frac{1}{p\!+q}}\!+\gamma\!-1},$
$\frac{\gamma}{\left(\frac{\gamma}{\eta\!+\mu}\!+1\!-\gamma\right)^{\frac{1}{p\!+q}}\!+\gamma\!-1}
\!-\frac{\gamma}{\left(\frac{\gamma}{\mu}\!+1\!-\gamma\right)^{\frac{1}{p\!+q}}\!+\gamma\!-1},
1\!-\frac{\gamma}{\left(\frac{\gamma}{1\!-\nu}
\!+1\!-\gamma\right)^{\frac{1}{p\!+q}}\!+\gamma\!-1}\Big\rangle$,
where
\begin{align*}
 \mu= 1\!-\frac{\gamma}{\left[\prod\limits_{i,j=1 \atop i\neq j}^{n}
 \left(\frac{\gamma^{2}}{\overline{\eta}_2\!-\gamma}\!+1\right)\right]
 ^{\frac{1}{n(n\!-1)}}\!+\gamma\!-1},
\end{align*}
\begin{align*}
\eta=&\frac{\gamma}{\left[\prod \limits_{i,j=1 \atop i\neq j}^{n}
\left(\frac{\gamma}{\frac{\gamma}{\overline{\eta}_1}
\!-\frac{\gamma}{\overline{\eta}_2}
\!+1\!-\frac{\gamma}{\overline{\eta}_3}}\!+1\!-\gamma\right)\right]^{\frac{1}{n(n\!-1)}}\!+\gamma\!-1}\\
& \!-\frac{\gamma}{\left[\prod\limits_{i,j=1 \atop i\neq j}^{n}
\left(\frac{\gamma^{2}}{\overline{\eta}_3\!-\gamma}\!+1\right)\right]
 ^{\frac{1}{n(n\!-1)}}\!+\gamma\!-1},
\end{align*}
\begin{align*}
\nu=\frac{\gamma}{\left[\prod\limits_{i,j=1 \atop i\neq j}^{n}
\left(\frac{\gamma^{2}}{\overline{\eta}_3\!-\gamma}\!+1\right)\right]
 ^{\frac{1}{n(n\!-1)}}\!+\gamma\!-1},
\end{align*}
and
$\overline{\eta}_1=\left(\frac{\gamma}{\eta_{\alpha_{i}}+\mu_{\alpha_{i}}}+1-\gamma\right)^{p}\left(\frac{
\gamma}{\eta_{\alpha_{j}}+\mu_{\alpha_{j}}}+1-\gamma\right)^{q}+\gamma-1$,
$\overline{\eta}_2=\left(\frac{\gamma}{\mu_{\alpha_{i}}}+1-\gamma\right)^{p}\left(\frac{
\gamma}{\mu_{\alpha_{j}}}+1-\gamma\right)^{q}+\gamma-1$,
$\overline{\eta}_3=\Big(\frac{\gamma}{1-\nu_{\alpha_{i}}}+1-\gamma\Big)^{p}\left(\frac{
\gamma}{1-\nu_{\alpha_{j}}}+1-\gamma\right)^{q}+\gamma-1$.

  (4) If we take $T$ as Frank t-norms $T_{\gamma}^{\mathbf{F}}$
  ($\gamma \in (0, 1) \cup (1, \!+\infty)$), then
$\mathrm{PFIBM}_{T_{\gamma}^{\mathbf{F}}}^{p,q}
(\alpha_{1}, \alpha_{2},\ldots ,\alpha_{n})
\!=\Bigg \langle \log_{\gamma}\left[\frac{\gamma\!-1}{\left(\frac{\gamma\!-1}
{\gamma^{\mu}\!-1}\right)^{\frac{1}{p\!+q}}}\!+1\right],
\log_{\gamma}\left[\frac{\gamma\!-1}{\left(\frac{\gamma\!-1}
{\gamma^{\eta\!+\mu}\!-1}\right)^{\frac{1}{p\!+q}}}\!+1\right]
\!-\log_{\gamma}\left[\frac{\gamma\!-1}{\left(\frac{\gamma\!-1}
{\gamma^{\mu}\!-1}\right)^{\frac{1}{p\!+q}}}\!+1\right],
1\!-\log_{\gamma}\left[\frac{\gamma\!-1}{\left(\frac{\gamma\!-1}
{\gamma^{1\!-\nu}\!-1}\right)^{\frac{1}{p\!+q}}}\!+1\right]\Bigg \rangle$,
where
$\mu\!= 1\!-\log_{\gamma}\left\{\frac{\gamma\!-1}{\left[\prod
 \limits_{i,j\!=1 \atop i\neq j}^{n}\left(\frac{\gamma\!-1}{\gamma^{1\!-\overline{\eta}_2}
 \!-1}\right)\right]^{\frac{1}{n(n\!-1)}}}\!+1\right\}$,
$\eta\!=\log_{\gamma}\left(\frac{\gamma\!-1}{\overline{\eta}}\!+1\right)
\!-\log_{\gamma}\left\{\frac{\gamma\!-1}{\left[\prod\limits_{i,j\!=1
 \atop i\neq j}^{n}\left(\frac{\gamma\!-1}{\gamma^{1\!-\overline{\eta}_3}
 \!-1}\right)\right]^{\frac{1}{n(n\!-1)}}}\!+1\right\}$,
$\nu \!= \log_{\gamma}\left\{\frac{\gamma\!-1}{\left[\prod\limits_{i,j\!=1
\atop i\neq j}^{n}\left(\frac{\gamma\!-1}{\gamma^{1\!-\overline{\eta}_3}
 \!-1}\right)\right]^{\frac{1}{n(n\!-1)}}}\!+1\right\}$,
and
$\overline{\eta}=\left[\prod\limits_{i,j=1 \atop i\neq j}^{n}
\left(\frac{\gamma-1}{\gamma^{\overline{\eta}_1-\overline{\eta}_2+1
-\overline{\eta}_3}-1}\right)\right]^{\frac{1}{n(n-1)}}$,
 $\overline{\eta}_1=\log_{\gamma}\left[\frac{\gamma-1}{\left(\frac{\gamma-1}
 {\gamma^{\eta_{\alpha_{i}}+\mu_{\alpha_{i}}}-1}\right)^{p}
 \left(\frac{\gamma-1}{\gamma^{\eta_{\alpha_{j}}
 +\mu_{\alpha_{j}}}-1}\right)^{q}}+1\right]$, $\overline{\eta}_2=
 \log_{\gamma}\left[\frac{\gamma-1}{\left(\frac{\gamma-1}
 {\gamma^{\mu_{\alpha_{i}}}-1}\right)^{p}\left(\frac{\gamma-1}
 {\gamma^{\mu_{\alpha_{j}}}-1}\right)^{q}}
 +1\right]$, $\overline{\eta}_3=\log_{\gamma}\left[\frac{\gamma-1}{\left(\frac{\gamma-1}
 {\gamma^{1-\nu_{\alpha_{i}}}-1}\right)^{p}\left(\frac{\gamma-1}
 {\gamma^{1-\nu_{\alpha_{j}}}-1}\right)^{q}}+1\right]$.

 (5) If we take $T$ as Dombi t-norms $T_{\gamma}^{\mathbf{D}}$
($\gamma \in (0, +\infty)$),
then
$\mathrm{PFIBM}_{T_{\gamma}^{\mathbf{D}}}^{p,q}(\alpha_{1}, \alpha_{2},\ldots ,\alpha_{n})=
 \Big\langle \frac{1}{1\!+\sqrt[\gamma]{\frac{1}{p\!+q}\left(\frac{1\!-\mu}{\mu}\right)^{\gamma}}},
 \frac{1}{1\!+\sqrt[\gamma]{\frac{1}{p\!+q}\left(\frac{1\!-\eta\!-\mu}{\eta\!+\mu}\right)^{\gamma}}}
\!-\frac{1}{1\!+\sqrt[\gamma]{\frac{1}{p\!+q}\left(\frac{1\!-\mu}{\mu}\right)^{\gamma}}},$
$
\frac{\sqrt[\gamma]{\frac{1}{p\!+q}\left(\frac{\nu}{1\!-\nu}\right)^{\gamma}}}{1\!+\sqrt[\gamma]
{\frac{1}{p\!+q}\left(\frac{\nu}{1\!-\nu}\right)^{\gamma}}}\Big\rangle,
$
where
\begin{align*}
 \mu= \frac{\sqrt[\gamma]{\frac{1}{n(n-1)}\sum\limits_{i,j=1 \atop i\neq j}^{n}{{\frac{1}{p\cdot \left(\frac{1\!-\mu_{\alpha_{i}}}{\mu_{\alpha_{i}}}\right)^{\gamma}\!+q\cdot \left(\frac{1\!-
 \mu_{\alpha_{j}}}{\mu_{\alpha_{j}}}\right)^{\gamma}}}}}}
 {1\!+\sqrt[\gamma]{\frac{1}{n(n\!-1)}\sum\limits_{i,j=1 \atop i\neq j}^{n}{{\frac{1}{p\cdot
 \left(\frac{1\!-\mu_{\alpha_{i}}}{\mu_{\alpha_{i}}}\right)^{\gamma}\!+q\cdot \left(\frac{1\!-
 \mu_{\alpha_{j}}}{\mu_{\alpha_{j}}}\right)^{\gamma}}}}}},
\end{align*}
\begin{align*}
\eta = & \frac{1}{1\!+\sqrt[\gamma]{\frac{1}{n(n\!-1)}\sum\limits_{i,j=1 \atop i\neq j}^{n}\left[\frac{1\!-\frac{1}{1\!+\sqrt[\gamma]{\overline{\eta}_1}}\!+\frac{1}{1\!+\sqrt[\gamma]
{\overline{\eta}_2}}\!-\frac{\sqrt[\gamma]
 {\overline{\eta}_3}}{1\!+\sqrt[\gamma]{\overline{\eta}_3}}}{\frac{1}{1\!+\sqrt[\gamma]
{\overline{\eta}_1}}\!-\frac{1}{1\!+\sqrt[\gamma]
{\overline{\eta}_2}}\!+\frac{\sqrt[\gamma]
 {\overline{\eta}_3}}{1\!+\sqrt[\gamma]{\overline{\eta}_3}}}\right]^{\gamma}}}\\
 &\!-\frac{1}{1\!+\sqrt[\gamma]{\frac{1}{n(n\!-1)}\sum\limits_{i,j=1 \atop i\neq j}^{n}
 \frac{1}{p\left(\frac{\nu_{\alpha_{i}}}{1\!-\nu_{\alpha_{i}}}\right){^{\gamma}}
 \!+q\left(\frac{\nu_{\alpha_{j}}}{1\!-\nu_{\alpha_{j}}}\right){^{\gamma}}}}},
\end{align*}
$\overline{\eta}_1=p\cdot\left(\frac{1\!-\eta_{\alpha_{i}}\!-\mu_{\alpha_{i}}}{\eta_{\alpha_{i}}\!+
\mu_{\alpha_{i}}}\right){^{\gamma}}\!+q\cdot \left(\frac{1\!-\eta_{\alpha_{j}}
\!-\mu_{\alpha_{j}}}{\eta_{\alpha_{j}}\!+\mu_{\alpha_{j}}}\right){^{\gamma}}$,
$\overline{\eta}_2=p\cdot \left(\frac{1\!-\mu_{\alpha_{i}}}{\mu_{\alpha_{i}}}\right){^{\gamma}}
\!+q\!\cdot \left(\frac{1\!-\mu_{\alpha_{j}}}{\mu_{\alpha_{j}}}\right){^{\gamma}}$,
$\overline{\eta}_3=p\cdot \left(\frac{\nu_{\alpha_{i}}}{1\!-\nu_{\alpha_{i}}}\right){^{\gamma}}
 \!+q\!\cdot \left(\frac{\nu_{\alpha_{j}}}{1\!-\nu_{\alpha_{j}}}\right){^{\gamma}}$,
and
\begin{align*}
\nu=\frac{1}{1\!+\sqrt[\gamma]{\frac{1}{n(n\!-1)}
\sum\limits_{i,j=1 \atop i\neq j}^{n}{\frac{1}{p\left(\frac{\nu_{\alpha_{i}}}
{1\!-\nu_{\alpha_{i}}}\right)^{\gamma}
\!+q\left(\frac{\nu_{\alpha_{j}}}{1\!-\nu_{\alpha_{j}}}\right)^{\gamma}}}}}.
\end{align*}
{In particular, if we take $\eta_{\alpha_i}=0$ for each $\alpha_i$, i.e., each $\alpha_i$ reduces to an
intuitionistic fuzzy number, then
\begin{equation}
\label{eq-Dom}
\begin{split}
&\quad \mathrm{PFIBM}_{T_{\gamma}^{\mathbf{D}}}^{p,q}
(\alpha_{1}, \alpha_{2},\ldots ,\alpha_{n})\\
& = \left\langle\frac{1}{1\!+\frac{1}{\sqrt[\gamma]{\frac{p\!+q}{n(n\!-1)}\sum\limits_{i,j=1 \atop i\neq j}^{n}{{\frac{1}{p\cdot\left(\frac{1\!-\mu_{\alpha_{i}}}{\mu_{\alpha_{i}}}\right)^{\gamma}\!+q\cdot
\left(\frac{1\!-\mu_{\alpha_{j}}}{\mu_{\alpha_{j}}}\right)^{\gamma}}}}}}}, ~ 0, \right.
\\
& \quad \quad \left. 1\!-\frac{1}{1\!+\frac{1}{\sqrt[\gamma]{\frac{p\!+q}{n(n\!-1)}
\sum\limits_{i,j=1 \atop i\neq j}^{n}{\frac{1}{p\cdot \left(\frac{\nu_{\alpha_{i}}}
{1\!-\nu_{\alpha_{i}}}\right)^{\gamma}\!+q\cdot \left(\frac{\nu_{\alpha_{j}}}
{1\!-\nu_{\alpha_{j}}}\right)^{\gamma}}}}}} \right\rangle,
\end{split}
\end{equation}
which is consistent with \cite[Theorem~1]{LLC2017}.
This means that \cite[Theorem~1]{LLC2017} is a direct corollary
of Theorem~\ref{PFIBM-Formula-Thm}.}

(6) If we take $T$ as Acz\'{e}l-Alsina t-norms $T_{\gamma}^{\mathbf{AA}}$
  ($\gamma \in (0, \!+\infty)$), then
$\mathrm{PFIBM}_{T_{\gamma}^{\mathbf{AA}}}^{p,q}(\alpha_{1}, \alpha_{2},\ldots,\alpha_{n})=
\Big\langle e^{\!-\sqrt[\gamma]{\frac{1}{p\!+q}\left(\!-\log\mu\right)^{\gamma}}},
e^{\!-\sqrt[\gamma]{\frac{1}{p\!+q}\left(\!-\log\left(\eta\!+\mu\right)\right)^{\gamma}}}
\!-e^{\!-\sqrt[\gamma]{\frac{1}{p\!+q}\left(\!-\log\mu\right)^{\gamma}}},$
$1\!-e^{\!-\sqrt[\gamma]{\frac{1}{p\!+q}\left(\!-\log\left(1\!-\nu\right)\right)^{\gamma}}}\Big\rangle$,
where
\begin{align*}
 \mu=& 1\!-e^{\!-\sqrt[\gamma]{\frac{1}{n(n\!-1)}\sum\limits_{i,j=1 \atop i\neq j}^{n}
 \left[\!-\log\left(1\!-\overline{\eta}_2\right)\right]^{\gamma}}},\\
\eta=& e^{\!-\sqrt[\gamma]{
\frac{1}{n(n\!-1)}\sum\limits_{i,j=1 \atop i\neq j}^{n}
\left[\!-\log\left(\overline{\eta}_1 \!- \overline{\eta}_2 \!+1\!-\overline{\eta}_3 \right)\right]
^{\gamma}}}
\\
 &\quad \!-e^{\!-\sqrt[\gamma]{\frac{1}{n(n\!-1)}\sum\limits_{i,j=1 \atop i\neq j}^{n}
 \left[\!-\log\left(1\!-\overline{\eta}_3\right)\right]^{\gamma}}},\\
\nu =& e^{\!-\sqrt[\gamma]{\frac{1}{n(n\!-1)}\sum\limits_{i,j=1 \atop i\neq j}^{n}
\left[\!-\log\left(1\!-\overline{\eta}_3\right)\right]^{\gamma}}},
\end{align*}
and
$\overline{\eta}_1=e^{-\sqrt[\gamma]{p\cdot \left(-\log\left(\eta_{\alpha_{i}}+\mu_{\alpha_{i}}\right)\right)^{\gamma}
+q\cdot \left(-\log\left(\eta_{\alpha_{j}}+\mu_{\alpha_{j}}\right)\right)^{\gamma}}}$, $\overline{\eta}_2=$
$e^{-\sqrt[\gamma]{p\cdot \left(-\log\mu_{\alpha_{i}}\right)^{\gamma}+q\cdot \left(-\log\mu_{\alpha_{j}}\right)^{\gamma}}}$,
$\overline{\eta}_3=e^{-\sqrt[\gamma]{p\cdot\left(-\log\left(1-\nu_{\alpha_{i}}\right)\right)^{\gamma}
+q\cdot\left(-\log\left(1-\nu_{\alpha_{j}}\right)\right)^{\gamma}}}$.

\begin{theorem}[{\textrm{\protect Monotonicity}}]
\label{PFIBM-Monotonicity-Thm}
Let $\alpha_{i}=\langle\mu_{\alpha_{i}},\eta_{\alpha_{i}},\nu_{\alpha_{i}}\rangle$ ($i=1,2,\ldots,n$), $\beta_{i}=\langle\mu_{\beta_{i}},\eta_{\beta_{i}},\nu_{\beta_{i}}\rangle$ ($i=1,2,\ldots,n$) be
two collections of PFNs such that $\mu_{\alpha_{i}} \leq \mu_{\beta_{i}}$,
$\eta_{\alpha_{i}} \leq \eta_{\beta_{i}}$, and $\nu_{\alpha_{i}} \geq \nu_{\beta_{i}}$.
Then, for $p$, $q>0$, we have
$$
\mathrm{PFIBM}_{T}^{p,q}(\alpha_{1}, \alpha_{2},\ldots, \alpha_{n})
\preccurlyeq_{_{\mathrm{W}}} \mathrm{PFIBM}_{T}^{p,q}(\beta_{1}, \beta_{2},\ldots, \beta_{n}).
$$
\end{theorem}
\begin{IEEEproof}
Let
$\alpha=\left\langle\mu_{\alpha},\eta_{\alpha},\nu_{\alpha}\right\rangle
=\mathrm{PFIBM}_{T}^{p,q}(\alpha_{1}, \alpha_{2}, \ldots,$ $\alpha_{n})$
and
$\beta=\left\langle\mu_{\beta},\eta_{\beta},\nu_{\beta}\right\rangle
=\mathrm{PFIBM}_{T}^{p,q}(\beta_{1}, \beta_{2},\ldots, \beta_{n}).$

Applying Theorem \ref{PFIBM-Formula-Thm}, we have
\begin{equation}
\label{eq-*}
\begin{split}
&\left\langle\mu_{\alpha},\eta_{\alpha},\nu_{\alpha}\right\rangle\\
=& \Bigg\langle\tau^{\!-1}\left(\frac{1}{p\!+q}\cdot\tau\left(\mu\right)\right),
 \tau^{\!-1}\left(\frac{1}{p\!+q}\cdot\tau\left(\eta\!+\mu\right)\right)\Bigg.\\
&\quad\Bigg.\!-\tau^{\!-1}\left(\frac{1}{p\!+q}\cdot\tau\left(\mu\right)\right),
\zeta^{\!-1}\left(\frac{1}{p\!+q}\cdot\zeta\left(\nu\right)\right)\Bigg\rangle,
\end{split}
\end{equation}
where
{\small \begin{align*}
\mu = \zeta^{\!-1}\Bigg(\frac{1}{n(n\!-1)}\sum_{i,j=1 \atop i\neq j}^{n}\zeta
\Big(\tau^{\!-1}\left(p\cdot\tau\left(\mu_{\alpha_{i}}\right)
\!+q\cdot\tau\left(\mu_{\alpha_{j}}\right)\right)\Big)\Bigg),
\end{align*}
\begin{align*}
\eta
=& \tau^{\!-1}\Bigg(\frac{1}{n(n\!-1)}\sum_{i,j=1 \atop i\neq j}^{n}\tau\Big(\tau^{\!-1}
\left(p\cdot\tau\left(\eta_{\alpha_{i}}\!+\mu_{\alpha_{i}}\right)\right.\Big.\Bigg.\\
&\quad \quad\left.\!+q\cdot\tau\left(\eta_{\alpha_{j}}\!+\mu_{\alpha_{j}}\right)\right)
\!-\tau^{\!-1}\left(p\cdot\tau\left(\mu_{\alpha_{i}}\right)\!+q\cdot
\tau\left(\mu_{\alpha_{j}}\right)\right)\\
& \Bigg.\Big. \quad \quad \!+\zeta^{\!-1}\left(p\cdot\zeta\left(\nu_{\alpha_{i}}\right)
\!+q\cdot\zeta\left(\nu_{\alpha_{j}}\right)\right)\Big)\Bigg)\\
&\!-\tau^{\!-1}\Bigg(\frac{1}{n(n\!-1)}\sum_{i,j=1 \atop i\neq j}^{n}
\tau\Big(\zeta^{\!-1}\left(p\cdot\zeta\left(\nu_{\alpha_{i}}\right)
\!+q\cdot\zeta\left(\nu_{\alpha_{j}}\right)\right)\Big)\Bigg),
\end{align*}
\begin{align*}
\nu = \tau^{\!-1}\Bigg(\frac{1}{n(n-1)}\sum_{i,j=1 \atop i\neq j}^{n}
\tau\Big(\zeta^{\!-1}\left(p\cdot\zeta\left(\nu_{\alpha_{i}}\right)
\!+q\cdot\zeta\left(\nu_{\alpha_{j}}\right)\right)\Big)\Bigg),
\end{align*}}
and
\begin{equation}
\label{eq-**}
\begin{split}
&\left\langle\mu_{\beta},\eta_{\beta},\nu_{\beta}\right\rangle\\
= &\Bigg\langle\tau^{-1}\left(\frac{1}{p\!+q}\cdot\tau\left(\mu^{\prime}\right)\right),
\tau^{\!-1}\left(\frac{1}{p\!+q}\cdot\tau\left(\eta^{\prime}\!+\mu^{\prime}\right)\right)\Bigg.\\
&\Bigg.\!-\tau^{\!-1}\left(\frac{1}{p\!+q}\cdot\tau\left(\mu^{\prime}\right)\right),
\zeta^{\!-1}\left(\frac{1}{p\!+q}\cdot\zeta\left(\nu^{\prime}\right)\right)\Bigg\rangle,
\end{split}
\end{equation}
where
{\small \begin{align*}
\mu^{\prime} = \zeta^{\!-1}\Bigg(\frac{1}{n(n\!-1)}\sum_{i,j=1 \atop i\neq j}^{n}
\zeta\Big(\tau^{\!-1}\left(p\cdot\tau\left(\mu_{\beta_{i}}\right)
\!+q\cdot\tau\left(\mu_{\beta_{j}}\right)\right)\Big)\Bigg),
\end{align*}
\begin{align*}
\eta^{\prime}
=& \tau^{\!-1}\Bigg(\frac{1}{n(n\!-1)}\sum_{i,j=1 \atop i\neq j}^{n}\tau\Big(\tau^{\!-1}
\left(p\cdot\tau\left(\eta_{\beta_{i}}\!+\mu_{\beta_{i}}\right)\right.\Big.\Bigg.\\
&\left.\quad\quad +q\cdot\tau\left(\eta_{\beta_{j}}+\mu_{\beta_{j}}\right)\right)
\!-\tau^{\!-1}\left(p\cdot\tau\left(\mu_{\beta_{i}}\right)
\!+q\cdot\tau\left(\mu_{\beta_{j}}\right)\right)\\
& \Bigg.\Big. \quad \quad \!+\zeta^{\!-1}\left(p\cdot\zeta\left(\nu_{\beta_{i}}\right)
\!+q\cdot\zeta\left(\nu_{\beta_{j}}\right)\right)\Big)\Bigg)\\
&\!-\tau^{\!-1}\Bigg(\frac{1}{n(n\!-1)}\sum_{i,j=1 \atop i\neq j}^{n}
\tau\Big(\zeta^{\!-1}\left(p\cdot\zeta\left(\nu_{\beta_{i}}\right)
\!+q\cdot\zeta\left(\nu_{\beta_{j}}\right)\right)\Big)\Bigg),
\end{align*}
\begin{align*}
\nu^{\prime} = \tau^{\!-1}\Bigg(\frac{1}{n(n\!-1)}\sum_{i,j=1 \atop i\neq j}^{n}
\tau\Big(\zeta^{\!-1}\left(p\cdot\zeta\left(\nu_{\beta_{i}}\right)
\!+q\cdot\zeta\left(\nu_{\beta_{j}}\right)\right)\Big)\Bigg).
\end{align*}}

Since $\tau$ is decreasing and $\zeta$ is increasing, from $\mu_{\alpha_{i}} \leq \mu_{\beta_{i}}$ and
$\nu_{\alpha_{i}} \geq \nu_{\beta_{i}}$, it follows that $\mu\leq \mu^{\prime}$ and
$\nu\geq \nu^{\prime}$, and thus $\mu_{\alpha} \leq \mu_{\beta}$ and
$\nu_{\alpha} \geq \nu_{\beta}$. To prove $\alpha \preccurlyeq_{_{\mathrm{W}}} \beta$,
let us consider the following three cases:

(1) If there exists $1 \leq i_{0} \leq n$ such that $\mu_{\alpha_{i_{0}}} < \mu_{\beta_{i_{0}}}$,
noting that $\tau$ is strictly decreasing, $\zeta$ is strictly increasing and $\tau(0) = +\infty$,
$\zeta(1) = +\infty$, from $\mu_{\alpha_{i}} \leq \mu_{\beta_{i}}$ and
$\nu_{\alpha_{i}} \geq \nu_{\beta_{i}}$ ($i = 1, 2,\ldots, n$),
it follows that
{\small \begin{align*}
\mu =& \zeta^{\!-1}\Bigg(\frac{1}{n(n\!-1)}\sum_{i,j=1 \atop i\neq j}^{n}\zeta
\Big(\tau^{\!-1}\left(p\cdot\tau\left(\mu_{\alpha_{i}}\right)\!+q\cdot
\tau\left(\mu_{\alpha_{j}}\right)\right)\Big)\Bigg)\\
=& \zeta^{\!-1}\Bigg(\frac{1}{n(n\!-1)}\sum_{i,j=1 \atop i\neq j, i\neq i_{0}}^{n}\zeta
\Big(\tau^{\!-1}\left(p\cdot\tau\left(\mu_{\alpha_{i}}\right)\!+q\cdot
\tau\left(\mu_{\alpha_{j}}\right)\right)\Big)\Bigg.\\
&\Bigg.\!+\frac{1}{n(n\!-1)}\sum_{j=1 \atop j\neq i_{0}}^{n}\zeta
\Big(\tau^{\!-1}\left(p\cdot\tau\left(\mu_{\alpha_{i_{0}}}\right)\!+q\cdot
\tau\left(\mu_{\alpha_{j}}\right)\right)\Big)\Bigg)\\
<& \zeta^{\!-1}\Bigg(\frac{1}{n(n\!-1)}\sum_{i,j=1 \atop i\neq j, i\neq i_{0}}^{n}\zeta
\Big(\tau^{\!-1}\left(p\cdot\tau\left(\mu_{\beta_{i}}\right)+q\cdot
\tau\left(\mu_{\beta_{j}}\right)\right)\Big)\Bigg.\\
&\Bigg.+\frac{1}{n(n\!-1)}\sum_{j=1 \atop j\neq i_{0}}^{n}\zeta
\Big(\tau^{\!-1}\left(p\cdot\tau\left(\mu_{\beta_{i_{0}}}\right)+q\cdot
\tau\left(\mu_{\beta_{j}}\right)\right)\Big)\Bigg)\\
 =& \zeta^{\!-1}\Bigg(\frac{1}{n(n\!-1)}\sum_{i,j=1 \atop i\neq j}^{n}\zeta
 \Big(\tau^{\!-1}\left(p\cdot\tau\left(\mu_{\beta_{i}}\right)\!+q\cdot
 \tau\left(\mu_{\beta_{j}}\right)\right)\Big)\Bigg)\\
 =& \mu^{\prime},
\end{align*}}
and
{\small
\begin{align*}
\nu=& \tau^{\!-1}\Bigg(\frac{1}{n(n\!-1)} \sum_{i,j=1 \atop i\neq j}^{n}
\tau\Big(\zeta^{\!-1}\left(p\cdot\zeta\left(\nu_{\alpha_{i}}\right)
\!+q\cdot\zeta\left(\nu_{\alpha_{j}}\right)\right)\Big)\Bigg)\\
\geq & \tau^{\!-1}\Bigg(\frac{1}{n(n\!-1)} \sum_{i,j=1 \atop i\neq j}^{n}
\tau\Big(\zeta^{\!-1}\left(p\cdot\zeta\left(\nu_{\beta_{i}}\right)
\!+q\cdot\zeta\left(\nu_{\beta_{j}}\right)\right)\Big)\Bigg)\\
 =&\nu^{\prime}.
\end{align*}}
Together with formulas~\eqref{eq-*} and \eqref{eq-**}, it holds that
$\mu_{\alpha}=\tau^{-1}\left(\frac{1}{p+q}\cdot\tau\left(\mu\right)\right)
<\tau^{-1}\left(\frac{1}{p+q}\cdot\tau\left(\mu^{\prime}\right)\right)=\mu_{\beta}$
and $\nu_{\alpha} = \zeta^{-1}\left(\frac{1}{p+q}\cdot\zeta\left(\nu\right)\right)\geq \zeta^{-1}\left(\frac{1}{p+q}\cdot\zeta\left({\nu^{\prime}}\right)\right)
= \nu_{\beta}$, implying that $S(\alpha)=\mu_{\alpha}-\nu_{\alpha}
 < \mu_{\beta}-\nu_{\beta} = S(\beta)$, and thus $\alpha \preccurlyeq_{_{\mathrm{W}}} \beta$.

(2) If there exists $1 \leq i_{0} \leq n$ such that $\nu_{\alpha_{i_{0}}}
> \nu_{\beta_{i_{0}}}$, similarly to the proof of (1), it can be verified that
$\mu_{\alpha} \leq \mu_{\beta}$ and $\nu_{\alpha}>\nu^{\beta}$, implying that
$S(\alpha)= \mu_{\alpha}-\nu_{\alpha} < \mu_{\beta}-\nu_{\beta} = S(\beta)$, and thus
$\alpha \preccurlyeq_{_{\mathrm{W}}} \beta$.

(3) If $\mu_{\alpha_{i}} = \mu_{\beta_{i}}$ and $\nu_{\alpha_{i}} =
\nu_{\beta_{i}}$ for all $1 \leq i \leq n$, then $\mu=\mu^{\prime}$,
$\mu_{\alpha} = \mu_{\beta}$, and $\nu_{\alpha} = \nu_{\beta}$, and thus
$S(\alpha) = S(\beta)$ and $H_{1}(\alpha) = H_{1}(\beta)$. Since $\tau$
is strictly decreasing and $\zeta$ is strictly increasing, from
$\eta_{\alpha_{i}} \leq \eta_{\beta_{i}}$ ($j = 1, 2,\ldots, n$), we have
{\small
\begin{align*}
\eta=& \tau^{\!-1}\Bigg(\frac{1}{n(n-1)}\sum_{i,j=1 \atop i\neq j}^{n}
\tau\Big(\tau^{\!-1}\left(p\cdot\tau\left(\eta_{\alpha_{i}}
\!+\mu_{\alpha_{i}}\right)\right.\Big.\Bigg.\\
&\left.\quad \quad \!+q\cdot\tau\left(\eta_{\alpha_{j}}+\mu_{\alpha_{j}}\right)\right)
\!-\tau^{\!-1}\left(p\cdot\tau\left(\mu_{\alpha_{i}}\right)
\!+q\cdot\tau\left(\mu_{\alpha_{j}}\right)\right)\\
& \Bigg.\Big.\quad \quad \!+\zeta^{\!-1}\left(p\cdot\zeta\left(\nu_{\alpha_{i}}\right)
\!+q\cdot\zeta\left(\nu_{\alpha_{j}}\right)\right)\Big)\Bigg)
\end{align*}
\begin{align*}
&\!-\tau^{\!-1}\Bigg(\frac{1}{n(n\!-1)}\sum_{i,j=1 \atop i\neq j}^{n}
\tau\Big(\zeta^{\!-1}\left(p\cdot\zeta\left(\nu_{\alpha_{i}}\right)
+q\cdot\zeta\left(\nu_{\alpha_{j}}\right)\right)\Big)\Bigg)\\
\leq & \tau^{\!-1}\Bigg(\frac{1}{n(n\!-1)}\sum_{i,j=1 \atop i\neq j}^{n}\tau\Big(\tau^{\!-1}
\left(p\cdot\tau\left(\eta_{\beta_{i}}+\mu_{\beta_{i}}\right)\right.\Big.\Bigg.\\
&\left.\quad \quad+q\cdot\tau\left(\eta_{\beta_{j}}\!+\mu_{\beta_{j}}\right)\right)
\!-\tau^{\!-1}\left(p\cdot\tau\left(\mu_{\beta_{i}}\right)
\!+q\cdot\tau\left(\mu_{\beta_{j}}\right)\right)
\\
& \Bigg.\Big.\quad \quad \quad
\!+\zeta^{\!-1}\left(p\cdot\zeta\left(\nu_{\beta_{i}}\right)\!+q\cdot
\zeta\left(\nu_{\beta_{j}}\right)\right)\Big)\Bigg)\\
&\!-\tau^{\!-1}\Bigg(\frac{1}{n(n\!-1)}\sum_{i,j=1 \atop i\neq j}^{n}
\tau\Big(\zeta^{\!-1}\left(p\cdot\zeta\left(\nu_{\beta_{i}}\right)
\!+q\cdot\zeta\left(\nu_{\beta_{j}}\right)\right)\Big)\Bigg)\\
=& \eta^{\prime}.
\end{align*}}
This, together with $\mu=\mu^{\prime}$, implies that
\begin{align*}
\eta_{\alpha}=&\tau^{\!-1}\left(\frac{1}{p\!+q}\cdot\tau(\eta+\mu)\right)
\!-\tau^{\!-1}\left(\frac{1}{p\!+q}\cdot\tau(\mu)\right)\\
\leq& \tau^{\!-1}\left(\frac{1}{p\!+q}\cdot\tau(\eta^{\prime}\!+\mu^{\prime})\right)
\!-\tau^{\!-1}\left(\frac{1}{p\!+q}\cdot\tau(\mu^{\prime})\right)\\
=& \eta_{\beta}.
\end{align*}
Then, $H_{2}(\alpha) = \mu_{\alpha}+\eta_{\alpha}+\nu_{\alpha}\leq
\mu_{\beta}+\eta_{\beta}+\nu_{\beta}=H_{2}(\beta)$, and thus $\alpha \preccurlyeq_{_{\mathrm{W}}} \beta$.
\end{IEEEproof}

\begin{theorem}[{\textrm{\protect Idempotency}}]
\label{PFIBM-Idempotency-Thm}
If all $\alpha_{i}$ ($i=1,2,\ldots,n$) are equal, i.e., $\alpha_{i}=\alpha$ for all $i=1,2,\ldots,n$,
then
$$
\mathrm{PFIBM}_{T}^{p,q}(\alpha_{1}, \alpha_{2},\ldots, \alpha_{n})=\alpha.
$$
\end{theorem}
\begin{IEEEproof}
By Theorem~\ref{Operation-Properties-Thm} (v), (vi), (ix), (x), it holds that
{\small \begin{align*}
&\mathrm{PFIBM}_{T}^{p,q}(\alpha_{1}, \alpha_{2},\ldots,\alpha_{n})\\
 =& \left[\frac {1}{n(n\!-1)}\cdot \bigoplus\limits_{i,j=1 \atop i\neq j}^{n}
 \left(\alpha^{p}\otimes\alpha^{q}\right)\right]^{\frac{1}{p\!+q}}\\
 =& \left[\frac {1}{n(n\!-1)}\cdot \bigoplus
\limits_{i,j=1 \atop i\neq j}^{n}\alpha^{p+q}\right]^{\frac{1}{p\!+q}} 
\quad \text{(by Theorem~\ref{Operation-Properties-Thm} (vi))}\\
 =& \left[\frac{1}{n(n\!-1)}\cdot\left(\left(n(n\!-1)\right)\cdot\alpha^{p+q}\right)\right]^{\frac{1}{p+q}}
 \quad \text{(by Theorem~\ref{Operation-Properties-Thm} (v))}\\
 =& \left(\alpha^{p\!+q}\right)^{\frac{1}{p\!+q}}  \quad \text{(by Theorem~\ref{Operation-Properties-Thm} (ix))}\\
 =& \alpha  \quad \text{(by Theorem~\ref{Operation-Properties-Thm} (x))}.
\end{align*}
}
\end{IEEEproof}

\begin{theorem}[{\textrm{\protect Boundedness}}]
\label{PFIBM-Boundedness-Thm}
Let $\alpha_{i}=\left\langle\mu_{\alpha_{i}},
\eta_{\alpha_{i}},\nu_{\alpha_{i}}\right\rangle$
($i=1, \ldots,n$) be a collection of PFNs, then
$\langle \min\limits_{i}\left\{\mu_{\alpha_{i}}\right\},$
$\min\limits_{i}\left\{\eta_{\alpha_{i}}\right\},
\max\limits_{i}\left\{\nu_{\alpha_{i}}\right\}\rangle
\preccurlyeq_{_{\mathrm{W}}} \mathrm{PFIBM}_{T}^{p,q}(\alpha_{1}, \alpha_{2},
\ldots ,\alpha_{n})
\preccurlyeq_{_{\mathrm{W}}}\langle \max\limits_{i}\left\{\mu_{\alpha_{i}}\right\},
 1\!- (\max\limits_{i}\left\{\mu_{\alpha_{i}}\right\}\!+\min\limits_{i}
 \left\{\nu_{\alpha_{i}}\right\}), \min\limits_{i}\left\{\nu_{\alpha_{i}}\right\} \rangle$.
\end{theorem}
\begin{IEEEproof}
Let
$\alpha^{\!-} \!= \langle \mu_{1}, \eta_{1}, \nu_{1} \rangle \!=
\langle \min\limits_{i}\left\{\mu_{\alpha_{i}}\right\},
\min\limits_{i}\left\{\eta_{\alpha_{i}}\right\},$ $
\max\limits_{i}\left\{\nu_{\alpha_{i}}\right\} \rangle$,
$\alpha^{\!+} \!=\langle \mu_{2}, \eta_{2}, \nu_{2}\rangle\!=
\langle \max\limits_{i}\left\{\mu_{\alpha_{i}}\right\},
1\!- (\max\limits_{i}\left\{\mu_{\alpha_{i}}\right\}\!+
\min\limits_{i}\left\{\nu_{\alpha_{i}}\right\}),
\min\limits_{i}\left\{\nu_{\alpha_{i}}\right\}\rangle$,
and
$\alpha=\langle\mu_{\alpha},\eta_{\alpha},$ $\nu_{\alpha}\rangle
=\mathrm{PFIBM}_{T}^{p,q}(\alpha_{1}, \alpha_{2},\ldots, \alpha_{n}).$

From $\mu_{1} \leq \mu_{\alpha_{i}}$, $\eta_{1} \leq \eta_{\alpha_{i}}$, and
$\nu_{1} \geq \nu_{\alpha_{i}}$ for all $1 \leq i \leq n$, by Theorems~\ref{PFIBM-Monotonicity-Thm} and
\ref{PFIBM-Idempotency-Thm}, it follows that
$\alpha^{\!-}=\mathrm{PFIBM}_{T}^{p,q}(\alpha^{\!-}, \alpha^{\!-},\ldots , \alpha^{\!-})
 \preccurlyeq_{_{\mathrm{W}}}\mathrm{PFIBM}_{T}^{p,q}(\alpha_{1}, \alpha_{2},\ldots, \alpha_{n})=\alpha.
$

Clearly, $\mu_{\alpha_{i}} \leq \mu_{2}$ and $\nu_{\alpha_{i}} \geq \nu_{2}$
holds for all $1 \leq i \leq n$. To prove $\alpha \leq \alpha^{+}$, we consider the following three cases:

(1) If there exists $1 \leq i_{0} \leq n$ such that $\mu_{\alpha_{i_{0}}} < \mu_{2}$, since $\tau$
is strictly decreasing with $\tau(0) = +\infty$, $\zeta$ is strictly increasing with $\zeta(1) = +\infty$,
and $\mu_{\alpha_{i}} \leq \mu_{2}$ ($j=1,2,\ldots, n$), by Theorem \ref{PFIBM-Formula-Thm}, we have
$\mu_{\alpha}=\tau^{\!-1}\big(\frac{1}{p+q}\cdot\tau (\zeta^{\!-1} (\frac{1}{n(n\!-1)}
\sum\limits_{i,j=1 \atop i\neq j}^{n}\zeta(\tau^{-1}(p\cdot\tau (\mu_{\alpha_{i}} )
\!+q\cdot\tau (\mu_{\alpha_{j}} ) ) ) ) ) \big)< \tau^{\!-1}\big(\frac{1}{p\!+q}\cdot
\tau (\zeta^{\!-1} (\frac{1}{n(n\!-1)}\sum\limits_{i,j=1 \atop i\neq j}^{n}\zeta (\tau^{\!-1}
 (p\cdot\tau (\mu_{2} )\!+q\cdot\tau (\mu_{2} ) ) ) ) ) \big)=\mu_{2}$,
and from $\nu_{\alpha_{i}} \geq \nu_{2}$, we have
$\nu_{\alpha} =\zeta^{\!-1}\big(\frac{1}{p\!+q}\cdot\zeta(\tau^{\!-1}
(\frac{1}{n(n\!-1)} \sum\limits_{i,j=1 \atop i\neq j}^{n}
\tau(\zeta^{\!-1}(p\cdot\zeta(\nu_{\alpha_{i}})\!+q\cdot\zeta(\nu_{\alpha_{j}}
)))))\big) \geq \zeta^{\!-1}\big(\frac{1}{p\!+q}\cdot\zeta(\tau^{\!-1}(\frac{1}{n(n\!-1)}
\sum\limits_{i,j=1 \atop i\neq j}^{n}\tau(\zeta^{\!-1}(p\cdot\zeta
(\nu_{2})\!+q\cdot\zeta(\nu_{2})))))\big)=\nu_{2}$,
implying that
$S(\alpha) = \mu_{\alpha}\!-\nu_{\alpha} \!< 
\mu_{2}\!-\nu_{2} = S(\alpha^{\!+})$,
and thus $\alpha \prec_{_{\mathrm{W}}} \alpha^{\!+}$.

(2) If there exists $1 \leq i_{0} \leq n$ such that $\nu_{\alpha_{i_{0}}} > \nu_{2}$,
similarly to the above proof, it can be verified that
$S(\alpha) = \mu_{\alpha}-\nu_{\alpha} < \mu_{2}-\nu_{2} = S(\alpha^{+})$,
and thus $\alpha \prec_{_{\mathrm{W}}} \alpha^{+}$.

(3) If $\mu_{\alpha_{i}} = \mu_{2}$ and $\nu_{\alpha_{i}} = \nu_{2}$
holds for all $1 \leq i \leq n$, from $\mu_{\alpha_{i}}+\eta_{\alpha_{i}}
+\nu_{\alpha_{i}} \leq 1$,  it follows that $\eta_{\alpha_{i}} \leq
\max\limits_{i}\{\eta_{\alpha_{i}}\} \leq 1
-(\mu_{\alpha_{i}}+\nu_{\alpha_{i}}) = \eta_{2}$.
Then by Theorems~\ref{PFIBM-Monotonicity-Thm} and
\ref{PFIBM-Idempotency-Thm}, we have $\alpha \preccurlyeq_{_{\mathrm{W}}} \alpha^{+}$.
\end{IEEEproof}

\begin{theorem}[{\textrm{\protect Commutativity}}]
\label{PFIBM-Commutativity-Thm}
Let $\alpha_{i}$ ($i=1,2,\ldots,n$) be a collection of PFNs. Then
$\mathrm{PFIBM}_{T}^{p,q}(\alpha_{1}, \alpha_{2},\ldots ,\alpha_{n})=
\mathrm{PFIBM}_{T}^{p,q}(\dot{\alpha}_{1}, \dot{\alpha}_{2}, \ldots,\dot{\alpha}_{n})$,
where $(\dot{\alpha}_{1}, \dot{\alpha}_{2},\ldots, \dot{\alpha}_{n})$ is
any permutation of $(\alpha_{1}, \alpha_{2},\ldots ,\alpha_{n})$.
\end{theorem}

\begin{IEEEproof}
It follows directly from the commutativity of the operations $\oplus$ and $\otimes$.
\end{IEEEproof}

\begin{remark}
By Definition~\ref{Wu-Order-Def}, Theorems~\ref{PFIBM-Monotonicity-Thm}, \ref{PFIBM-Idempotency-Thm}, and \ref{PFIBM-Commutativity-Thm},
and formulas~\eqref{eq-T-formula}, \eqref{eq-Tp}, and \eqref{eq-Dom},
it follows that \cite[Theorems~2, 3, and 6]{Ni2016}, \cite[Idempotency, Monotonicity, and Commutativity]{XC2012,XY2011},
and \cite[Theorems 2 and 3]{{LLC2017}} hold trivially.
\end{remark}

\section{PFIWBM and PFINWM operators}\label{Sec-5}

The PFIBM operator can capture the interrelationship between each pair of input parameters;
however, it has a shortcoming that it does not consider the essentiality of parameters
(i.e., the weights of criteria or experts for MCDM or MCGDM problems).
In other words, the PFIBM is only applicable for equal-weighted MCDM problems.
To overcome this shortcoming, we shall introduce the picture fuzzy interactional weighted
Bonferroni mean (PFIWBM) and picture fuzzy interactional normalized weighted Bonferroni mean
(PFINWM) operators in this section.

\begin{definition}
\label{PFIWBM-Def}
Let $\alpha_{i}$ ($i=1,2,\ldots,n$) be a collection of PFNs,
$\omega=(\omega_{1}, \omega_{2},\ldots, \omega_{n})^{\mathrm{T}}$ be the weight vector
of $\alpha_{i}$ ($i=1,2,\ldots,n$) such that $\omega_{i} \in (0, 1]$ and $\sum_{i=1}^{n}
\omega_{i}=1$, and $T$ be a strict t-norm. For any $p$, $q> 0$, define the \textit{picture
fuzzy interactional weighted Bonferroni mean} (PFIWBM) induced by $T$ as follows:
{\small \begin{equation}\label{19}
\begin{aligned}
&\mathrm{PFIWBM}_{T}^{p,q}(\alpha_{1}, \alpha_{2},\ldots ,\alpha_{n})\\
= &\left[\frac {1}{n(n\!-1)}\cdot \bigoplus\limits_{i,j=1 \atop i\neq j}^{n}
\Big(\left(\omega_{i}\cdot\alpha_{i}\right)^{p}\otimes
\left(\omega_{j}\cdot\alpha_{j}\right)^{q}\Big)\right]^{\frac{1}{p\!+q}}.
\end{aligned}
\end{equation}}
\end{definition}

\begin{definition}
\label{PFINWBM-Def}
Let $\alpha_{i}$ ($i=1,2,\ldots,n$) be a collection of PFNs,
$\omega=(\omega_{1}, \omega_{2},\ldots, \omega_{n})^{\mathrm{T}}$ be the weight vector
of $\alpha_{i}$ ($i=1,2,\ldots,n$) such that $\omega_{i} \in (0, 1]$ and $\sum_{i=1}^{n}\omega_{i}=1$,
and $T$ be a strict t-norm. For any $p$, $q > 0$, define the \textit{picture fuzzy interactional
normalized weighted Bonferroni mean} (PFINWBM) induced by $T$ as follows:
{\small \begin{equation}\label{24}
\mathrm{PFINWBM}_{T}^{p,q}(\alpha_{1}, \ldots ,\alpha_{n})
=\left[\bigoplus\limits_{i,j\!=1 \atop i\neq j}^{n}
\frac {\omega_{i}\omega_{j}}{1\!-\omega_{i}}\cdot\left(\alpha_{i}^{p}\otimes
\alpha_{j}^{q}\right)\right]^{\frac{1}{p\!+q}}.
\end{equation}}
\end{definition}

Similarly to the proof of Theorem~\ref{PFIBM-Formula-Thm},
the following results can be obtained.

\begin{theorem}
\label{PFIWBM-Formula-Thm}
Let $\alpha_{i}=\left\langle\mu_{\alpha_{i}},\eta_{\alpha_{i}},\nu_{\alpha_{i}}\right\rangle$
($i=1,2,\ldots,n$) be a collection of PFNs, $\omega=(\omega_{1}, \omega_{2},\ldots, \omega_{n})^{\mathrm{T}}$
be the weight vector of $\alpha_{i}$ ($i=1,2,\ldots,n$) such that $\omega_{i} \in (0, 1]$ and
$\sum_{i=1}^{n}\omega_{i}=1$, and $T$ be a strict t-norm with an additive generator $\tau$.
Then, for $p$, $q > 0$, the aggregated value by using the PFIWBM induced by $T$ is also an PFN, and
\begin{equation}\label{20}
\begin{aligned}
&\mathrm{PFIWBM}_{T}^{p,q}(\alpha_{1}, \alpha_{2},\ldots,\alpha_{n})\\
=&\left\langle\tau^{\!-1}\left(\frac{1}{p\!+q}\cdot\tau\left(\bar{\mu}\right)\right)\right.,\\
& \quad \tau^{\!-1}\left(\frac{1}{p\!+q}\cdot\tau(\bar{\eta}\!+\bar{\mu})\right)\!-\tau^{\!-1}
\left(\frac{1}{p\!+q}\cdot\tau
(\bar{\mu})\right),\\
& \quad \left.\zeta^{\!-1}\left(\frac{1}{p\!+q}\cdot\zeta(\bar{\nu})\right)\right\rangle,
\end{aligned}
\end{equation}
where
$\bar{\mu} = \zeta^{\!-1}\bigg(\frac{1}{n(n\!-1)}\cdot\sum\limits_{i,j=1 \atop i\neq j}^{n}\zeta\big(\tau^{\!-1}(p\cdot\tau(\zeta^{\!-1}(\omega_{i}\cdot\zeta(\mu_{\alpha_{i}})))
\!+q\cdot\tau(\zeta^{\!-1}(\omega_{j}\cdot\zeta(\mu_{\alpha_{j}}))))\big)\bigg)$,
$\bar{\eta}
=\tau^{\!-1}\bigg(\frac{1}{n(n\!-1)}\cdot\sum\limits_{i,j=1 \atop i\neq j}^{n}\tau\bigg(\tau^{\!-1}\Big(p\cdot\tau(\tau^{\!-1}(\omega_{i}\cdot\tau(\eta_{\alpha_{i}}
\!+\nu_{\alpha_{i}}))\!-\tau^{\!-1}(\omega_{i}\cdot\tau(\nu_{\alpha_{i}}))
\!+\zeta^{\!-1}(\omega_{i}\cdot\zeta(\mu_{\alpha_{i}})))
\!+q\cdot\tau (\tau^{\!-1}(\omega_{j}\cdot\tau(\eta_{\alpha_{j}}\!+\nu_{\alpha_{j}}))
\!-\tau^{\!-1}(\omega_{j}\cdot\tau({\nu_{\alpha_{j}}}))
\!+\zeta^{\!-1}(\omega_{j}\cdot\zeta(\mu_{\alpha_{j}})))\Big)
\!-\tau^{\!-1}
\Big(p\cdot\tau(\zeta^{\!-1}(\omega_{i}\cdot\zeta
(\mu_{\alpha_{i}})))\!+q\cdot\tau(\zeta^{\!-1}
(\omega_{j}\cdot\zeta(\mu_{\alpha_{j}})))\Big)
\!+\zeta^{\!-1}\Big(p\cdot\zeta(\tau^{\!-1}(\omega_{i}\cdot\tau(\nu_{\alpha_{i}})))
\!+q\cdot\zeta(\tau^{\!-1}(\omega_{j}\cdot\tau(\nu_{\alpha_{j}})))\Big)\bigg)\bigg)
\!-\tau^{\!-1}\bigg(\frac{1}{n(n\!-1)}\cdot\sum\limits_{i,j=1 \atop i\neq j}^{n}\tau\bigg(\zeta^{\!-1}\Big(p\cdot\zeta(\tau^{\!-1}(\omega_{i}\cdot\tau(\nu_{\alpha_{i}})))
\!+q\cdot\zeta(\tau^{\!-1}(\omega_{j}\cdot\tau(\nu_{\alpha_{j}})))\Big)\bigg)\bigg)$,
$\bar{\nu}= \tau^{\!-1}\bigg(\frac{1}{n(n\!-1)}\cdot\sum\limits_{i,j=1 \atop i\neq j}^{n}
\tau\Big(\zeta^{\!-1}(p\cdot\zeta(\tau^{\!-1}(\omega_{i}
\cdot\tau(\nu_{\alpha_{i}}))))\!+q\cdot\zeta(\tau^{\!-1}
(\omega_{j}\cdot\tau(\nu_{\alpha_{j}}))))\Big)\bigg)$,
and $\zeta(x)=\tau(1\!-x)$.
\end{theorem}

\begin{theorem}
\label{PFINWBM-Formula-Thm}
Let $\alpha_{i}=\langle\mu_{\alpha_{i}},\eta_{\alpha_{i}},\nu_{\alpha_{i}}\rangle$
($i=1,2,\ldots,n$) be a collection of PFNs, $\omega=(\omega_{1}, \omega_{2},\ldots, \omega_{n})^{\mathrm{T}}$
be the weight vector of $\alpha_{i}$ ($i=1,2,\ldots,n$) such that $\omega_{i} \in (0, 1]$ and
$\sum_{i=1}^{n}\omega_{i}=1$, and $T$ be a strict t-norm with an additive generator $\tau$.
Then, for $p$, $q > 0$, the aggregated value by using the PFINWBM operator induced by $T$ is also an PFN, and
\begin{equation}\label{25}
\begin{aligned}
&\mathrm{PFINWBM}_{T}^{p,q}(\alpha_{1}, \alpha_{2},\ldots,\alpha_{n})\\
=&\left\langle\tau^{\!-1}\left(\frac{1}{p\!+q}\cdot\tau(\hat{\mu})\right)\right.,\\
& \tau^{\!-1}\left(\frac{1}{p\!+q}\cdot\tau(\hat{\eta}\!+\hat{\mu})\right)
\!-\tau^{\!-1}\left(\frac{1}{p\!+q}\cdot\tau(\hat{\mu})\right),\\
& \left.\zeta^{\!-1}\left(\frac{1}{p\!+q}\cdot\zeta(\hat{\nu})\right)\right\rangle,
\end{aligned}
\end{equation}
where
$\hat{\mu}= \zeta^{\!-1}\bigg(\sum\limits_{i,j=1 \atop i\neq j}^{n}
\frac{\omega_{i}\omega_{j}}{1\!-\omega_{i}}\cdot\zeta(\tau^{\!-1}
(p\cdot\tau(\mu_{\alpha_{i}})\!+q \cdot\tau(\mu_{\alpha_{j}})))\bigg)$,
$\hat{\eta}= \tau^{\!-1}\bigg(\sum\limits_{i,j=1 \atop i\neq j}^{n}
\frac{\omega_{i}\omega_{j}}{1\!-\omega_{i}}\cdot \tau
(\tau^{\!-1}(p\cdot\tau(\eta_{\alpha_{i}}\!+\mu_{\alpha_{i}})
\!+q\cdot\tau(\eta_{\alpha_{j}}\!+\mu_{\alpha_{j}}))
\!-\tau^{\!-1}(p\cdot\tau(\mu_{\alpha_{i}})\!+q\cdot\tau(\mu_{\alpha_{j}}))
\!+\zeta^{\!-1}(p\cdot\zeta(\nu_{\alpha_{i}})\!+q\cdot\zeta(\nu_{\alpha_{j}})))\bigg)
\!-\tau^{\!-1}\bigg(\sum\limits_{i,j=1 \atop i\neq j}^{n}
\frac{\omega_{i}\omega_{j}}{1\!-\omega_{i}}
\cdot\tau (\zeta^{\!-1} (p\cdot\zeta(\nu_{\alpha_{i}})
\!+q\cdot\zeta(\nu_{\alpha_{j}})))\bigg)$,
$\hat{\nu} = \tau^{\!-1}\Big(\sum\limits_{i,j=1 \atop i\neq j}^{n}
\frac{\omega_{i}\omega_{j}}{1\!-\omega_{i}}
\cdot\tau\Big(\zeta^{\!-1}(p\cdot\zeta(\nu_{\alpha_{i}})
\!+q\cdot\zeta(\nu_{\alpha_{j}}))\Big)\Big)$,
and $\zeta(x)=\tau(1\!-x)$.
\end{theorem}

\begin{theorem}[{\textrm{\protect Commutativity}}]
\label{PFIWBM-Commutativity-Thm}
Let $\alpha_{i}$ ($i=1,2,\ldots,n$) be a collection of PFNs. Then
\begin{align*}
\mathrm{PFIWBM}_{T}^{p,q}(\alpha_{1}, \ldots ,\alpha_{n})\!=
\mathrm{PFIWBM}_{T}^{p,q}(\dot{\alpha}_{1}, \ldots,\dot{\alpha}_{n}),
\end{align*}
and
\begin{align*}
\mathrm{PFINWBM}_{T}^{p,q}(\alpha_{1},\ldots ,\alpha_{n})\!=
\mathrm{PFINWBM}_{T}^{p,q}(\dot{\alpha}_{1}, \ldots,\dot{\alpha}_{n}),
\end{align*}
where $(\dot{\alpha}_{1}, \ldots, \dot{\alpha}_{n})$ is
any permutation of $(\alpha_{1}, \ldots ,\alpha_{n})$.
\end{theorem}

\begin{theorem}[{\textrm{\protect Monotonicity}}]
\label{PFIWBM-Monotonicity-Thm}
Let $\alpha_{i}=\langle\mu_{\alpha_{i}},\eta_{\alpha_{i}},\nu_{\alpha_{i}}\rangle$ ($i=1,2,\ldots,n$), $\beta_{i}=\langle\mu_{\beta_{i}},\eta_{\beta_{i}},\nu_{\beta_{i}}\rangle$ ($i=1,2,\ldots,n$) be
two collections of PFNs, $\omega=(\omega_{1}, \omega_{2},\ldots, \omega_{n})^{\mathrm{T}}$
be the weight vector such that $\omega_{i} \in (0, 1]$ and $\sum_{i=1}^{n}\omega_{i}=1$,
and $T$ be a strict t-norm with an additive generator $\tau$. If $\mu_{\alpha_{i}} \leq \mu_{\beta_{i}}$,
$\eta_{\alpha_{i}} \leq \eta_{\beta_{i}}$, and $\nu_{\alpha_{i}} \geq \nu_{\beta_{i}}$, then, for $p$, $q>0$, we have
$$
\mathrm{PFIWBM}_{T}^{p,q}(\alpha_{1}, \ldots, \alpha_{n})
\preccurlyeq_{_{\mathrm{W}}} \mathrm{PFIWBM}_{T}^{p,q}(\beta_{1}, \ldots, \beta_{n}),
$$
and
$$
\mathrm{PFINWBM}_{T}^{p,q}(\alpha_{1}, \ldots, \alpha_{n})
\preccurlyeq_{_{\mathrm{W}}} \mathrm{PFINWBM}_{T}^{p,q}(\beta_{1}, \ldots, \beta_{n}).
$$
\end{theorem}

\begin{IEEEproof}
Similarly to the proof of Theorem~\ref{PFIBM-Monotonicity-Thm},
by using Theorems~\ref{PFIWBM-Formula-Thm} and \ref{PFINWBM-Formula-Thm},
it can be verified that this theorem holds. 
\end{IEEEproof}

\begin{theorem}[{\textrm{\protect Idempotency}}]
\label{PFINWBM-Idempotency-Thm}
Let $\alpha_{i}$ ($i=1,2,\ldots,n$) be a collection of PFNs,
$\omega=(\omega_{1}, \omega_{2},\ldots, \omega_{n})^{\mathrm{T}}$ be the weight vector
of $\alpha_{i}$ ($i=1,2,\ldots,n$) such that $\omega_{i} \in (0, 1]$ and $\sum_{i=1}^{n}\omega_{i}=1$,
and $T$ be a strict t-norm. If all $\alpha_{i}$ ($i=1,2,\ldots,n$) are equal, i.e.,
$\alpha_{i}=\alpha$ for all $i=1, 2, \ldots, n$,
then, for any $p$, $q > 0$, we have
$
\mathrm{PFINWBM}_{T}^{p,q}(\alpha_{1}, \alpha_{2},\ldots, \alpha_{n})=\alpha.
$
\end{theorem}

\begin{IEEEproof}
It follows directly from formula~\eqref{25} by direct calculation.
\end{IEEEproof}

\begin{remark}
It should be noted that, by formula~\eqref{20}, the PFIWBM operator $\mathrm{PFIWBM}_{T}^{p,q}$
does not have the idempotency. For example, choose $\alpha_1=\alpha_2=\langle \frac{1}{3}, \frac{1}{3},
\frac{1}{3}\rangle$, $\omega=(\frac{1}{3}, \frac{2}{3})^{\mathrm{T}}$, $p=q=1$, and $T=T_{\mathbf{P}}$. Applying formula~\eqref{20},
by direct calculation, it follows that
$
\mathrm{PFIWBM}_{T}^{1,1}(\alpha_{1}, \alpha_{2}) \thickapprox \langle 0.17304, 0.22599,
0.60097\rangle
\neq \langle \frac{1}{3}, \frac{1}{3},
\frac{1}{3}\rangle.
$
\end{remark}

\begin{theorem}[{\textrm{\protect Boundedness}}]
\label{PFINWBM-Boundedness-Thm}
Let $\alpha_{i}=\left\langle\mu_{\alpha_{i}},\eta_{\alpha_{i}},\nu_{\alpha_{i}}\right\rangle$
($i=1,2,\ldots,n$) be a collection of PFNs, then
$\langle \min\limits_{i}\{\mu_{\alpha_{i}}\}, $
$\min\limits_{i}\{\eta_{\alpha_{i}}\},$
$\max\limits_{i}\{\nu_{\alpha_{i}}\}\rangle$
$\preccurlyeq_{_{\mathrm{W}}} \mathrm{PFINWBM}_{T}^{p,q}(\alpha_{1}, \alpha_{2}, \ldots ,\alpha_{n})$
$\preccurlyeq_{_{\mathrm{W}}} \langle \max\limits_{i}\{\mu_{\alpha_{i}}\},
1\!-(\max\limits_{i}\{\mu_{\alpha_{i}}\}\!+\min\limits_{i}
\{\nu_{\alpha_{i}}\}), \min\limits_{i}\{\nu_{\alpha_{i}}\}\rangle$.
\end{theorem}

\begin{IEEEproof}
Similarly to the proof of Theorem~\ref{PFIBM-Boundedness-Thm},
by using Theorem~\ref{PFINWBM-Formula-Thm}, it can be verified that this is true. 
\end{IEEEproof}

\section{A model for MCDM using picture fuzzy information}\label{Sec-6}

Because the PFINWBM operator not only has some expected properties, including the monotonicity,
idempotency, boundedness, and commutativity (see Theorems~\ref{PFIWBM-Commutativity-Thm}
--\ref{PFINWBM-Boundedness-Thm}), but also take into account the weights for criteria or
experts, this section use it for the MCDM problems under picture fuzzy environment.

For a multi-criteria decision making (MCDM) under picture fuzzy
environment, let $A=\{A_1, A_2, \ldots, A_m\}$ be a set of alternatives
to be selected, and $G=\{G_1, G_2,\ldots, G_n\}$ be a set of criteria to be evaluated
whose weight vector is $\omega=(\omega_1, \omega_2, \ldots, \omega_n)^{\mathrm{T}}$
such that $\omega_{j}\in (0, 1]$ and $\sum_{j=1}^{n}\omega_j=1$.
Assume that The performance of the alternative $A_{i}$ with respect to the criterion $G_{j}$
is measured by a PFN $\alpha_{ij}=\langle \mu_{ij}, \eta_{ij}, \nu_{ij}\rangle$,
where $\mu_{ij}$ is the degree of the positive membership; $\eta_{ij}$ be
the degree of neutral membership, and $\nu_{ij}$ indicates the degree that
the alternative $A_{i}$ does not satisfy the attribute $G_{j}$.
To rank the alternatives, the following steps are given:

\textbf{Step~1:}  (Construct the decision matrix) It is given that the
decision-maker gave their preference in the form of PFNs
$\alpha_{ij}=\langle \mu_{ij}, \eta_{ij}, \nu_{ij}\rangle$ towards the alternative
$A_{i}$ with respect to the criterion $G_{j}$ and hence construct a picture fuzzy decision
matrix $D=(\alpha_{ij})_{m\times n}$ as
\begin{table}[H]
\centering
\resizebox{\columnwidth}{!}{
	\begin{tabular}{ccccc}
		\toprule
		$$ & $G_{1}$ &  $G_{2}$ & $\cdots$ & $G_{n}$ \\
		\midrule
		$A_{1}$ & $\langle \mu_{1 1}, \eta_{1 1}, \nu_{1 1}\rangle$ &
        $\langle \mu_{1 2}, \eta_{1 2}, \nu_{1 2}\rangle$ &
           $\cdots$ & $\langle \mu_{1 n}, \eta_{1 n}, \nu_{1 n}\rangle$\\
		$A_{2}$ & $\langle \mu_{2 1}, \eta_{2 1}, \nu_{2 1}\rangle$ &
        $\langle \mu_{2 2}, \eta_{2 2}, \nu_{2 2}\rangle$ &
           $\ldots$ & $\langle \mu_{2 n}, \eta_{2 n}, \nu_{2 n}\rangle$ \\
        $\vdots$ & $\vdots$ & $\vdots$ &
           $\ddots$ & $\vdots$ \\
        $A_{m}$ & $\langle \mu_{m 1}, \eta_{m 1}, \nu_{m 1}\rangle$ &
        $\langle \mu_{m 2}, \eta_{m 2}, \nu_{m 2}\rangle$ &
           $\ldots$ & $\langle \mu_{m n}, \eta_{m n}, \nu_{m n}\rangle$ \\
        \bottomrule
	\end{tabular}
}
\end{table}

\textbf{Step~2:} (Normalize the decision matrix) Transform the picture fuzzy decision
matrix $D=(\alpha_{ij})_{m\times n}$ into the normalized picture fuzzy decision
matrix $R=(r_{ij})_{m\times n}$ as follows:
$$
r_{ij}=
\begin{cases}
\alpha_{ij}, & \text{for benefit criteria } G_{j}, \\
\alpha_{ij}^{\complement}, & \text{for cost criteria } G_{j},
\end{cases}
$$
where $\alpha_{ij}^{\complement}=\langle \nu_{ij}, \eta_{ij}, \mu_{ij}\rangle$
is the complement of $\alpha_{ij}$.

\textbf{Step~3:} (Calculate the aggregated value) Based
on the normalized decision matrix $R$ obtained from step~2,
the overall aggregated value of every alternative $A_{i}$ ($i=1,2,\ldots, m$),
under the different criterias $G_{1}, G_{2}, \ldots, G_{n}$,
is obtained by using PFINWBM operator $\mathrm{PFINWBM}_{T}^{p,q}$
(in general, we can take $p=q=1$) for some strict t-norm $T$,
and hence get the collective value $r_{i}$ for each alternative $A_{i}$:
$$
r_{i}=
\mathrm{PFINWBM}_{T}^{p,q}(r_{i1}, r_{i2}, \ldots, r_{in}).
$$

\textbf{Step~4:} (Rank the alternative) Rank the alternatives $A_{1}, A_{2},
  \ldots, A_{m}$ by using the total order defined in Definition~\ref{Wu-Order-Def}
  and select the most desirable alternative.

\section{An illustrative example}\label{Sec-7}

In this section, we utilize a practical MCDM problem to illustrate the application of our developed approach in Section~\ref{Sec-6}.
Suppose an organization plans to implement enterprise resource planning (ERP) system (adapted from~\cite{LXW2007}).
The first step is to form a project team consisting of Chief Information Officer and two senior representatives
from related departments. By collecting all possible information about ERP vendors and systems, the project team
chooses five potential ERP systems $A_{i}$ ($i = 1, 2, \ldots, 5$) as candidates. The project team selects four
criteria to evaluate the alternatives: (1) $G_1$ is function and technology; (2) $G_2$ is strategic fitness;
(3) $G_3$ is vendors ability; (4) $G_4$ is vendors reputation. The five possible ERP systems
$A_{i}$ ($i=1, 2, \ldots, 5$) are to be evaluated using the PFNs by the decision makers under the above
four criteria whose weighting vector is $\omega = (0.2, 0.1, 0.3, 0.4)^{\mathrm{T}}$. Since the problem being addressed
does not have any cost criteria, picture fuzzy decision matrix $D=(\alpha_{ij})_{5\times 4}$ is same as the
normalized picture fuzzy decision matrix $R=(r_{ij})_{5\times4}$ shown in Table~\ref{tab1}.
\begin{table}[H]
\centering
\caption{The picture fuzzy decision matrix $D$}
\label{tab1}
\resizebox{\columnwidth}{!}{
\begin{tabular}{lcccc}
\toprule
 ~ &$G_{1}$ &$G_{2}$ &$G_{3}$ &$G_{4}$\\
 \midrule
$A_{1}$ & $\langle0.53, 0.33, 0.09\rangle$ & $\langle0.89, 0.08, 0.03\rangle$ & $\langle0.42, 0.35, 0.18\rangle$ & $\langle0.08, 0.89, 0.02\rangle$\\
$A_{2}$ & $\langle0.73, 0.12, 0.08\rangle$ & $\langle0.13, 0.64, 0.21\rangle$ & $\langle0.03, 0.82, 0.13\rangle$ & $\langle0.73, 0.15, 0.08\rangle$\\
$A_{3}$ & $\langle0.91, 0.03, 0.02\rangle$ & $\langle0.07, 0.09, 0.05\rangle$ & $\langle0.04, 0.85, 0.10\rangle$ & $\langle0.68, 0.26, 0.06\rangle$\\
$A_{4}$ & $\langle0.85, 0.09, 0.05\rangle$ & $\langle0.74, 0.16, 0.10\rangle$ & $\langle0.02, 0.89, 0.05\rangle$ & $\langle0.08, 0.84, 0.06\rangle$\\
$A_{5}$ & $\langle0.90, 0.05, 0.02\rangle$ & $\langle0.68, 0.08, 0.21\rangle$ & $\langle0.05, 0.87, 0.06\rangle$ & $\langle0.13, 0.75, 0.09\rangle$\\
\bottomrule
\end{tabular}
}
\end{table}

In the following, we use the PFINWBM operator $\mathrm{PFINWBM}^{p,q}_{T_{\gamma}^{\mathbf{H}}}$
induced by the Hamacher t-norms $T_{\gamma}^{\mathbf{H}}$ ($\gamma \in (0, +\infty)$) to
select the best ERP system.

\textbf{Step 1.}  According to Table \ref{tab1}, aggregate all PFNs by using
$\mathrm{PFINWBM}^{p,q}_{T_{\gamma}^{\mathbf{H}}}$ ($\gamma \in (0, +\infty)$)
to obtain the overall PFNs $r_{i}$ ($i=1, 2, \ldots, 5$) of the alternative $A_{i}$ ($i=1, 2, \ldots, 5$),
where we take $T$ as Hamacher t-norms $T_{\gamma}^{\mathbf{H}}$ ($\gamma \in (0, +\infty)$)
as follows:
$r_{i}\!=\mathrm{PFINWBM}_{T_{\gamma}^{\mathbf{H}}}^{p,q}(r_{i1}, r_{i2}, r_{i3}, r_{i4})
\!=\Big \langle\frac{\gamma}{\left(\frac{\gamma}{\hat{\mu}}\!+1
\!-\gamma\right)^{\frac{1}{p\!+q}}\!+\gamma\!-1},
\frac{\gamma}{\left(\frac{\gamma}{{\hat{\eta}\!+\hat{\mu}}}
\!+1\!-\gamma\right)^{\frac{1}{p\!+q}}\!+\gamma\!-1}
\!-\frac{\gamma}{\left(\frac{\gamma}{\hat{\mu}}
\!+1\!-\gamma\right)^{\frac{1}{p\!+q}}\!+\gamma\!-1},
1\!-\frac{\gamma}{\left(\frac{\gamma}{1\!-\hat{\nu}}
\!+1\!-\gamma\right)^{\frac{1}{p\!+q}}\!+\gamma\!-1}\Big\rangle$,
where
$\hat{\mu}\!= 1\!-\frac{\gamma}{\prod\limits_{j,k=1 \atop j\neq k}^{4}\left[\frac{\gamma^{2}}
{\left(\frac{\gamma}{\mu_{ij}}\!+1\!-\gamma\right)^{p}\left(\frac{\gamma}{\mu_{ik}}\!+1\!-\gamma
\right)^{q}\!-1}\!+1\right]^{\frac{\omega_{j}\omega_{k}}{1\!-\omega_{j}}}\!+\gamma\!-1}$,
\begin{align*}
\hat{\eta}=&\!+\frac{\gamma}{\prod\limits_{j,k=1 \atop j\neq k}^{4}
\left[\frac{1}{\frac{1}{\overline{\eta}_1}
\!-\frac{1}{\overline{\eta}_2}
\!+\frac{1}{\gamma}\!-\frac{1}{\overline{\eta}_3}}
\!+1\!-\gamma\right]^{\frac{\omega_{j}\omega_{k}}{1\!-\omega_{j}}}\!+\gamma\!-1}\\
&\!-\frac{\gamma}{\prod\limits_{j,k=1 \atop j\neq k}^{4}
\left[\frac{\gamma^{2}}{\left(\frac{\gamma}{1\!-\nu_{ij}}\!+1\!-\gamma\right)^{p}
\left(\frac{\gamma}{1\!-\nu_{ik}}\!+1\!-\gamma\right)^{q}\!-1}\!+
1\right]^{\frac{\omega_{j}\omega_{k}}{1\!-\omega_{j}}}\!+\gamma\!-1},
\end{align*}
$\overline{\eta}_1\!=\left(\frac{\gamma}{\eta_{ij}
\!+\mu_{ij}}\!+1\!-\gamma\right)^{p}\left(\frac{\gamma}{\eta_{ik}
\!+\mu_{ik}}\!+1\!-\gamma\right)^{q}\!+\gamma\!-1$,
$\overline{\eta}_2=\left(\frac{\gamma}{\mu_{ij}}\!+1\!-\gamma\right)^{p}
\left(\frac{\gamma}{\mu_{ik}}\!+1\!-\gamma\right)^{q}\!+\gamma\!-1$, $\overline{\eta}_3=\left(\frac{\gamma}{1\!-\nu_{ij}}\!+1\!-\gamma\right)^{p}
\left(\frac{\gamma}{1\!-\nu_{ik}}\!+1\!-\gamma\right)^{q}\!+\gamma\!-1$,
and
$\hat{\nu}\!=\frac{\gamma}{\prod\limits_{j,k\!=1 \atop j\neq k}^{4}
\left[\frac{\gamma^{2}}{\left(\frac{\gamma}{1\!-\nu_{ij}}\!+1\!-\gamma\right)^{p}
\left(\frac{\gamma}{1\!-\nu_{ik}}\!+1\!-\gamma\right)^{q}\!-1}
\!+1\right]^{\frac{\omega_{j}\omega_{k}}{1\!-\omega_{j}}}\!+\gamma\!-1}$.

Let $p=q = 1$. When $\gamma=2$, $T_{2}^{\mathbf{H}}$ is the Einstein operator, and the aggregated results
are shown in Table~\ref{tab2}.

\begin{table}[H]
\centering
\caption{The aggregated values of the ERP systems}
\label{tab2}
\resizebox{\columnwidth}{!}{
\begin{tabular}{p{5cm} cp{cm}}
\toprule
 ~ & $\mathrm{PFINWBM}_{T_{2}^{\mathbf{H}}}^{1,1}$\\
 \midrule
$r_{1}$ &$\langle0.3749, 0.5173, 0.0774\rangle$\\
$r_{2}$ &$\langle0.4403, 0.4039, 0.1079\rangle$\\
$r_{3}$ &$\langle0.4789, 0.3357, 0.0598\rangle$\\
$r_{4}$ &$\langle0.2901, 0.6297, 0.0587\rangle$\\
$r_{5}$ &$\langle0.3295, 0.5697, 0.0732\rangle$\\
\bottomrule
\end{tabular}
}
\end{table}

\textbf{Step 2.} According to the aggregated values shown in Table~\ref{tab2}, the score
of the ERP systems are shown in Table~\ref{tab3}.
\begin{table}[H]
\centering
\caption{The scores of the ERP systems}\label{tab3}
\resizebox{\columnwidth}{!}{
\begin{tabular}{p{6cm} cp{cm}}
\toprule
 ~ & $\mathrm{PFINWBM}_{T_{2}^{\mathbf{H}}}^{1,1}$ \\
 \midrule
$S_{r_{1}}$ &0.2975\\
$S_{r_{2}}$ &0.3324\\
$S_{r_{3}}$ &0.4191\\
$S_{r_{4}}$ &0.2313\\
$S_{r_{5}}$ &0.2563\\
\bottomrule
\end{tabular}
}
\end{table}

\textbf{Step 3.} According to the scores in Table~\ref{tab3}, we have
$$
S_{r_{3}}> S_{r_{2}}> S_{r_{1}} > S_{r_{5}} > S_{r_{4}},
$$
and thus
$$
A_{3}\succ_{_{\mathrm{W}}} A_{2}\succ_{_{\mathrm{W}}} A_{1}\succ_{_{\mathrm{W}}}
A_{5}\succ_{_{\mathrm{W}}} A_{4},
$$
i.e., $A_3$ is the best EPR system.

\subsection{The influence of the parameters $p$ and $q$ for MCDM results}

To illustrate the influence of the parameters $p$ and $q$ in the above example, we use different
values of parameters $p$ and $q$ in the PFINWBM operator $\mathrm{PFINWBM}_{T_{2}^{\mathbf{H}}}^{p,q}$.
The ranking results are shown in Table~\ref{tab4}. As we can see from Table~\ref{tab4},
the ranking of the ERP systems using different values of parameters $p$ and $q$ in aggregation process
is slightly different. However, the best ERP system is $A_{3}$ for all different combination of parameters.
\begin{table}[H]
\centering
\caption{Ranking results of different parameters $p$ and $q$
obtained by $\mathrm{PFINWBM}_{T_{2}^{\mathbf{H}}}^{p,q}$}
\label{tab4}
\resizebox{\columnwidth}{!}{
\begin{tabular}{p{5cm} cp{cm}}
\toprule
 ~ & ranking \\
 \midrule
$p=2$, $q=1$ & $A_{3}\succ_{_{\mathrm{W}}} A_{2}\succ_{_{\mathrm{W}}} A_{1}
\succ_{_{\mathrm{W}}} A_{5}\succ_{_{\mathrm{W}}} A_{4}$ \\
$p=1$, $q=2$ & $A_{3}\succ_{_{\mathrm{W}}} A_{2}\succ_{_{\mathrm{W}}} A_{1}
\succ_{_{\mathrm{W}}} A_{5}\succ_{_{\mathrm{W}}} A_{4}$ \\
$p=3$, $q=1$ & $A_{3}\succ_{_{\mathrm{W}}} A_{5}\succ_{_{\mathrm{W}}} A_{4}
\succ_{_{\mathrm{W}}} A_{2}\succ_{_{\mathrm{W}}} A_{1}$ \\
$p=7$, $q=2$ & $A_{3}\succ_{_{\mathrm{W}}} A_{5}\succ_{_{\mathrm{W}}} A_{4}
\succ_{_{\mathrm{W}}} A_{1}\succ_{_{\mathrm{W}}} A_{2}$ \\
$p=8$, $q=1$ & $A_{3}\succ_{_{\mathrm{W}}} A_{5}\succ_{_{\mathrm{W}}} A_{1}
\succ_{_{\mathrm{W}}} A_{4}\succ_{_{\mathrm{W}}} A_{2}$ \\
$p=6$, $q=9$ & $A_{3}\succ_{_{\mathrm{W}}} A_{5}\succ_{_{\mathrm{W}}} A_{4}
\succ_{_{\mathrm{W}}} A_{2}\succ_{_{\mathrm{W}}} A_{1}$ \\
$p=10$, $q=10$ & $A_{3}\succ_{_{\mathrm{W}}} A_{5}\succ_{_{\mathrm{W}}} A_{4}
\succ_{_{\mathrm{W}}} A_{2}\succ_{_{\mathrm{W}}} A_{1}$ \\
\bottomrule
\end{tabular}
}
\end{table}

\begin{figure}[H]
\centering
\subfigure[Scores for $A_{1}$]{\scalebox{0.35}{\includegraphics[]{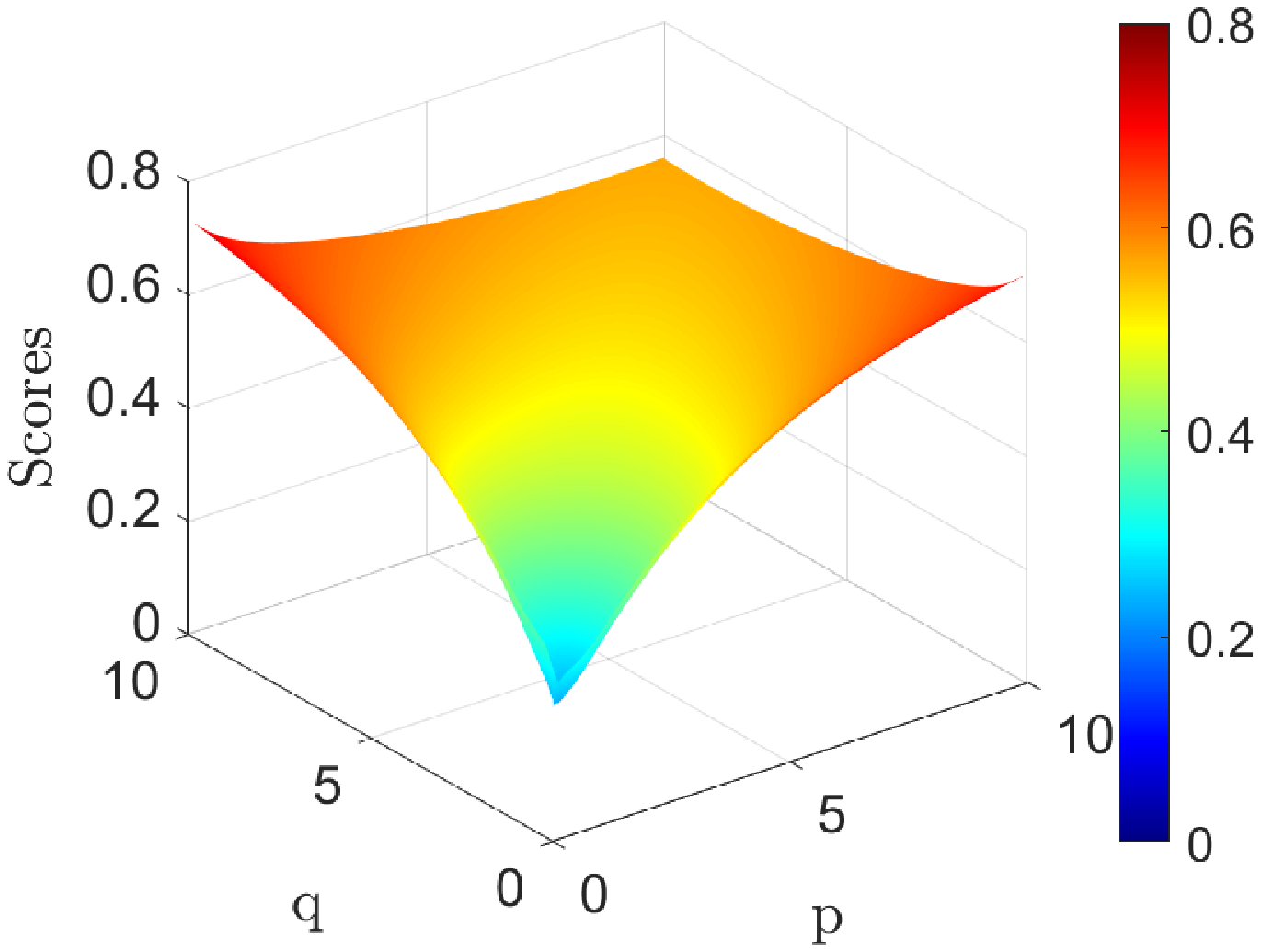}}}
\subfigure[Scores for $A_{2}$]{\scalebox{0.35}{\includegraphics[]{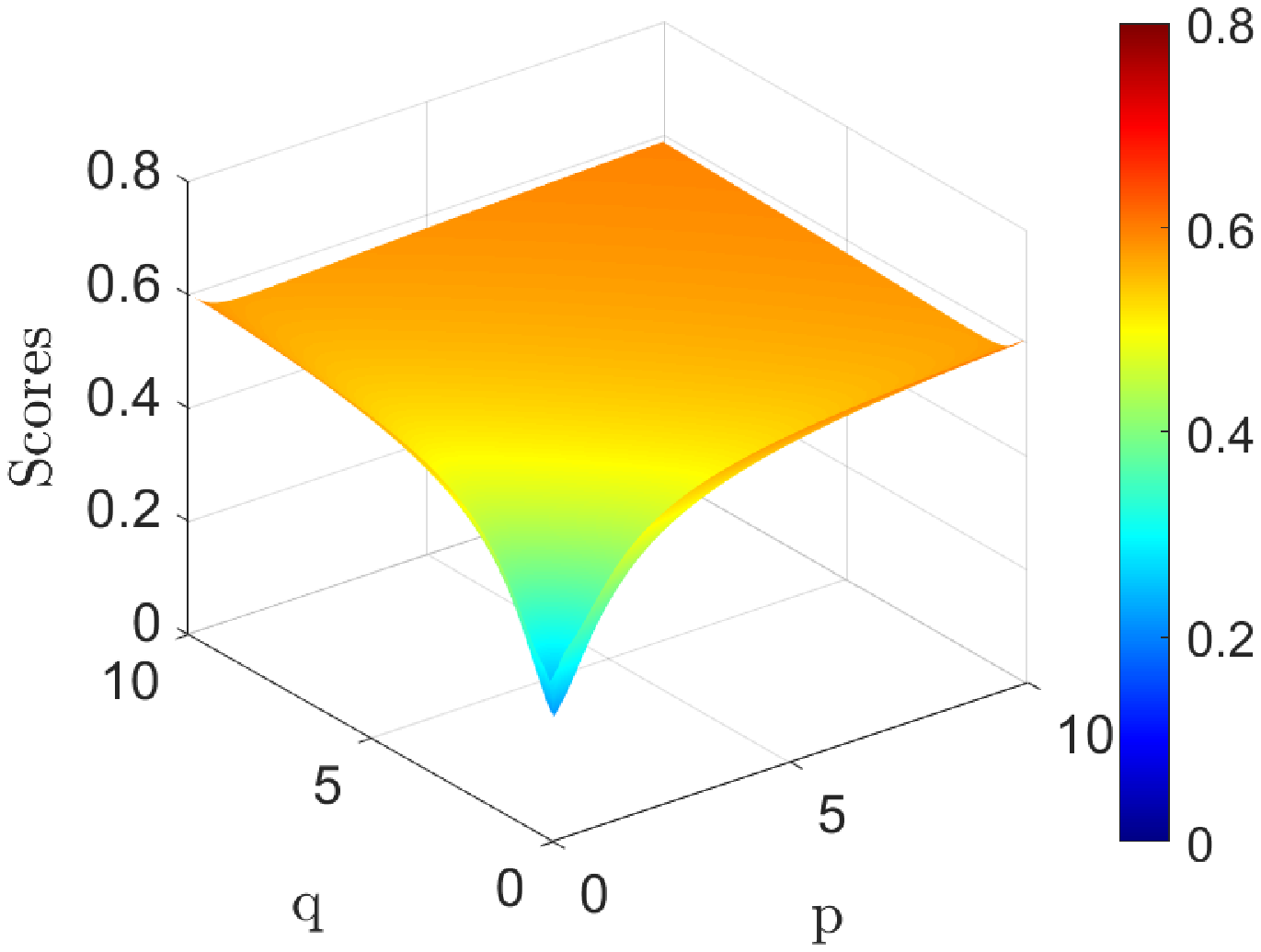}}}
\subfigure[Scores for $A_{3}$]{\scalebox{0.35}{\includegraphics[]{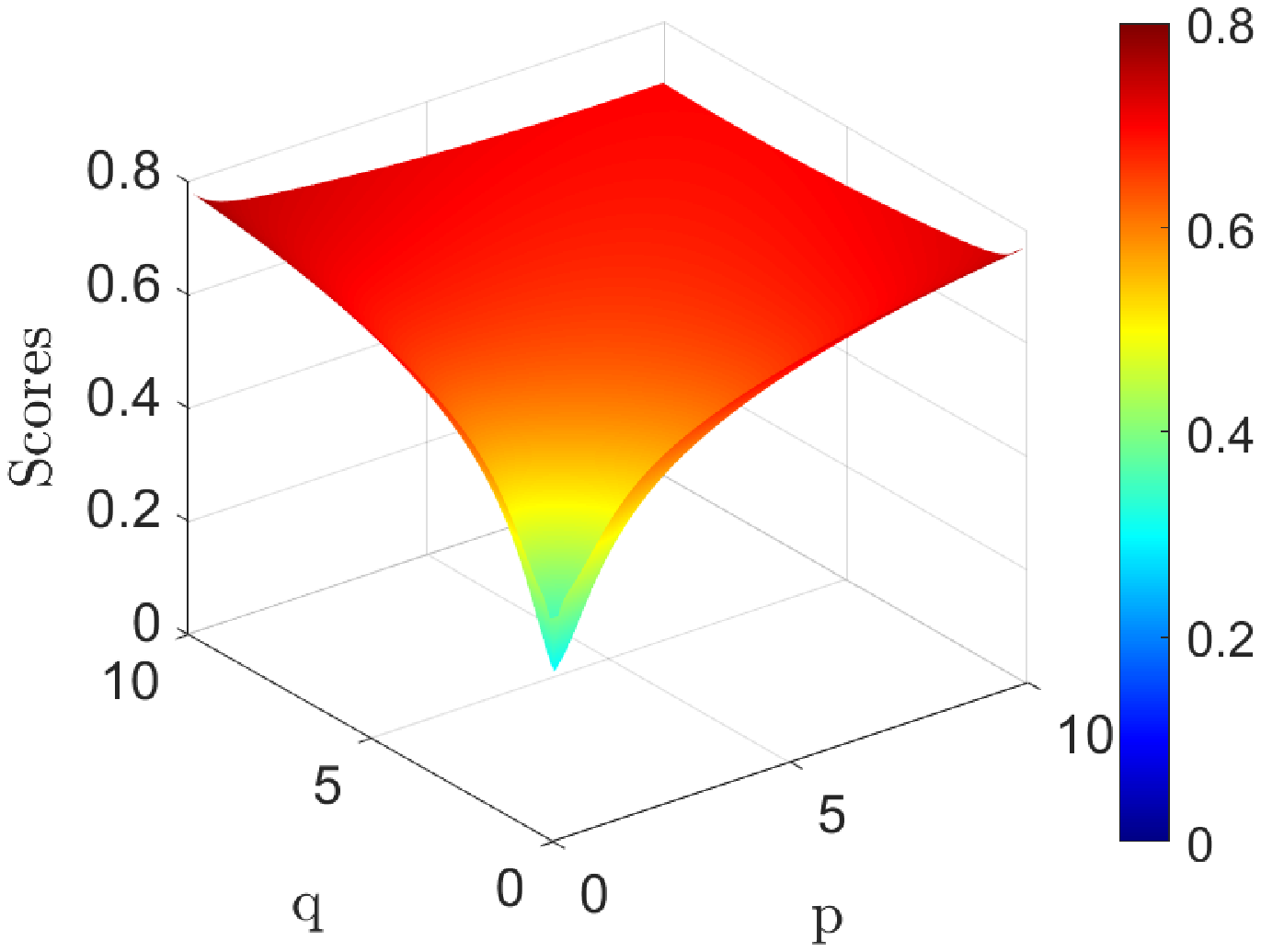}}}
\subfigure[Scores for $A_{4}$]{\scalebox{0.35}{\includegraphics[]{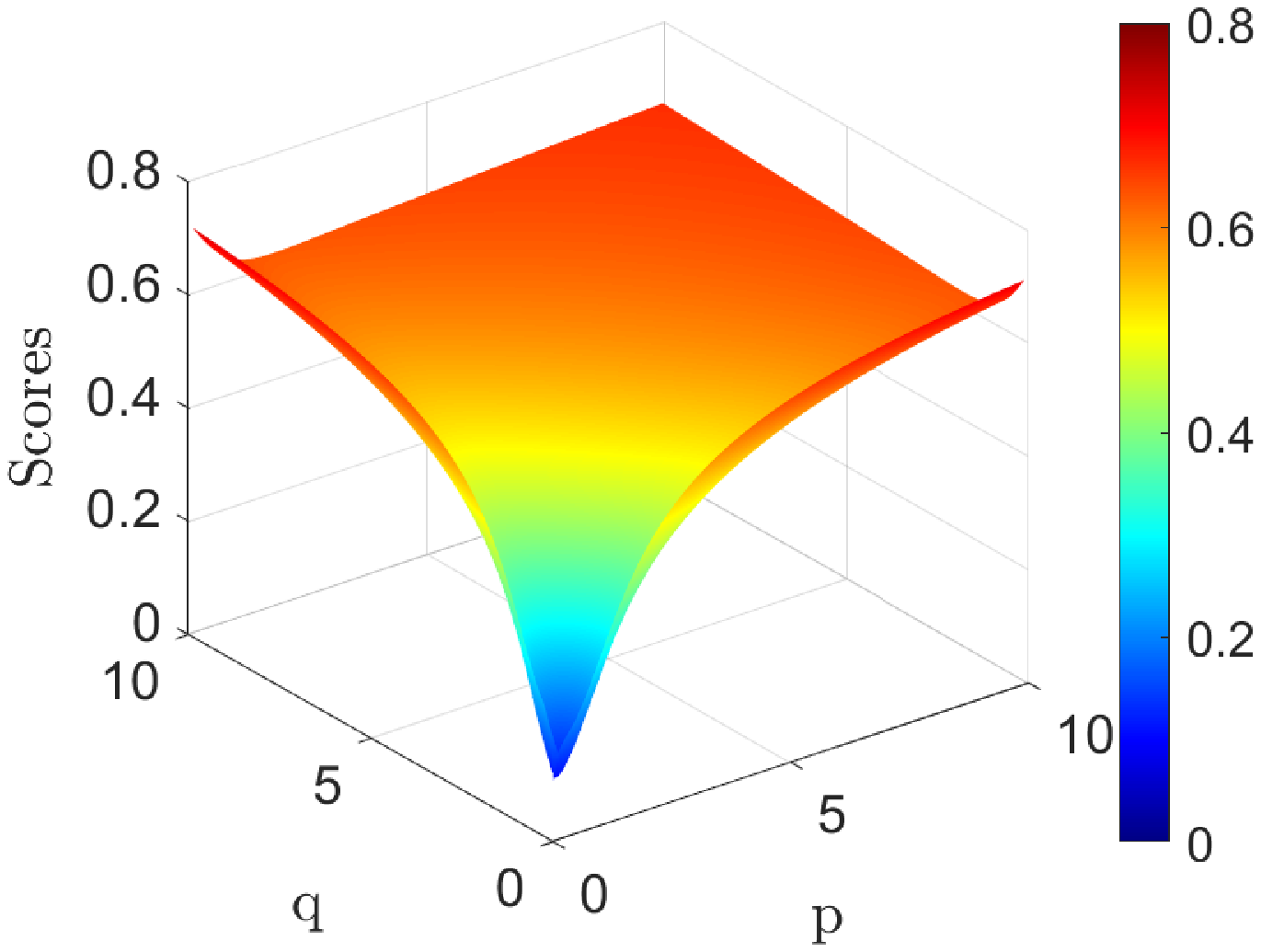}}}
\subfigure[Scores for $A_{5}$]{\scalebox{0.35}{\includegraphics[]{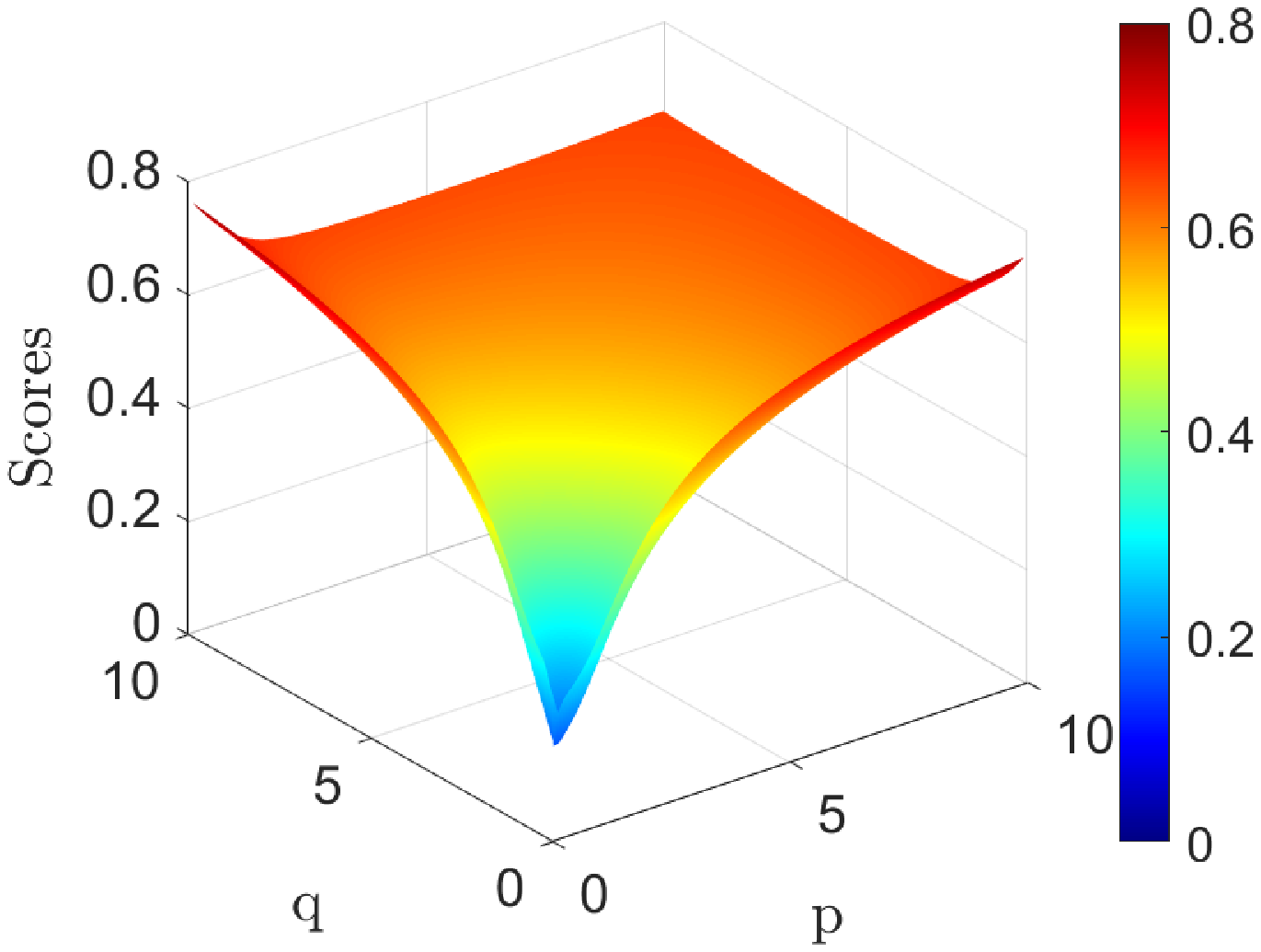}}}
\caption{Scores for alternatives in different values of $p$, $q$
obtained by $\mathrm{PFINWBM}_{T_{2}^{\mathbf{H}}}^{p,q}$}\label{fig-H}
\end{figure}

\begin{figure}[H]
\centering
\subfigure[Scores for $A_{1}$]
{\scalebox{0.35}{\includegraphics[]{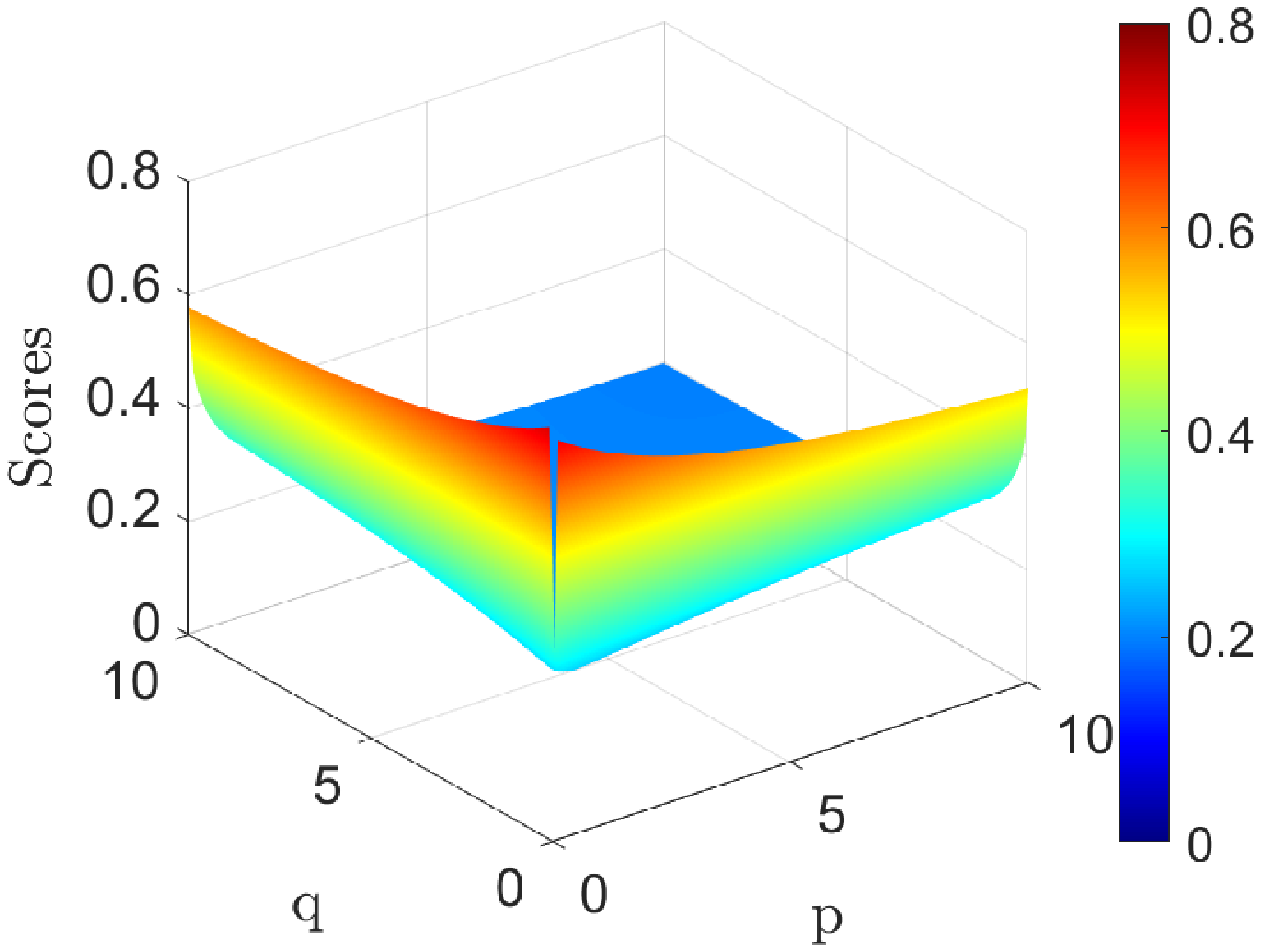}}}
\subfigure[Scores for $A_{2}$]
{\scalebox{0.35}{\includegraphics[]{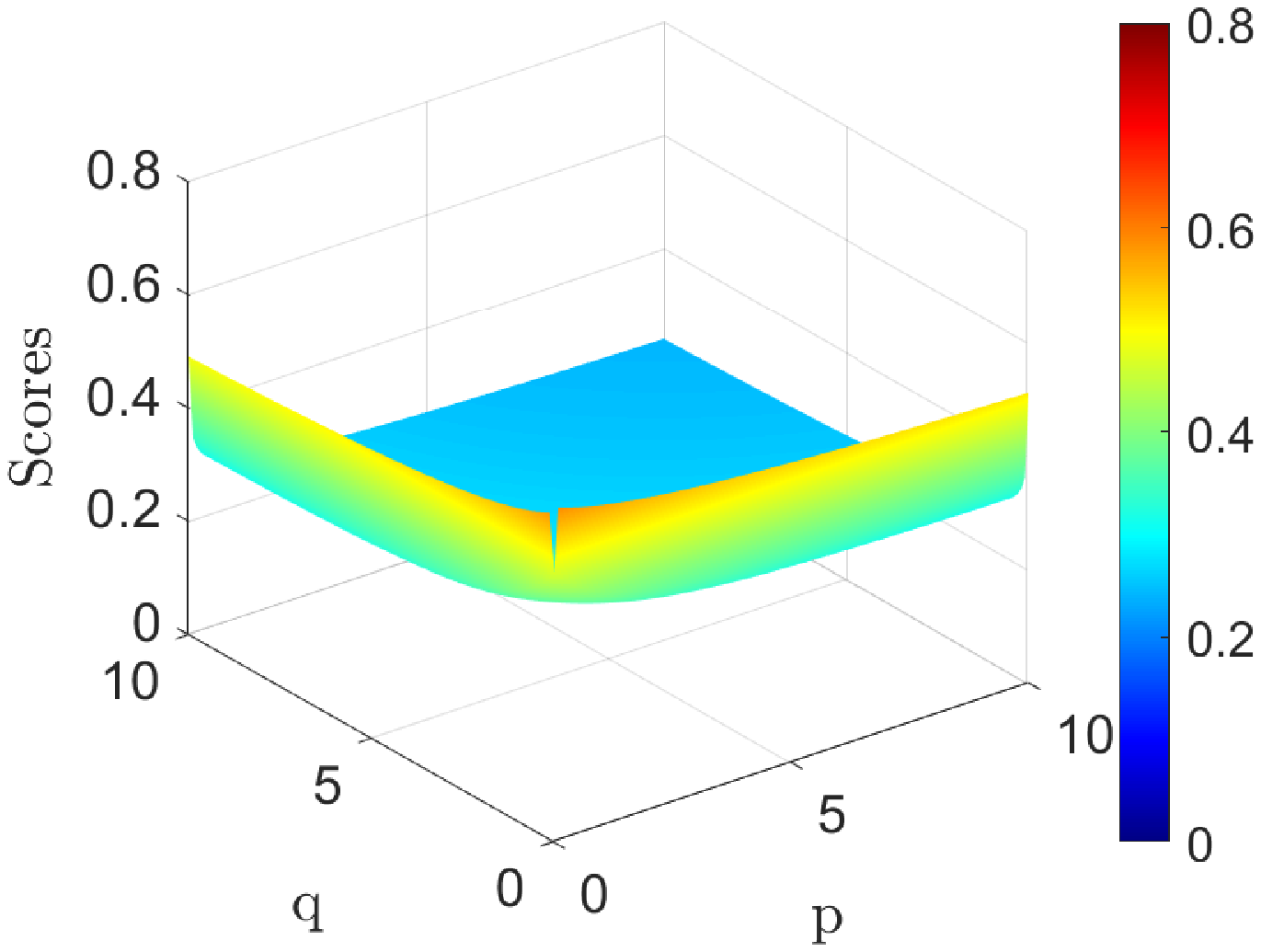}}}
\subfigure[Scores for $A_{3}$]
{\scalebox{0.35}{\includegraphics[]{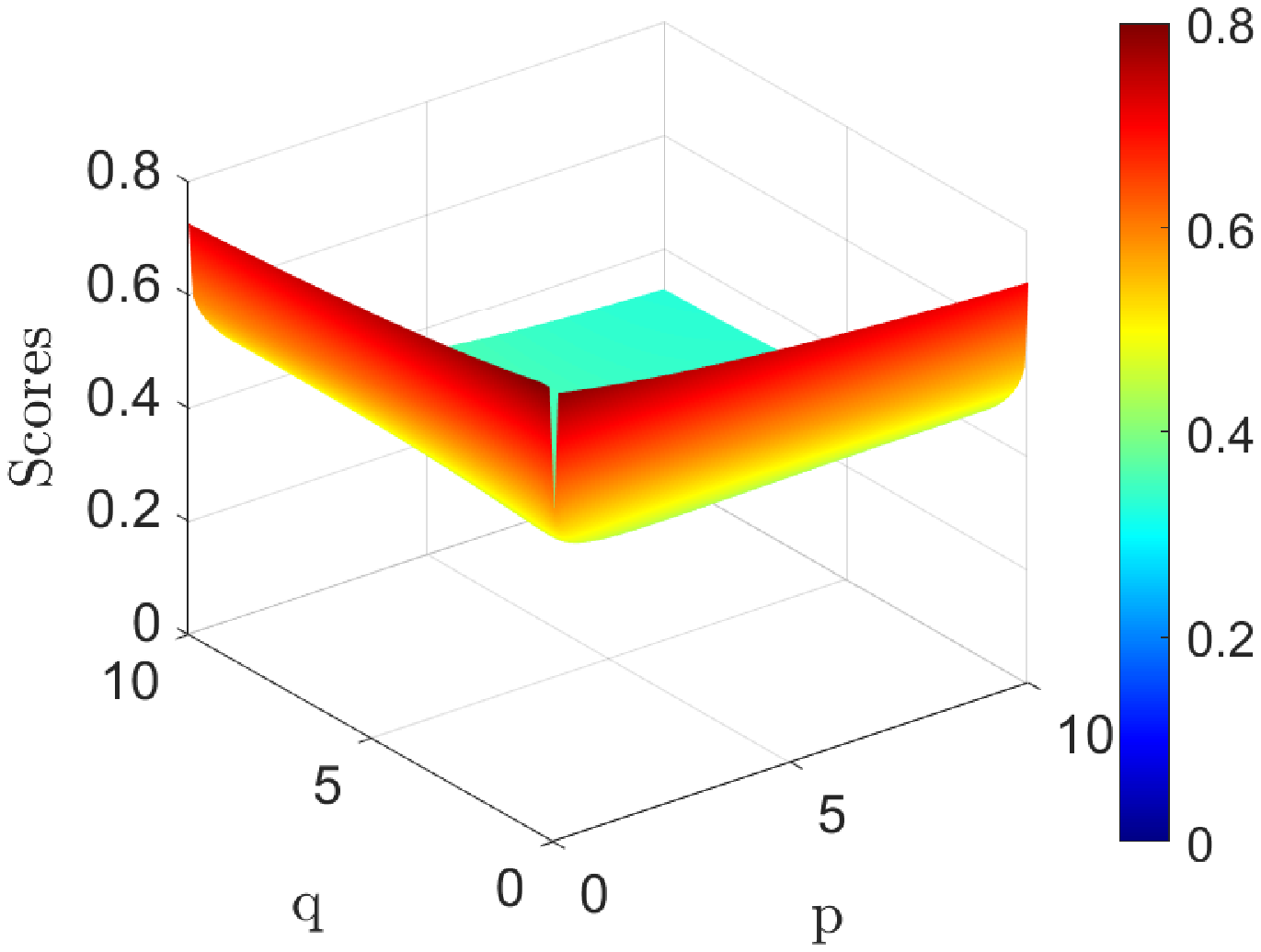}}}
\subfigure[Scores for $A_{4}$]
{\scalebox{0.35}{\includegraphics[]{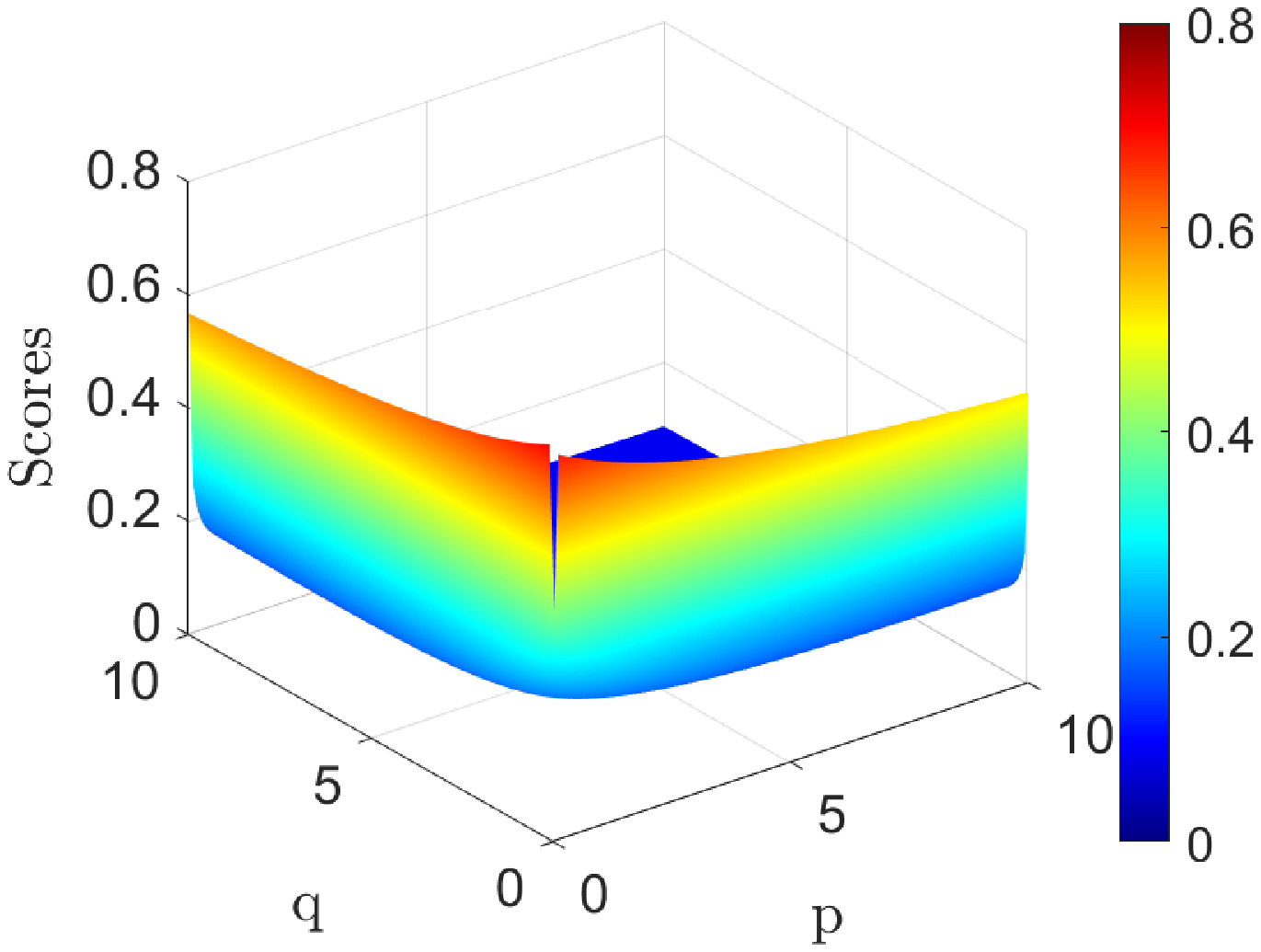}}}
\subfigure[Scores for $A_{5}$]
{\scalebox{0.35}{\includegraphics[]{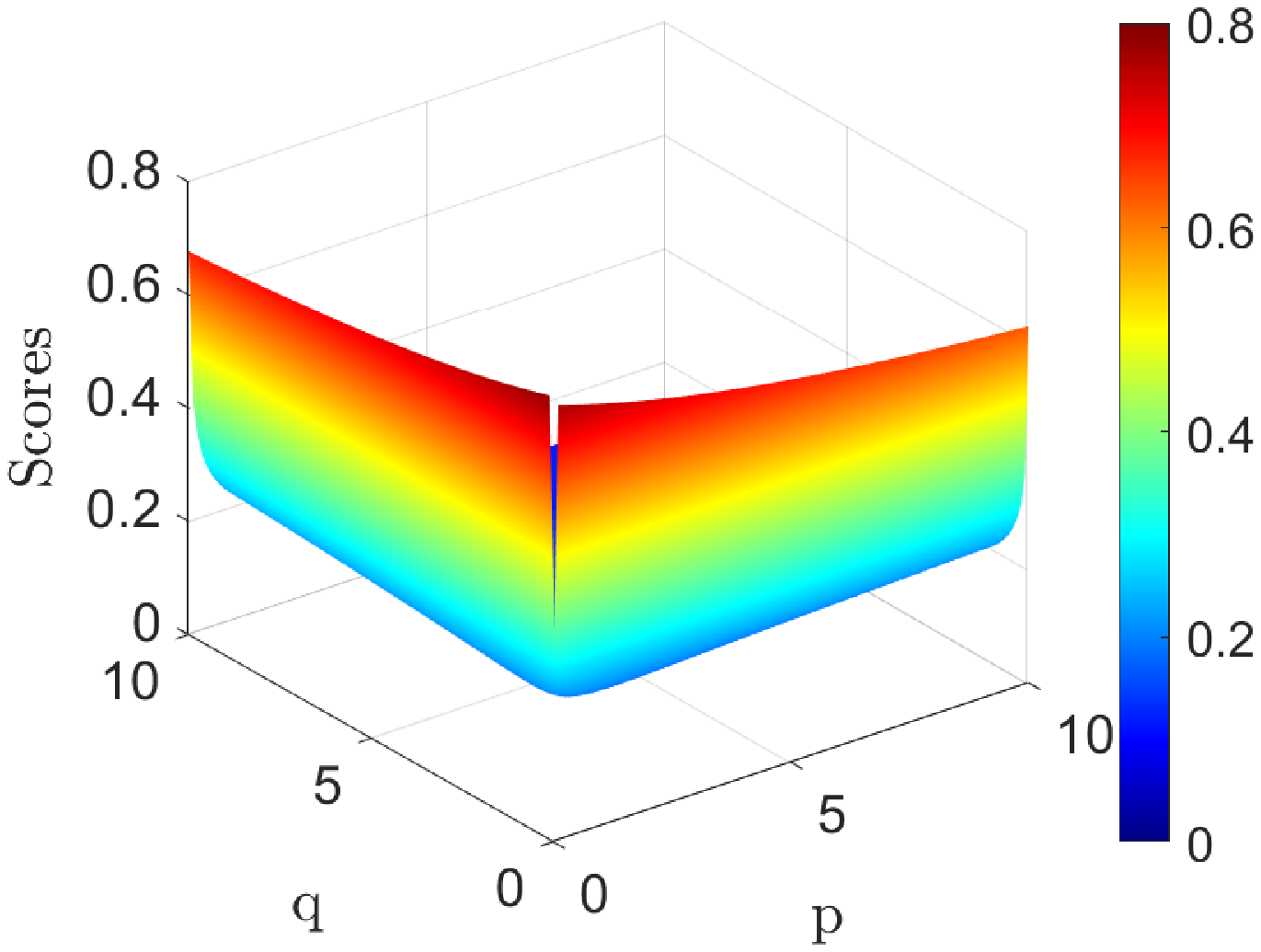}}}
\caption{Scores for alternatives in different values of $p$, $q$
obtained by $\mathrm{PFINWBM}_{T_{-2}^{\mathbf{SS}}}^{p,q}$}\label{fig-SS}
\end{figure}

\begin{figure}[H]
\centering
\subfigure[Scores for $A_{1}$]
{\scalebox{0.35}{\includegraphics[]{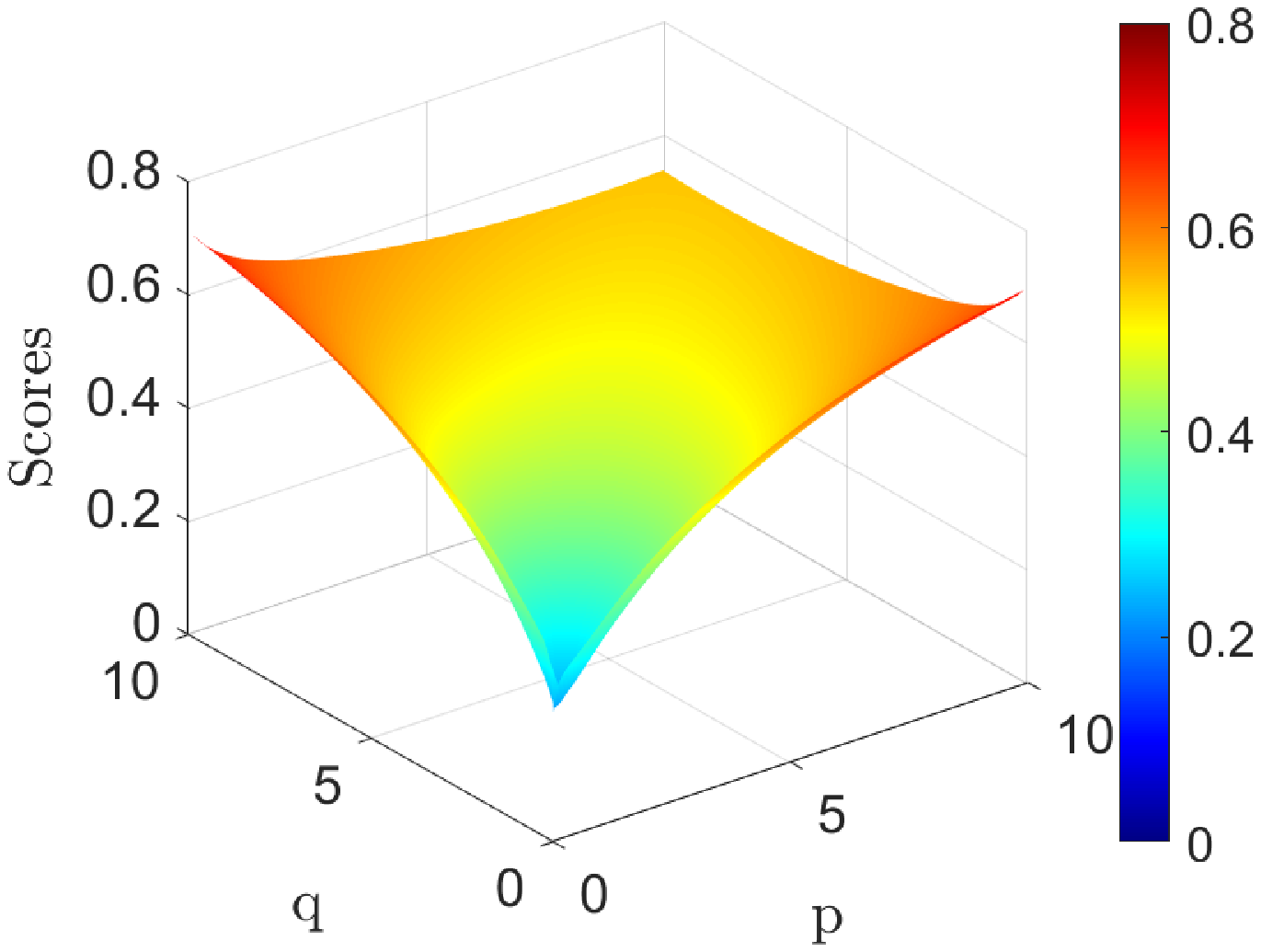}}}
\subfigure[Scores for $A_{2}$]
{\scalebox{0.35}{\includegraphics[]{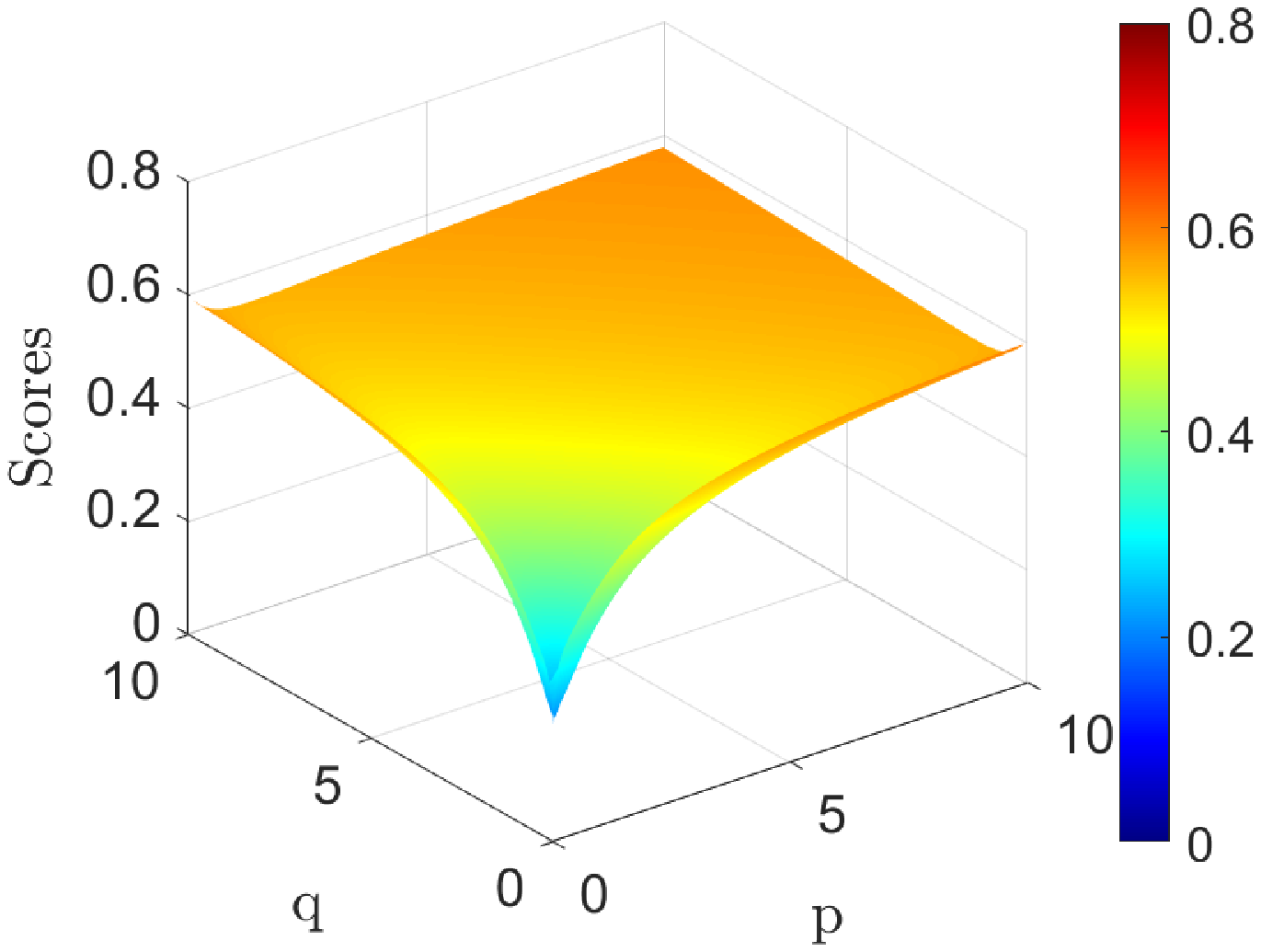}}}
\subfigure[Scores for $A_{3}$]
{\scalebox{0.35}{\includegraphics[]{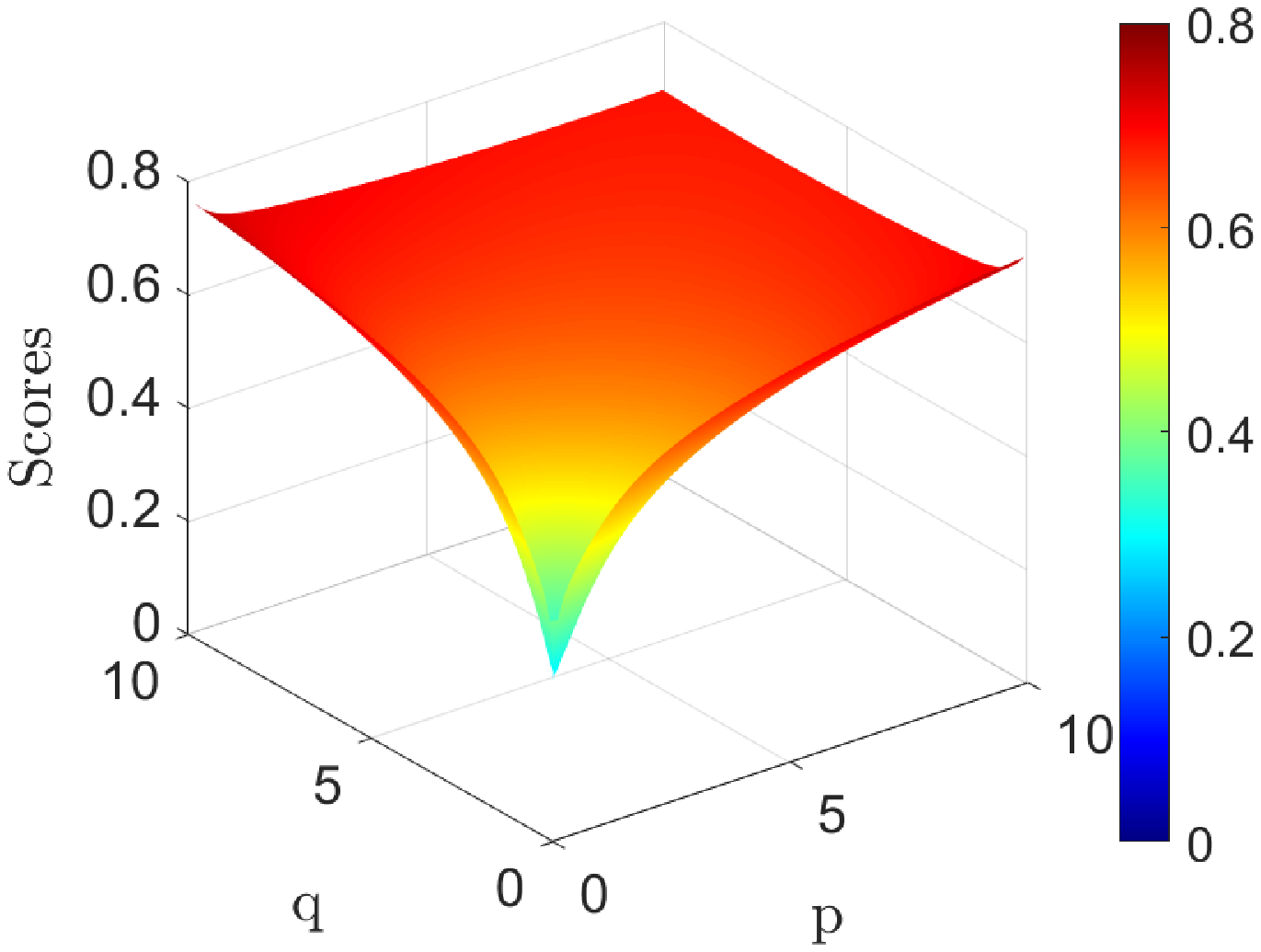}}}
\subfigure[Scores for $A_{4}$]
{\scalebox{0.35}{\includegraphics[]{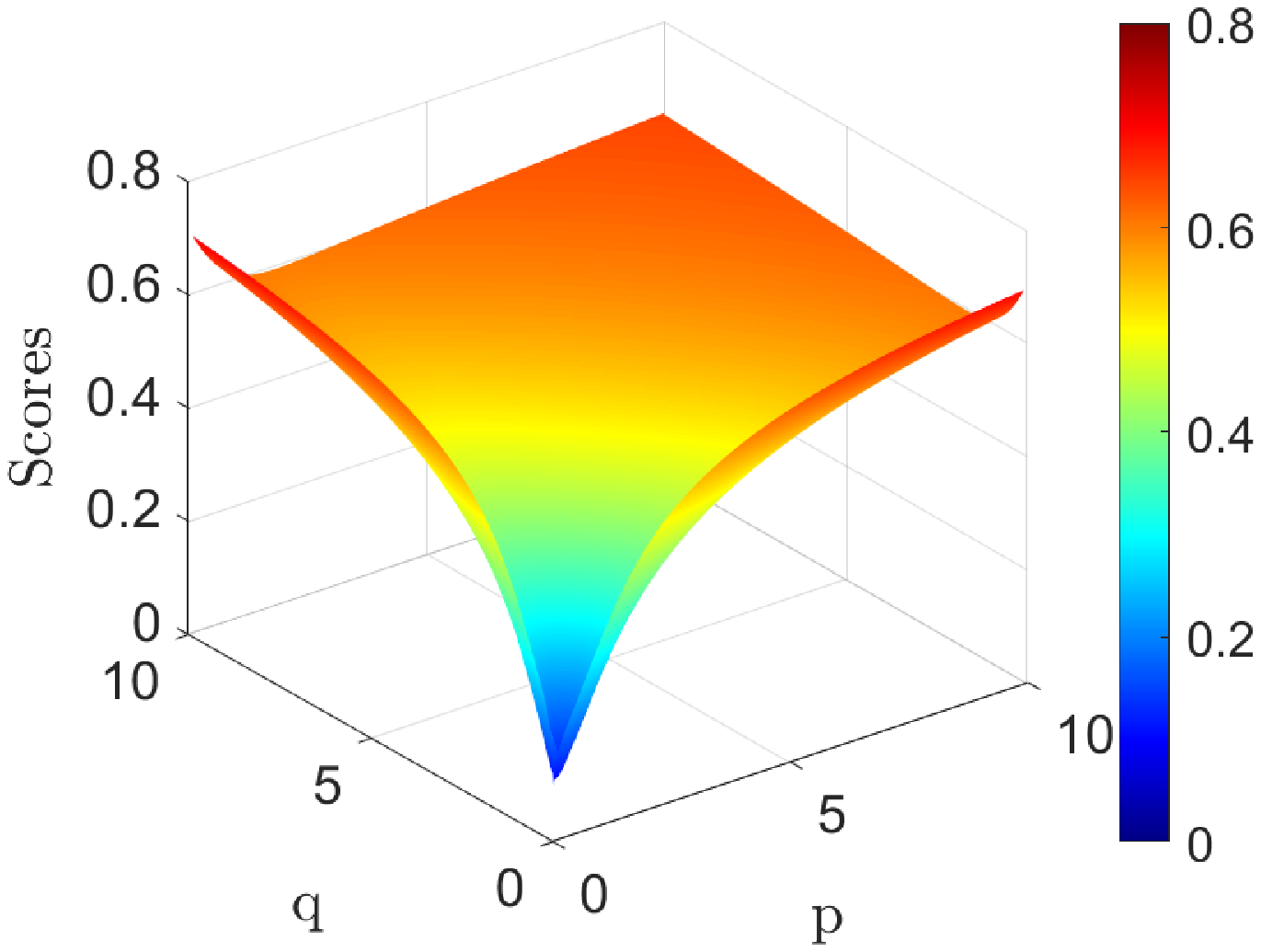}}}
\subfigure[Scores for $A_{5}$]
{\scalebox{0.35}{\includegraphics[]{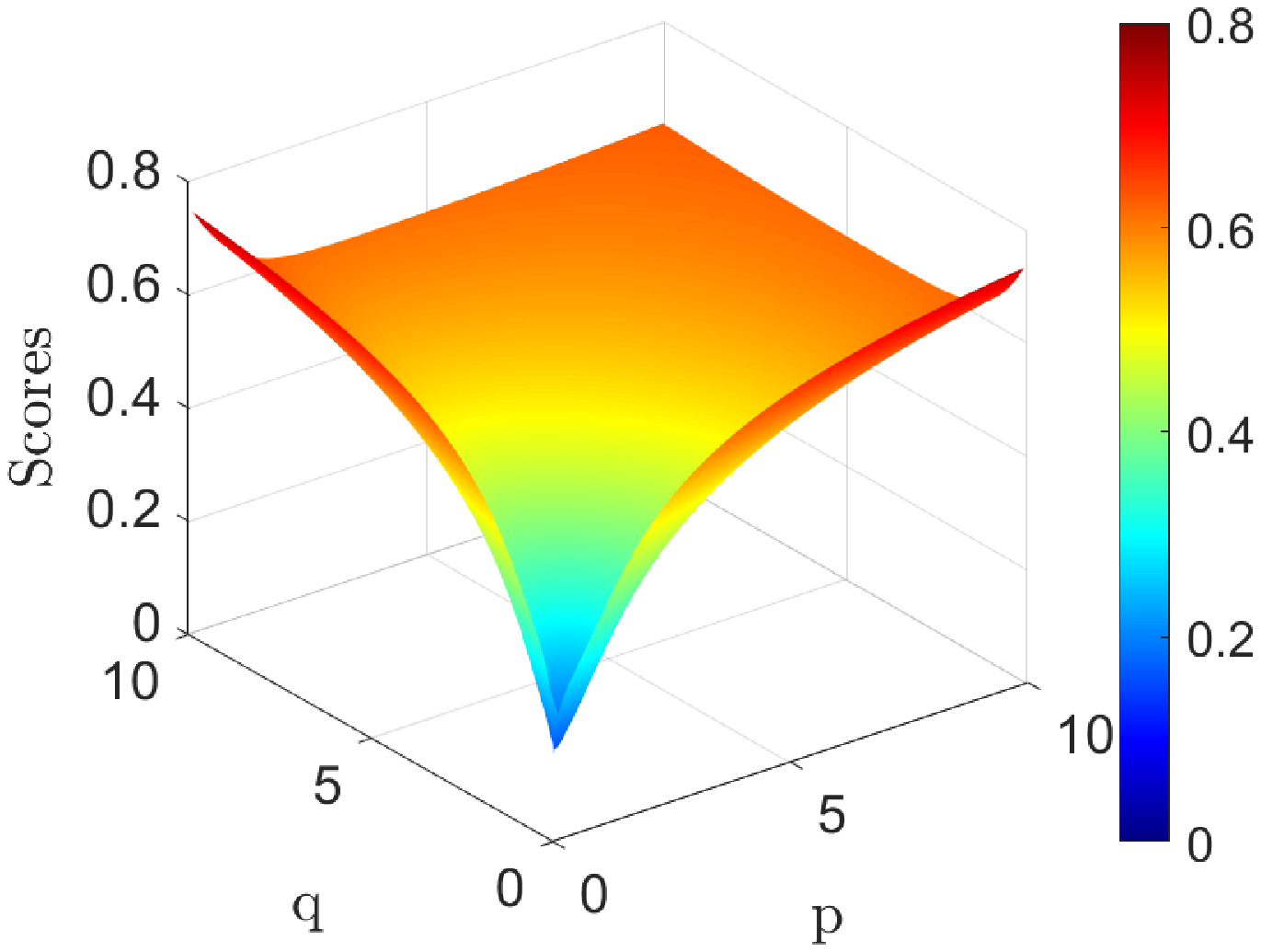}}}
\caption{Scores for alternatives in different values of $p$, $q$
obtained by $\mathrm{PFINWBM}_{T_{2}^{\mathbf{F}}}^{p,q}$}\label{fig-F}
\end{figure}

\begin{figure}[H]
\centering
\subfigure[Scores for $A_{1}$]
{\scalebox{0.35}{\includegraphics[]{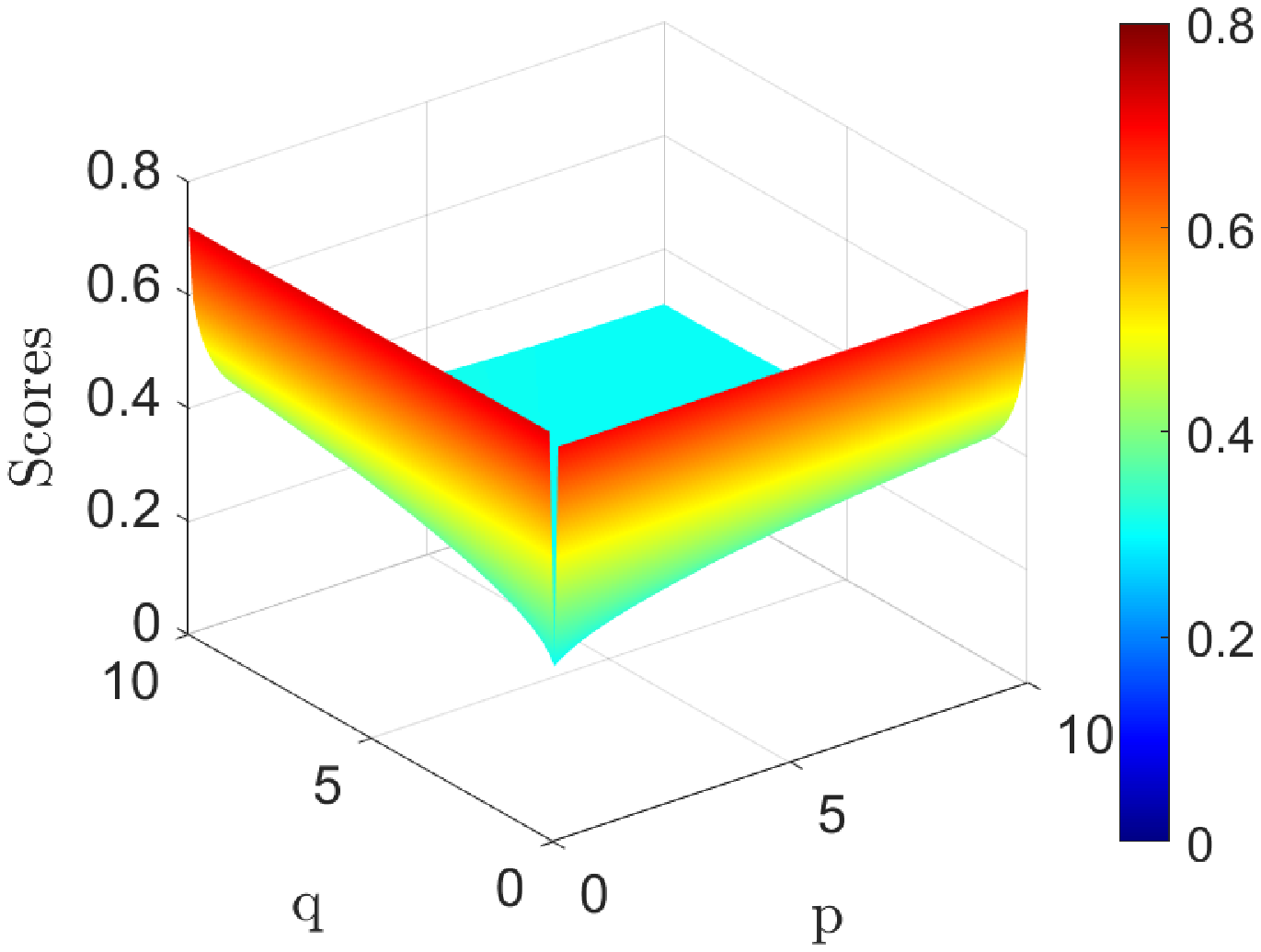}}}
\subfigure[Scores for $A_{2}$]
{\scalebox{0.35}{\includegraphics[]{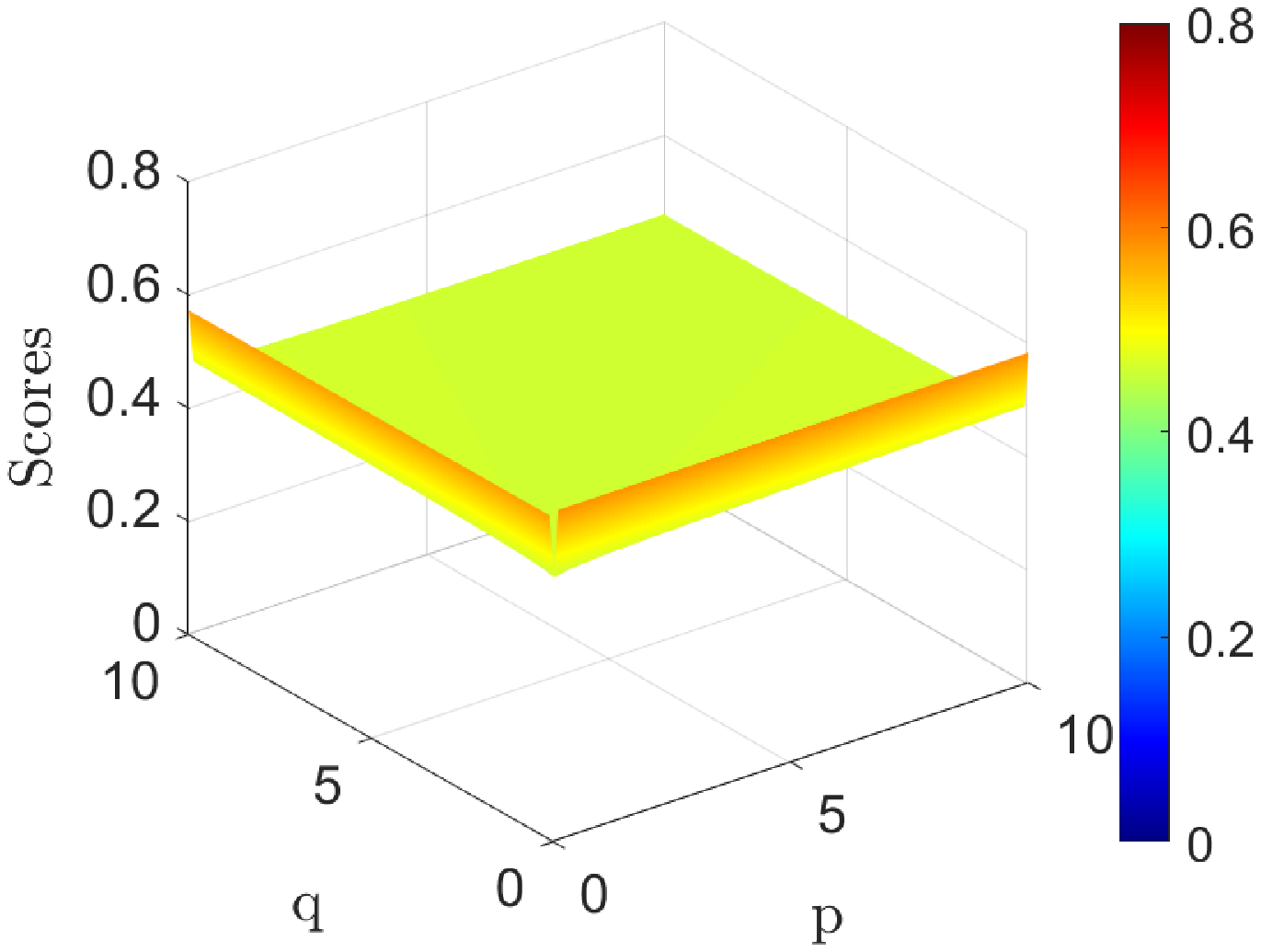}}}
\subfigure[Scores for $A_{3}$]
{\scalebox{0.35}{\includegraphics[]{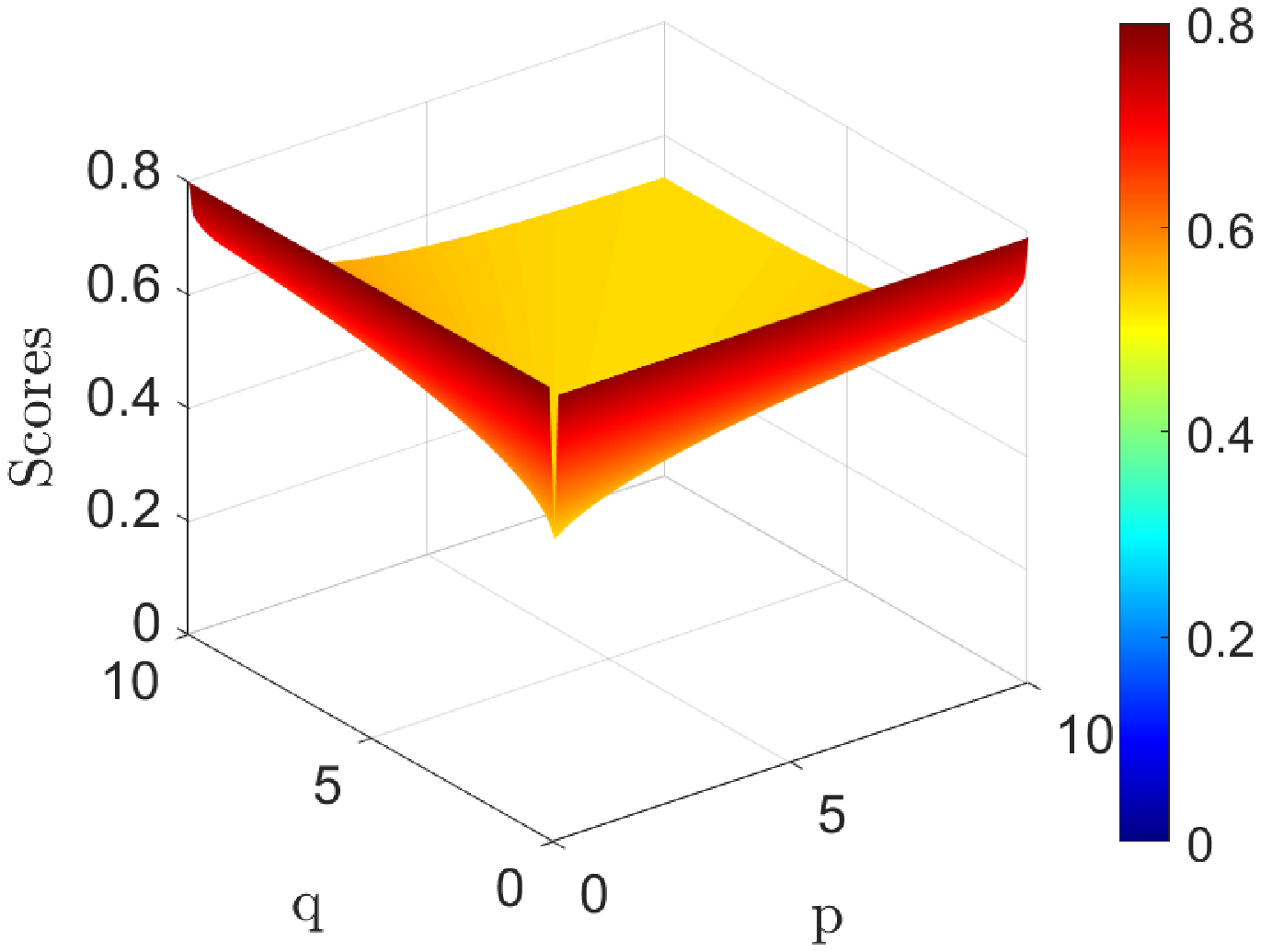}}}
\subfigure[Scores for $A_{4}$]
{\scalebox{0.35}{\includegraphics[]{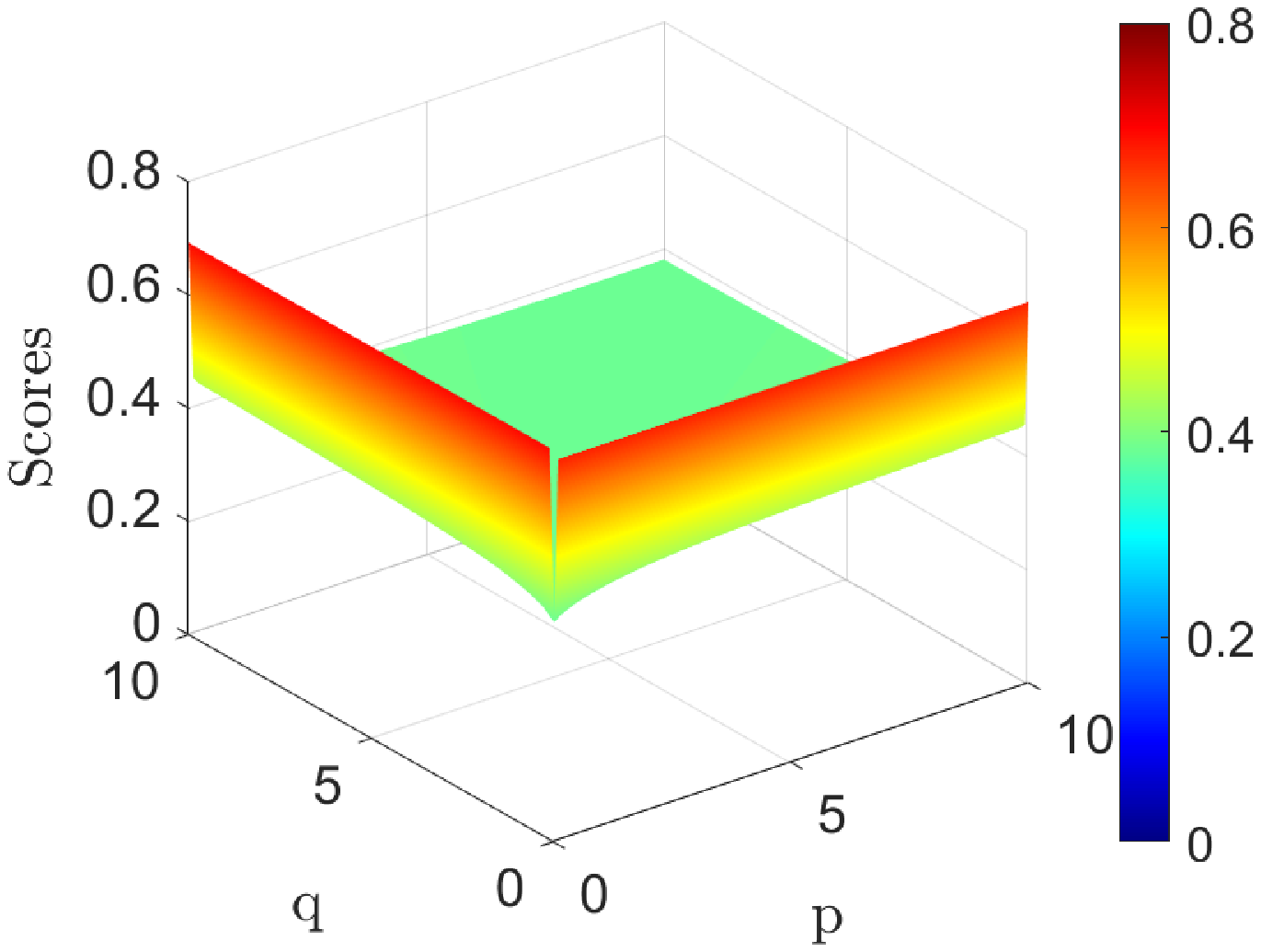}}}
\subfigure[Scores for $A_{5}$]
{\scalebox{0.35}{\includegraphics[]{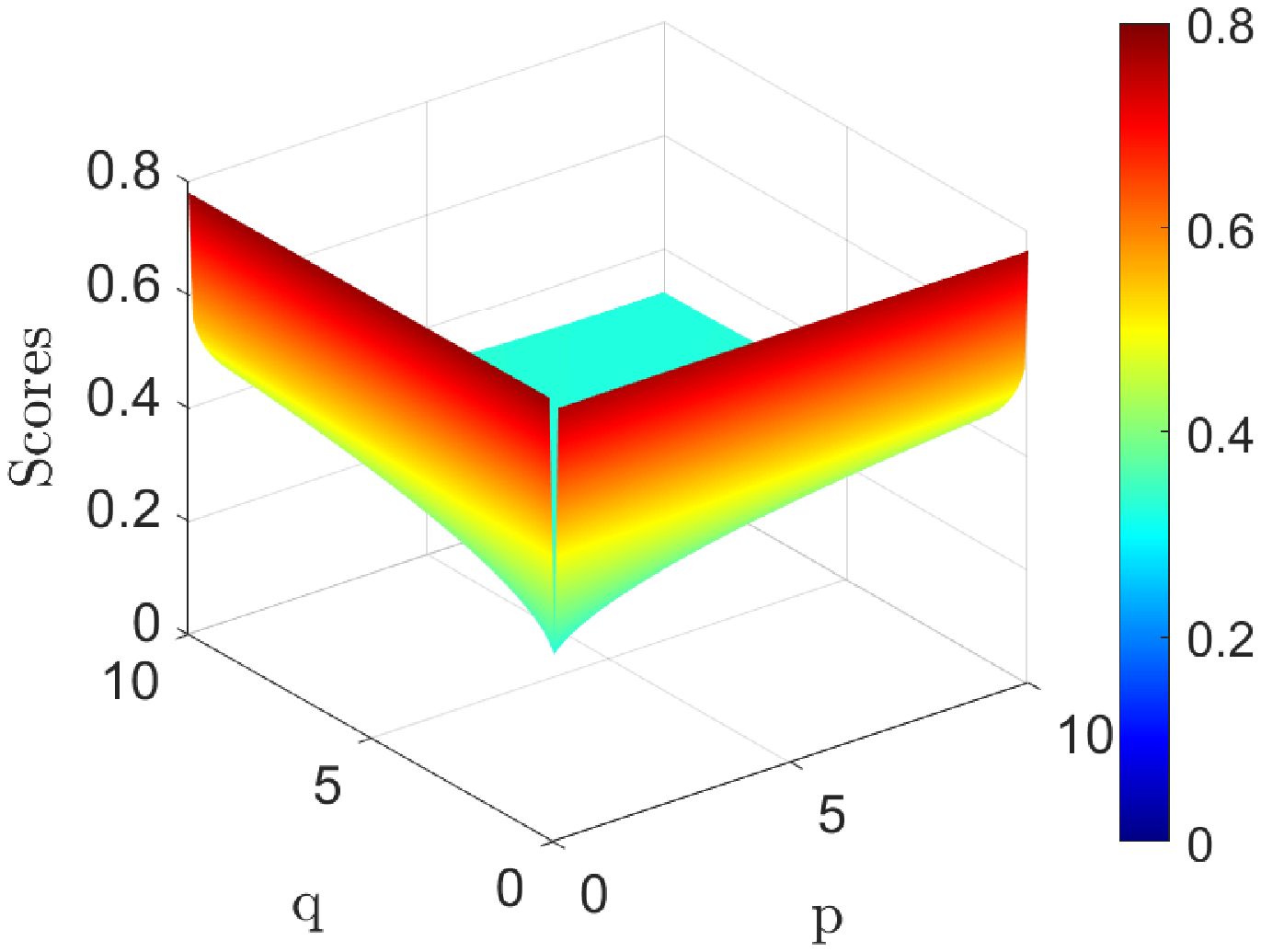}}}
\caption{Scores for alternatives in different values of $p$, $q$
obtained by $\mathrm{PFINWBM}_{T_{2}^{\mathbf{D}}}^{p,q}$}\label{fig-D}
\end{figure}

\begin{figure}[H]
\centering
\subfigure[Scores for $A_{1}$]
{\scalebox{0.35}{\includegraphics[]{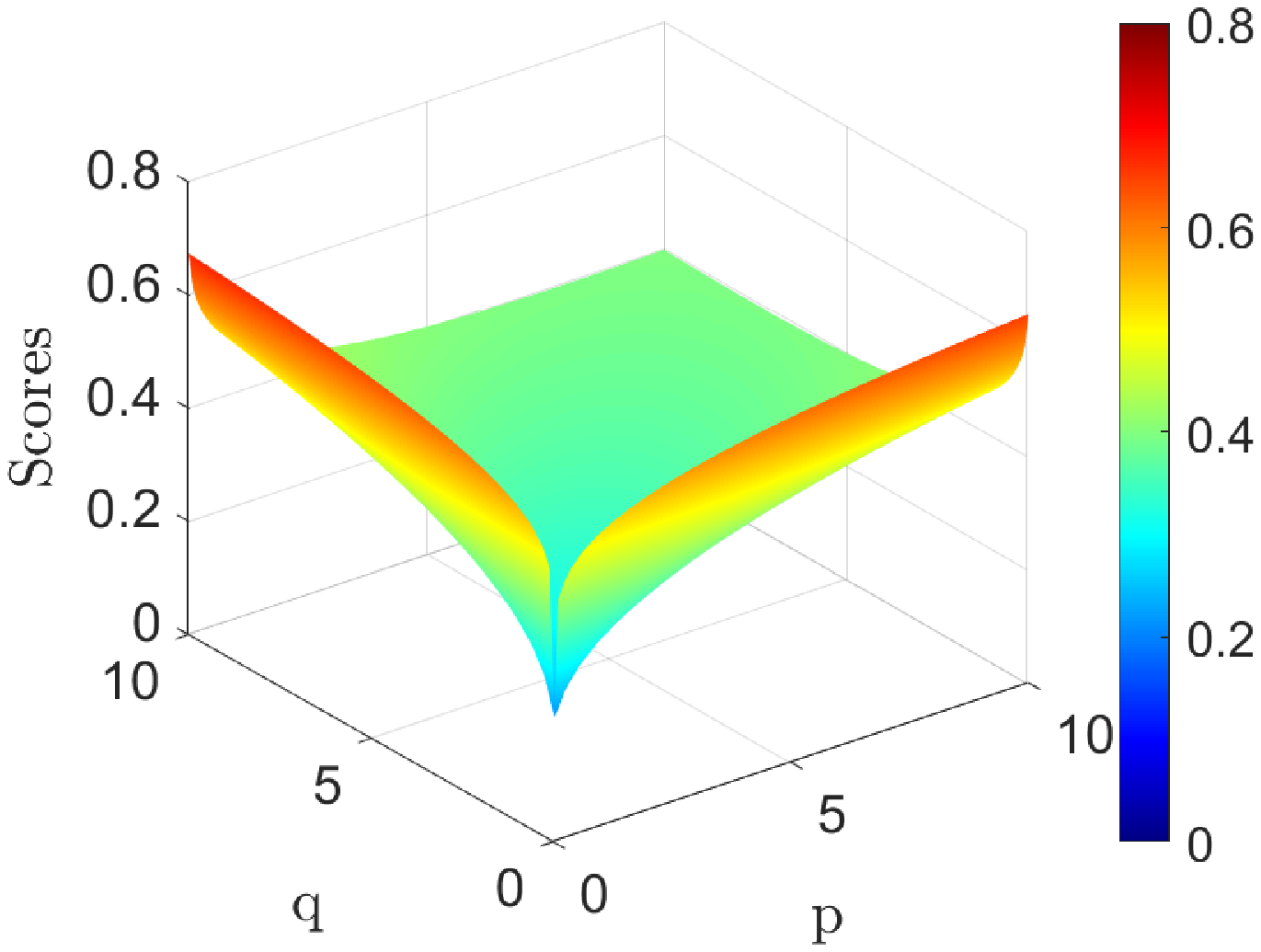}}}
\subfigure[Scores for $A_{2}$]
{\scalebox{0.35}{\includegraphics[]{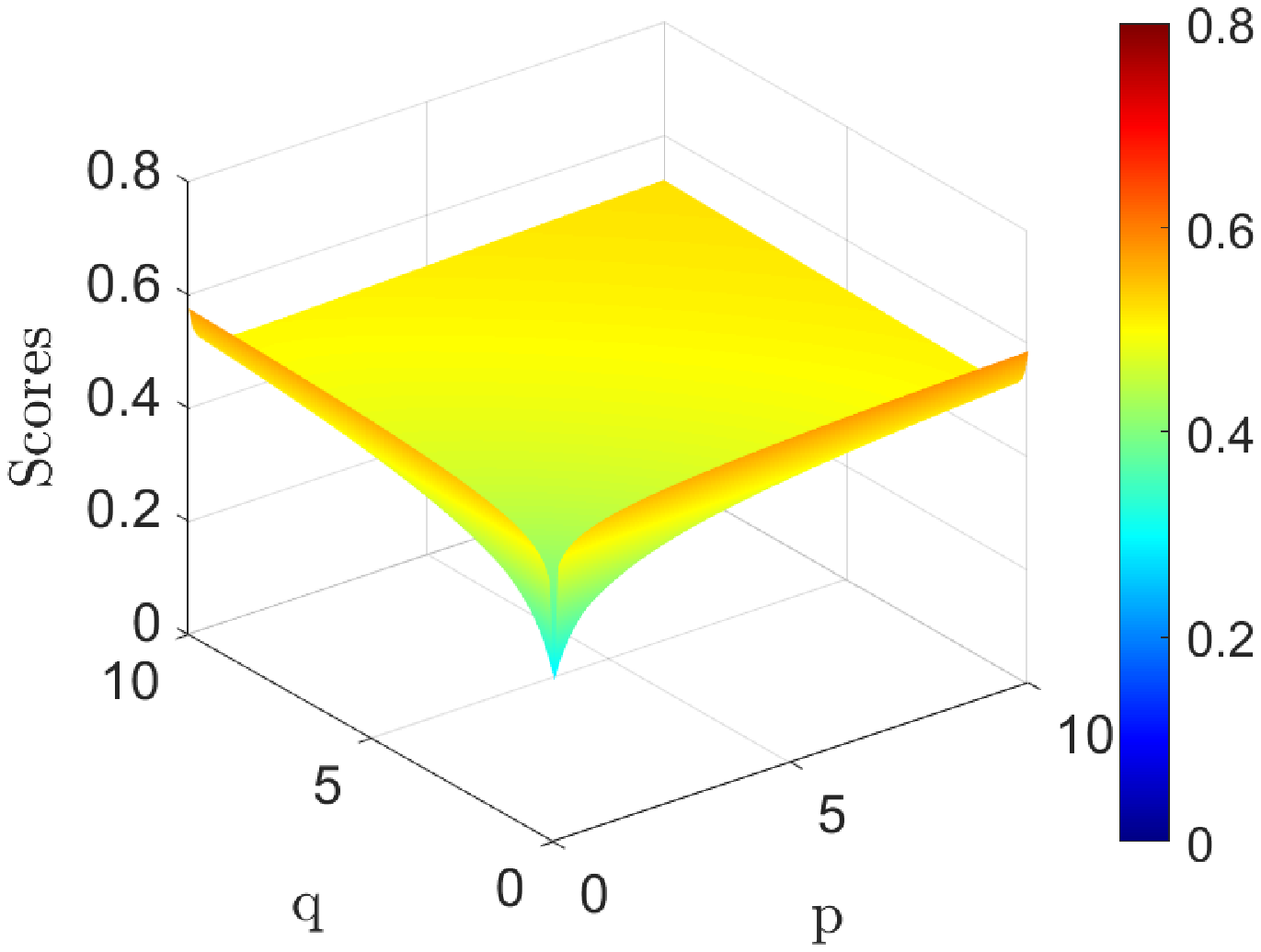}}}
\subfigure[Scores for $A_{3}$]
{\scalebox{0.35}{\includegraphics[]{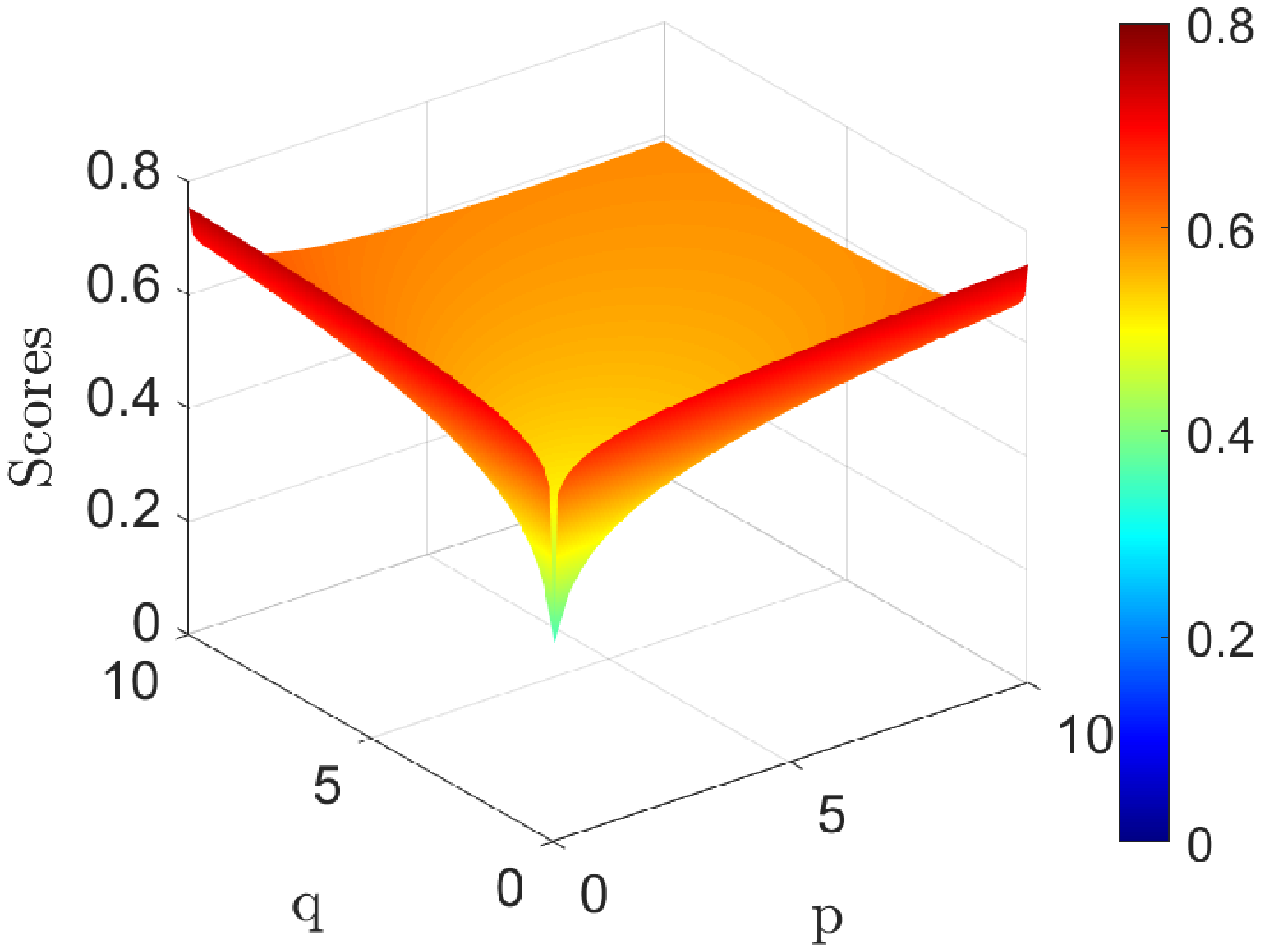}}}
\subfigure[Scores for $A_{4}$]
{\scalebox{0.35}{\includegraphics[]{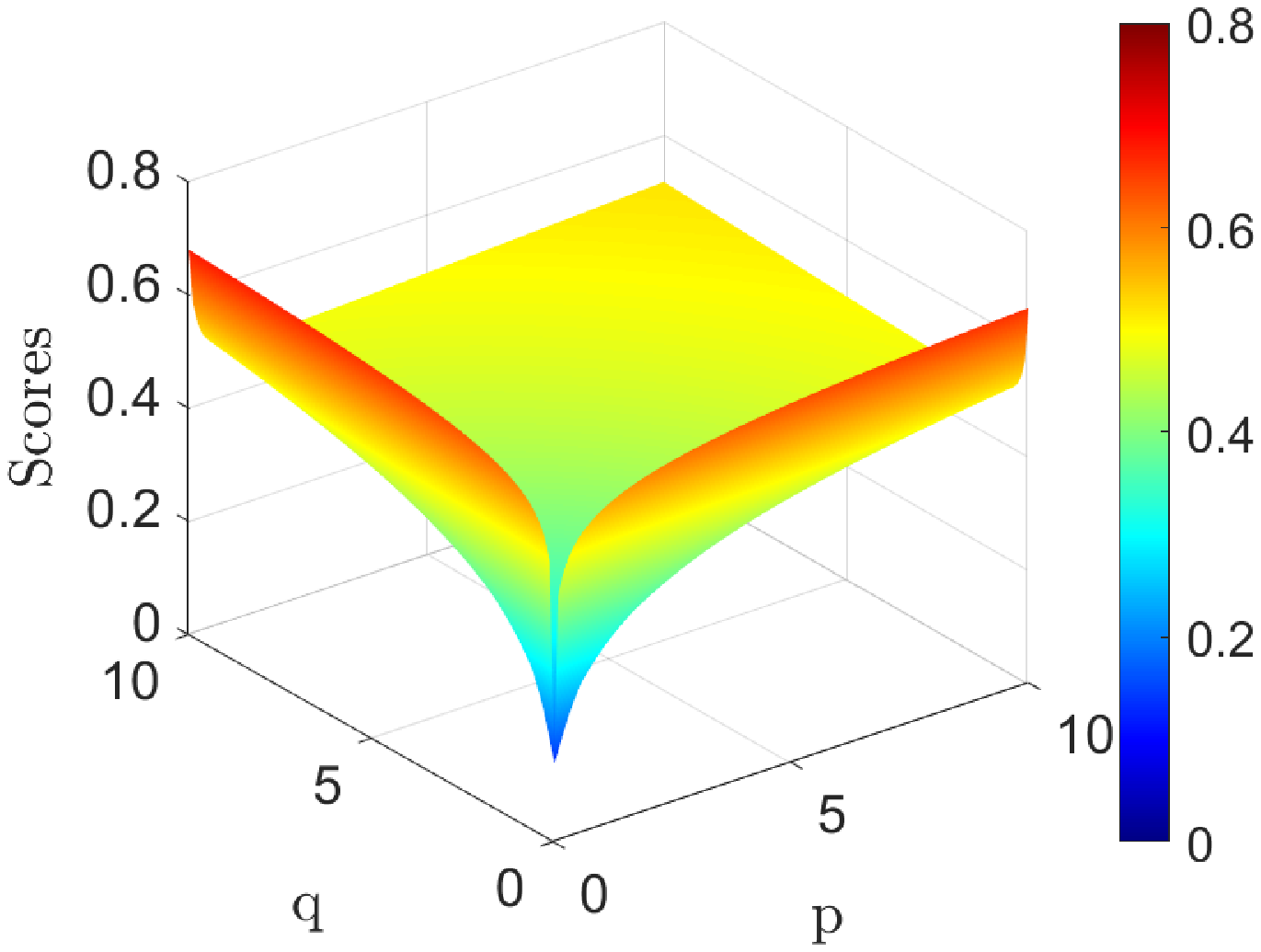}}}
\subfigure[Scores for $A_{5}$]
{\scalebox{0.35}{\includegraphics[]{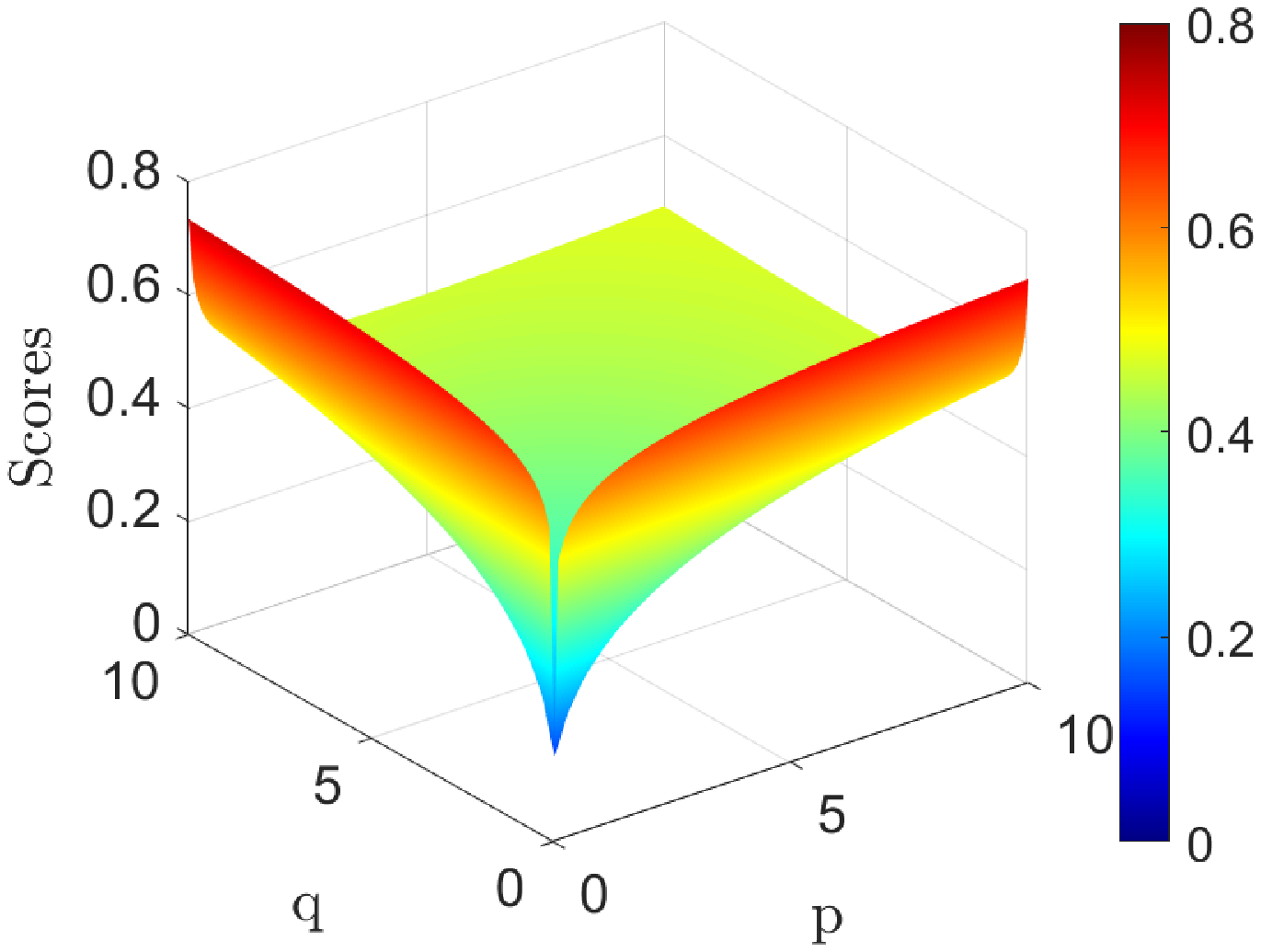}}}
\caption{Scores for alternatives in different values of $p$, $q$
obtained by $\mathrm{PFINWBM}_{T_{2}^{\mathbf{AA}}}^{p,q}$}\label{fig-AA}
\end{figure}

To illustrate the detailed influence of the parameters
$p$ and $q$ in the above example by using different PFINWBM operators induced by
 $T_{\gamma}^{\mathbf{H}}$, $T_{\gamma}^{\mathbf{SS}}$,
 $T_{\gamma}^{\mathbf{F}}$, $T_{\gamma}^{\mathbf{D}}$, and
 $T_{\gamma}^{\mathbf{AA}}$, the scores for alternatives $A_1$, $A_2$, $A_3$, $A_4$, $A_5$
 are shown in Figs.~\ref{fig-H}--\ref{fig-AA},
 respectively.

\subsection{The influence of the parameter $\gamma$ for MCDM results}

To study the changing trend of the scores and the rankings
of the alternatives $A_1, A_2, A_3, A_4, A_5$ with the change of
the t-norm $T$ and the parameter $\gamma$, we use the following
to illustrate these issues.

(1) Let $p=q=1$ and $\gamma=2$. If we use $\mathrm{PFINWBM}_{T_{\gamma}^{\mathbf{F}}}^{p,q}$,
$\mathrm{PFINWBM}_{T_{\gamma}^{\mathbf{D}}}^{p,q}$, $\mathrm{PFINWBM}_{T_{\gamma}^{\mathbf{AA}}}^{p,q}$
to aggregate the above PFNs in Table~\ref{tab1}, then the aggregated values, the scores, and the
ranking results are shown in Tables~\ref{tab5}--\ref{tab7}, respectively.
\begin{table}[H]
\centering
\caption{The aggregated values of the ERP systems}\label{tab5}
\resizebox{\columnwidth}{!}{
\begin{tabular}{lccc}
\toprule
~ &$\mathrm{PFINWBM}_{T_{\gamma}^{\mathbf{F}}}^{1,1}$ &$\mathrm{PFINWBM}_{T_{\gamma}^{\mathbf{D}}}^{1,1}$ &$\mathrm{PFINWBM}_{T_{\gamma}^{\mathbf{AA}}}^{1,1}$\\
\midrule
$r_{1}$ &$\langle0.3672, 0.5241, 0.0779\rangle$ &$\langle0.3639, 0.5385, 0.0606\rangle$ &$\langle0.3718, 0.5102, 0.0852\rangle$\\
$r_{2}$ &$\langle0.4369, 0.4114, 0.1080\rangle$ &$\langle0.5669, 0.2601, 0.1053\rangle$ &$\langle0.5160, 0.3238, 0.1097\rangle$\\
$r_{3}$ &$\langle0.4750, 0.3296, 0.0599\rangle$ &$\langle0.5871, 0.1807, 0.0599\rangle$ &$\langle0.5393, 0.2353, 0.0637\rangle$\\
$r_{4}$ &$\langle0.2806, 0.6396, 0.0587\rangle$ &$\langle0.4399, 0.4979, 0.0579\rangle$ &$\langle0.3640, 0.5723, 0.0591\rangle$\\
$r_{5}$ &$\langle0.3171, 0.5820, 0.0733\rangle$ &$\langle0.3939, 0.5094, 0.0699\rangle$ &$\langle0.3493, 0.5542, 0.0775\rangle$\\
\bottomrule
\end{tabular}
}
\end{table}

\begin{table}[H]
\centering 
\caption{The scores of the ERP systems}\label{tab6}
\resizebox{\columnwidth}{!}{
\begin{tabular}{lccc}
\toprule
 ~  & $\mathrm{PFINWBM}_{T_{\gamma}^{\mathbf{F}}}^{1,1}$
 &$\mathrm{PFINWBM}_{T_{\gamma}^{\mathbf{D}}}^{1,1}$ &$\mathrm{PFINWBM}_{T_{\gamma}^{\mathbf{AA}}}^{1,1}$\\
 \midrule
$S_{r_{1}}$ &0.2893 &0.3033 &0.2865\\
$S_{r_{2}}$ &0.3289 &0.4616 &0.4062\\
$S_{r_{3}}$ &0.4151 &0.5272 &0.4756\\
$S_{r_{4}}$ &0.2219 &0.3819 &0.3049\\
$S_{r_{5}}$ &0.2438 &0.3240 &0.2717\\
\bottomrule
\end{tabular}
}
\end{table}

\begin{table}[H]
\centering
\caption{The ranking results of the ERP systems}\label{tab7}
\resizebox{\columnwidth}{!}{
\begin{tabular}{p{5cm} cp{cm}}
\toprule
 ~  & Ranking \\
 \midrule
$\mathrm{PFINWBM}_{T_{\gamma}^{\mathbf{F}}}^{1,1}$ &$A_{3} \succ_{_{\mathrm{W}}} A_{2}
\succ_{_{\mathrm{W}}} A_{1} \succ_{_{\mathrm{W}}} A_{5} \succ_{_{\mathrm{W}}} A_{4}$\\
$\mathrm{PFINWBM}_{T_{\gamma}^{\mathbf{D}}}^{1,1}$ &$A_{3} \succ_{_{\mathrm{W}}} A_{2}
\succ_{_{\mathrm{W}}} A_{4} \succ_{_{\mathrm{W}}} A_{5} \succ_{_{\mathrm{W}}} A_{1}$\\
$\mathrm{PFINWBM}_{T_{\gamma}^{\mathbf{AA}}}^{1,1}$ &$A_{3} \succ_{_{\mathrm{W}}} A_{2}
\succ_{_{\mathrm{W}}} A_{4} \succ_{_{\mathrm{W}}} A_{1} \succ_{_{\mathrm{W}}} A_{5}$\\
\bottomrule
\end{tabular}
}
\end{table}

\begin{figure}[h]
\centering
\subfigure[Scores for alternatives in different values of $\gamma$ obtained by $\mathrm{PFINWBM}_{T_{\gamma}^{\mathbf{H}}}^{1,1}$]{\scalebox{0.294}{\includegraphics[]{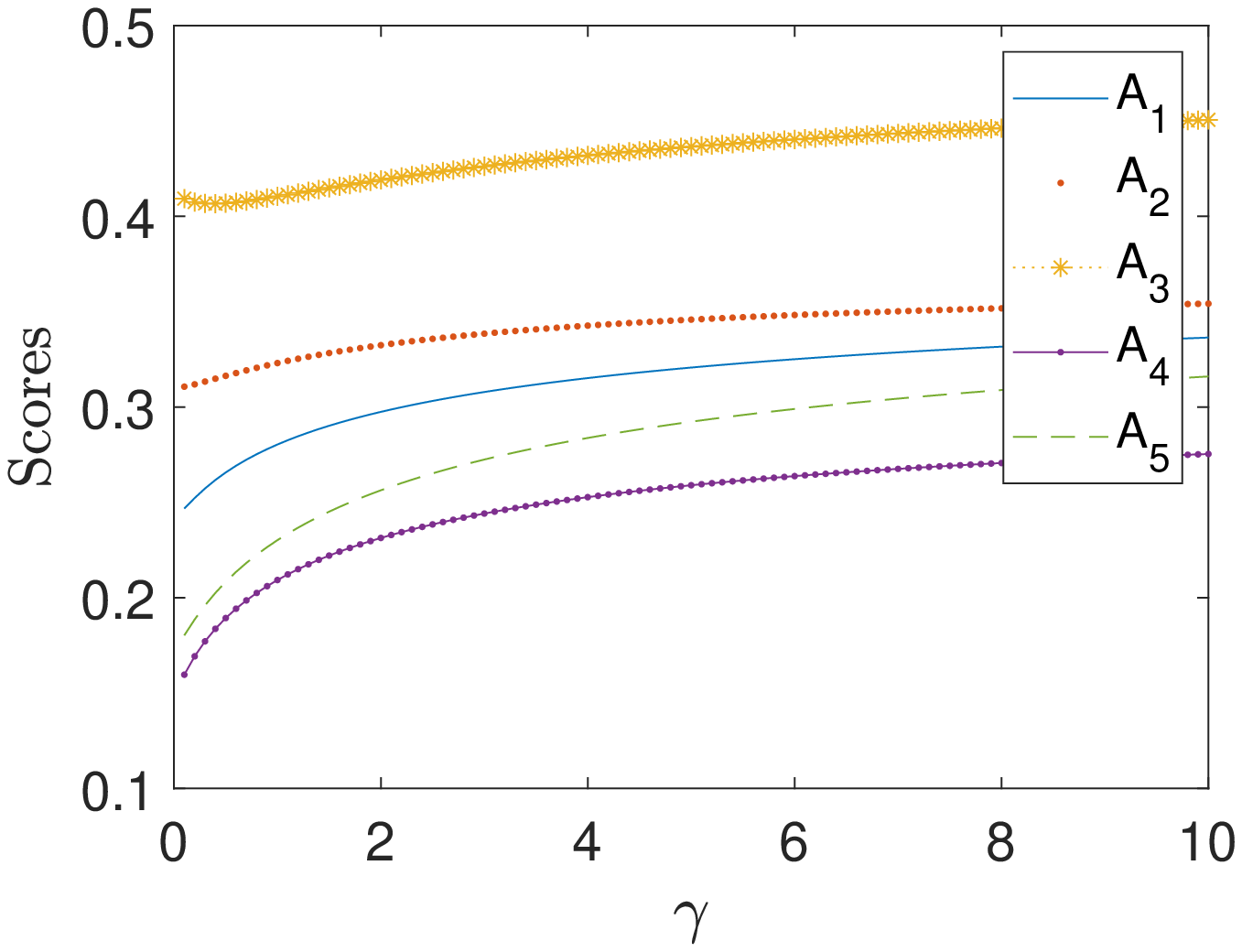}}}
\subfigure[Scores for alternatives in different values of $\gamma$ obtained by
$\mathrm{PFINWBM}_{T_{\gamma}^{\mathbf{SS}}}^{1,1}$]{\scalebox{0.294}{\includegraphics[]{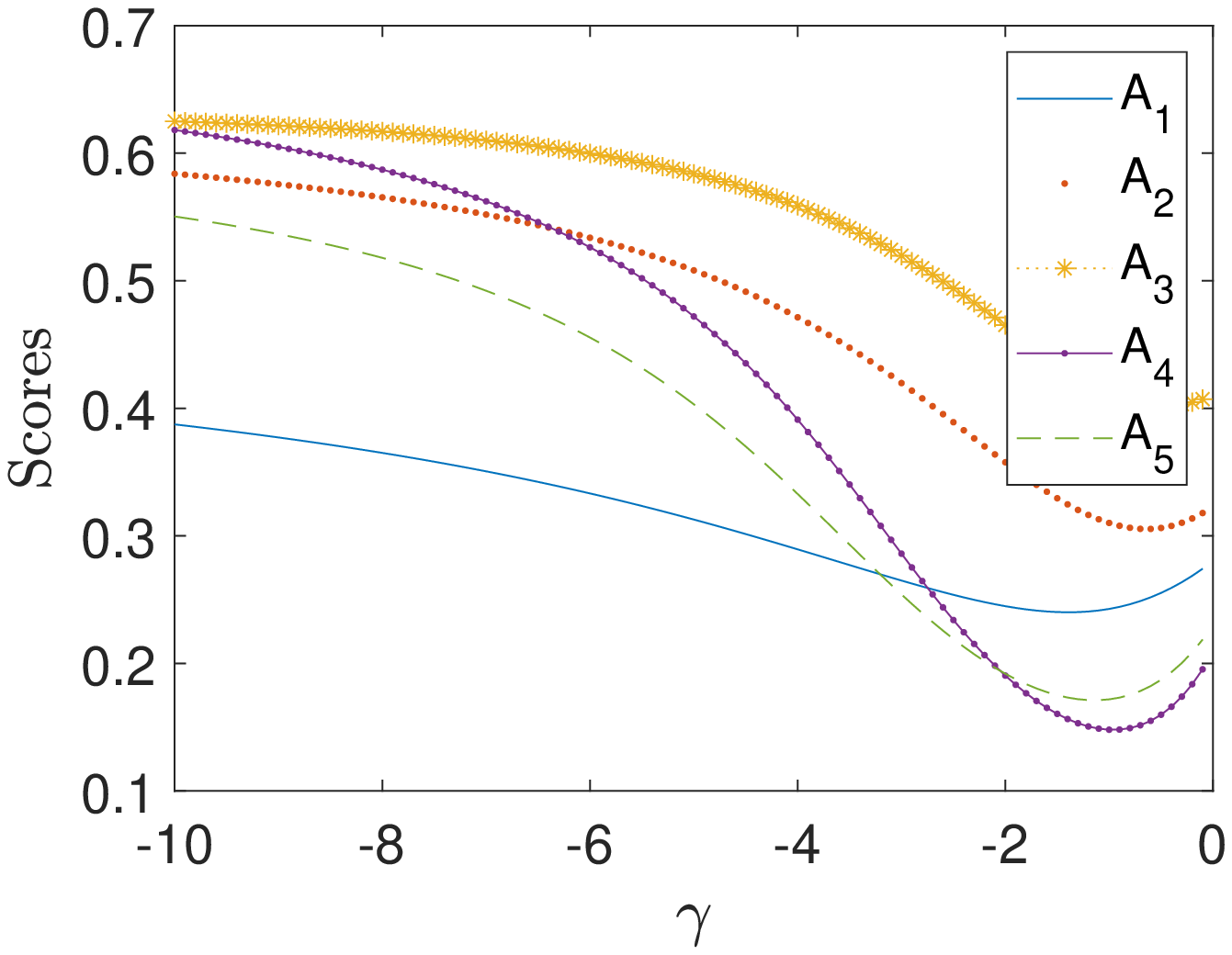}}}
\subfigure[Scores for alternatives in different values of $\gamma$ obtained by
$\mathrm{PFINWBM}_{T_{\gamma}^{\mathbf{F}}}^{1,1}$]{\scalebox{0.294}{\includegraphics[]{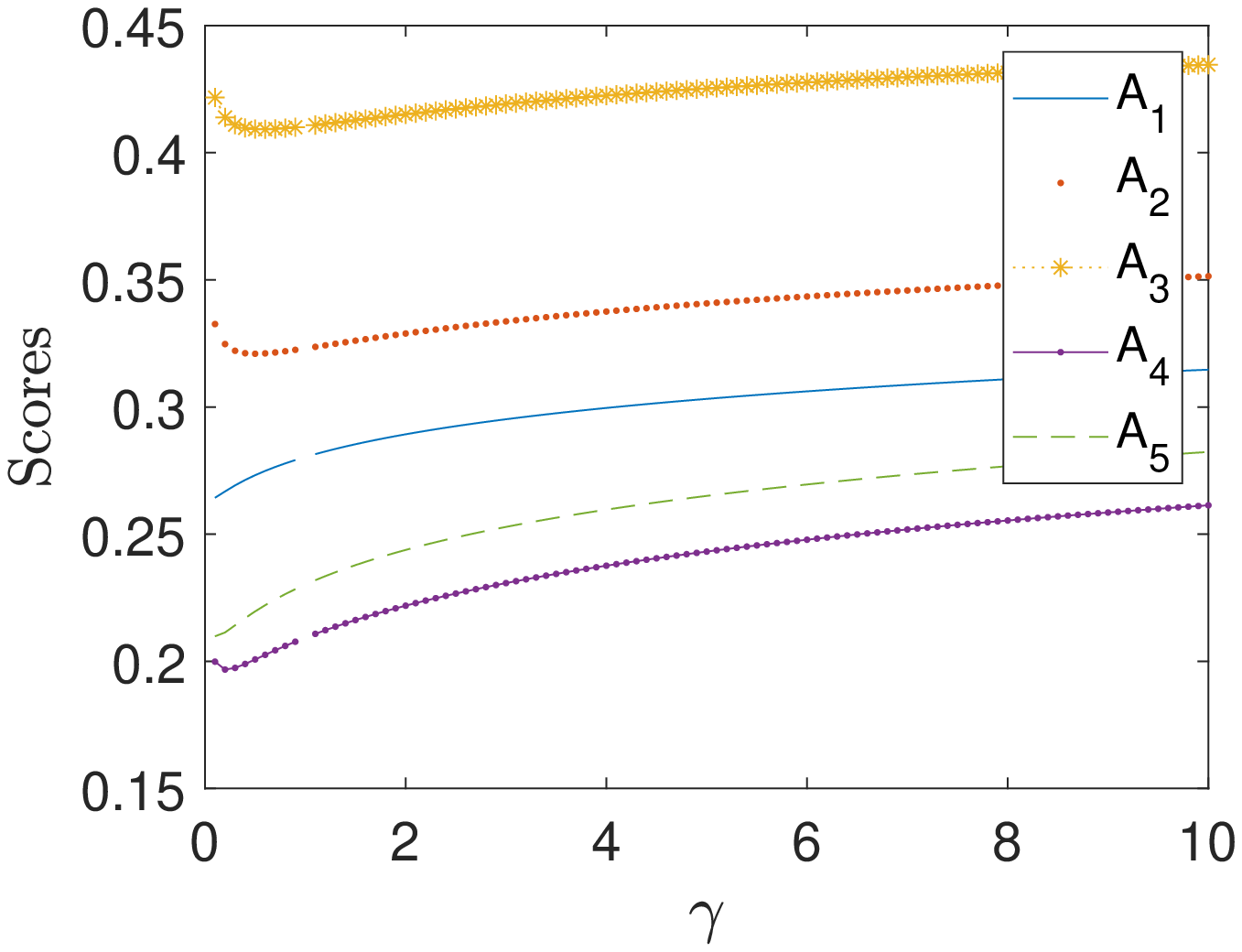}}}
\subfigure[Scores for alternatives in different values of $\gamma$ obtained by
$\mathrm{PFINWBM}_{T_{\gamma}^{\mathbf{D}}}^{1,1}$]{\scalebox{0.294}{\includegraphics[]{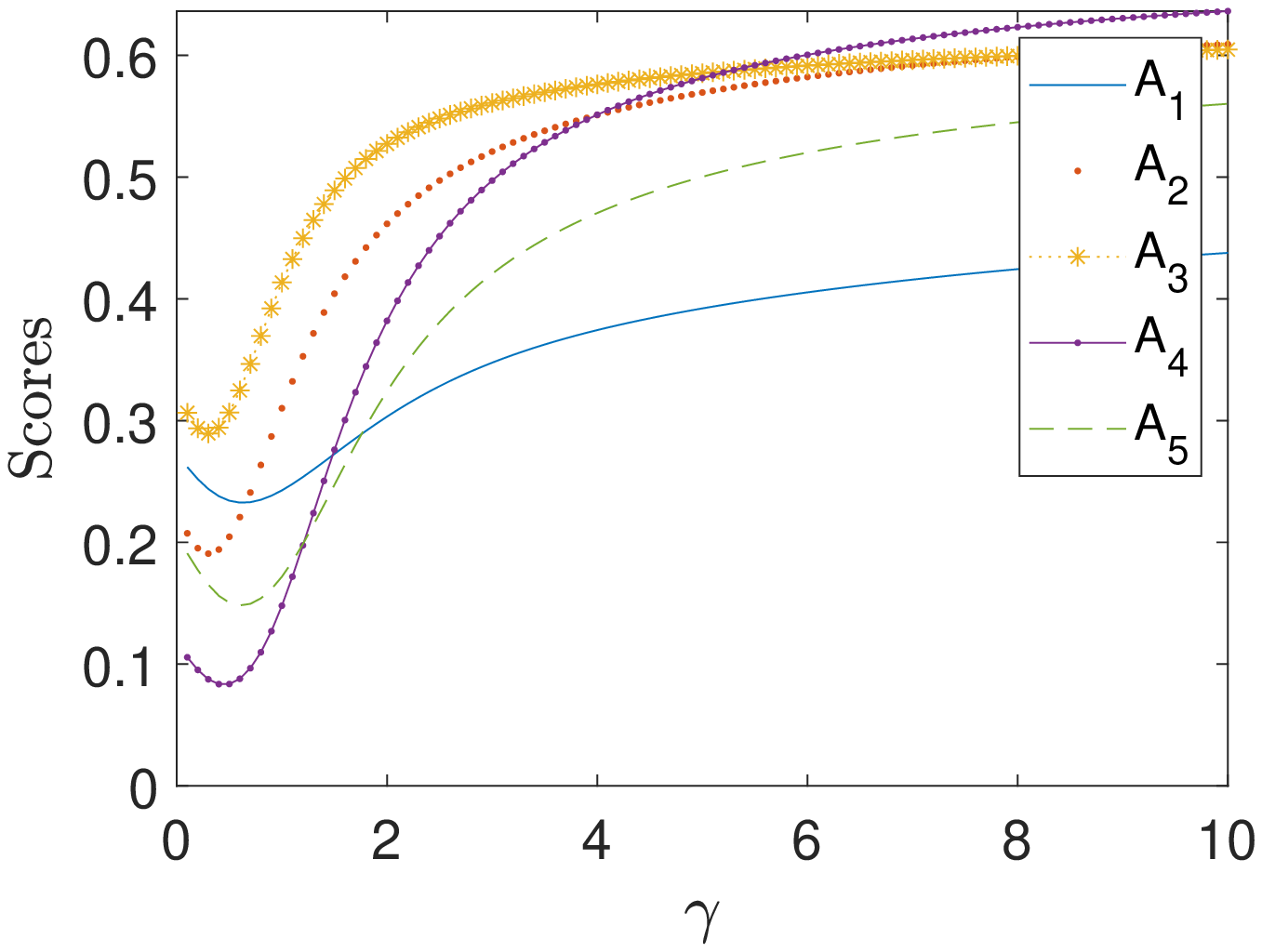}}}
\subfigure[Scores for alternatives in different values of $\gamma$ obtained by
$\mathrm{PFINWBM}_{T_{\gamma}^{\mathbf{AA}}}^{1,1}$]{\scalebox{0.294}{\includegraphics[]{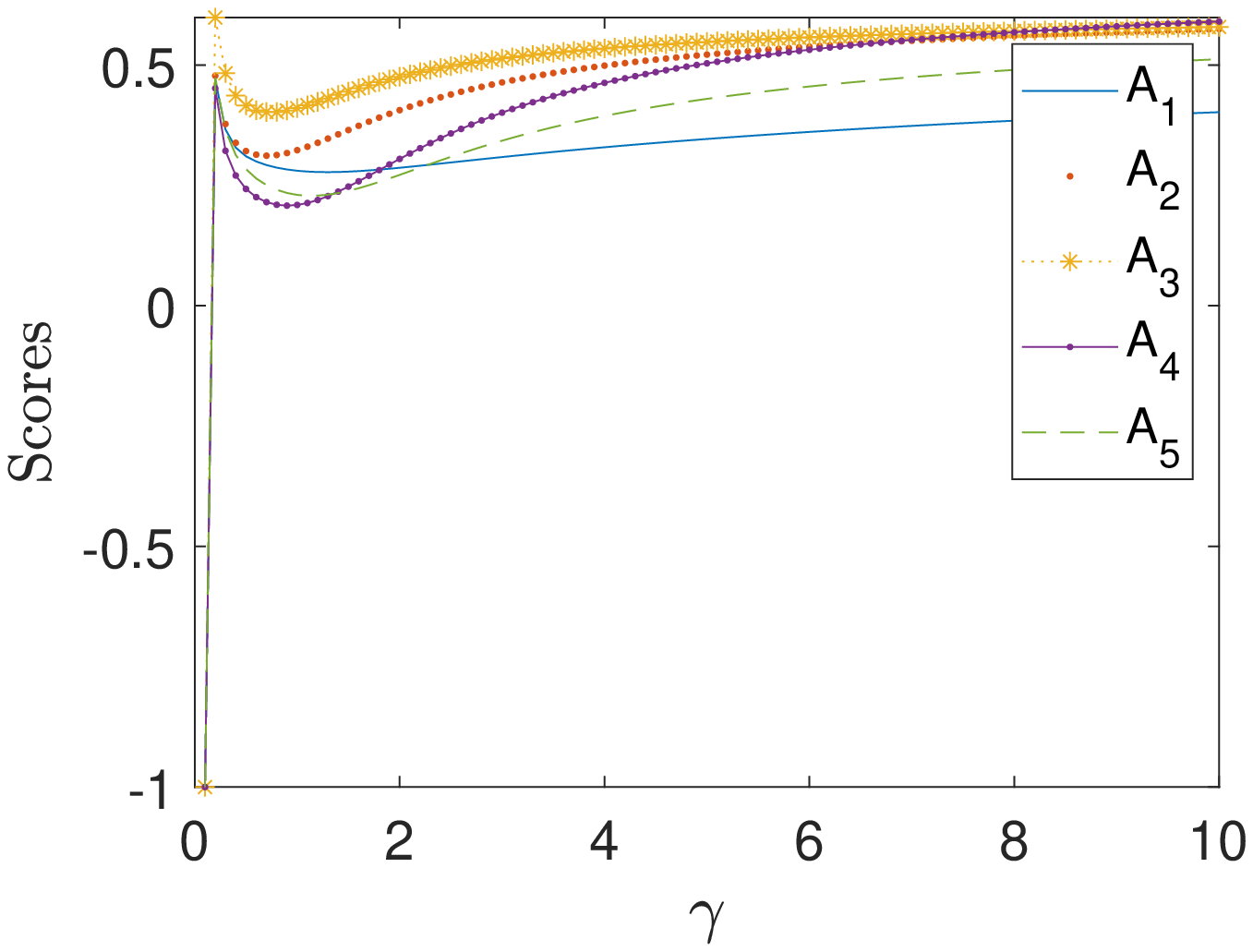}}}
\caption{Scores for alternatives in different values of $\gamma$ obtained by PFINWBM}\label{fig-PFINWBM}
\end{figure}

(2) Let $p=q=1$. If we change the values of the parameter $\gamma$ in PFINWBM operators induced by
 $T_{\gamma}^{\mathbf{H}}$, $T_{\gamma}^{\mathbf{SS}}$, $T_{\gamma}^{\mathbf{F}}$, $T_{\gamma}^{\mathbf{D}}$, and
 $T_{\gamma}^{\mathbf{AA}}$, the scores for alternatives $A_1, A_2, A_3, A_4, A_5$ are shown in
 Fig.~\ref{fig-PFINWBM} (a)--(e), respectively.

In general, from the above analysis, we observe that the parameter $\gamma$ can
be considered as a reflection of the decision makers' preferences,
as the parameter $\gamma$ changes in a certain range, although the scores of the
alternatives are different, and the rankings of the alternatives are
also different, the best ERP system is always $A_3$ except the aggregated results
obtained by $\mathrm{PFINWBM}_{T_{\gamma}^{\mathbf{D}}}^{1,1}$.
Therefore, we have every reason to conclude that the best ERP system is $A_3$.

\section{Conclusion}\label{Sec-8}
Duo to the special three-dimensional degree structure of PFNs, many
existing operators for PFNs are not closed in PFNs. However,
the closeness is very important for ensuring the fairness of
decision-making, which ensure that the aggregation output is
still a PFN; and so, the evaluation criteria are under a unified
framework. For this reason, Wu et al.~\cite{WZCLZY2021} introduced
four basic operations for PFNs, including addition, product,
scalar multiplication, and power, which are proved to be closed
in PFNs, monotonous, idempotent, bounded, shift-invariant, and homogeneous.
Based on these four basic operations for PFNs, we introduce the
picture fuzzy interactional Bonferroni mean (PFIBM), picture fuzzy interactional
weighted Bonferroni mean (PFIWBM), and picture fuzzy interactional normalized weighted
Bonferroni mean (PFINWBM) operators for PFNs. Furthermore, we prove that
PFIBM, PFIWBM, and PFINWBM operators are monotonous under the linear order
$\preccurlyeq_{_{\mathrm{W}}}$ in \cite{WZCLZY2021}, idempotent, bounded, and
commutative. To this end, we propose a novel MCDM method under the picture fuzzy
environment by using PFINWBM operator, which is applied to the enterprise
resource planning systems selection. By using six classes of well-known
triangular norms, including the algebraic product $T_{\mathbf{P}}$, Schweizer-Sklar
t-norm $T_{\gamma }^{\mathbf{SS}}$, Hamacher t-norm $T_{\gamma }^{\mathbf{H}}$,
Frank t-norm $T_{\gamma }^{\mathbf{F}}$, Dombi t-norm $T_{\gamma }^{\mathbf{D}}$,
and Acz\'{e}l-Alsina t-norm $T_{\gamma }^{\mathbf{AA}}$, the best ERP system
is always the same one, demonstrating the effectiveness of our method.

\bibliographystyle{IEEEtran}
\bibliography{IEEEexample}


\end{document}